\documentclass{article}

\usepackage{arxiv}
\usepackage{amsmath}
\usepackage{amsthm}
 
\usepackage[T1]{fontenc}    % use 8-bit T1 fonts
\usepackage{hyperref}       % hyperlinks
\usepackage{url}            % simple URL typesetting

\usepackage{booktabs}       % professional-quality tables
\usepackage{amsfonts}       % blackboard math symbols
\usepackage{nicefrac}       % compact symbols for 1/2, etc.
\usepackage{bm}             % bold math symbols
\usepackage{xcolor}         % color support for \textcolor
\usepackage{microtype}      % microtypography
\usepackage{lipsum}		% Can be removed after putting your text content
\usepackage{graphicx}
\usepackage{natbib}
\usepackage{doi}

\usepackage{graphicx} 
\newtheorem{theorem}{Theorem}[section]
\newtheorem{lemma}[theorem]{Lemma}

\newtheorem{proposition}[theorem]{Proposition}

\newtheorem{remark}{Remark}

\title{Parameter Estimation for Fractional Autoregressive Process with Periodic Structure}

\author{
 Cai Chunhao \\
  School of mathematics(Zhuhai) \\
  Sun Yat-sen University \\
  Zhuhai,People’s Republic of China\\
  \texttt{caichh9@mail.sysu.edu.cn} \\
   \And
  Shang Yiwu \\
  School of mathematical sciences \\
  Nankai University \\
  Tianjin,People’s Republic of China\\
  \texttt{shangyiwu@mail.nankai.edu.cn} \\
}
\renewcommand{\shorttitle}{}

\hypersetup{
pdftitle={Parameter Estimation for Fractional Autoregressive Process with Periodic Structure},
pdfsubject={q-bio.NC, q-bio.QM},
pdfauthor={David S.~Hippocampus, Elias D.~Striatum},
pdfkeywords={First keyword, Second keyword, More},
}

\begin{document}
\maketitle
\begin{abstract}
This paper introduces a new periodic fractional autoregressive process (PFAR) driven by fractional Gaussian noise (fGn) to model time series of precipitation evapotranspiration. Compared with the similar model in [\emph{Water Resources Research}, \textbf{20} (1984) 1898--1908], the new model incorporates a periodic structure via specialized varying coefficients and captures long memory and rough voltality through fGn for $0<H<1$, rather than via fractional differencing. In this work, Generalized Least Squares  Estimation (GLSE) and the GPH method are employed to construct an initial estimator for the joint estimation of model parameters. A One-Step procedure is then used to obtain a more asymptotically efficient estimator. The paper proves that both estimators are consistent and asymptotically normal, and their performance is demonstrated via Monte Carlo simulations with finite-size samples. Simulation studies suggest that, while both estimation methods can accurately estimate the model parameters, the One-Step estimator outperforms the initial estimator.
\end{abstract}

\section{Introduction}
In many fields of time series analysis, such as finance, meteorology, and engineering, the phenomena of long memory and periodicity have received extensive attention. Long memory implies that there is still a strong correlation in time series data over a relatively long time interval. For example, the precipitation during the dry season will remain at a relatively low level, while the precipitation during the rainy season will remain at a relatively high level; the precipitation of each month or quarter is not independent. Periodicity indicates that the data will exhibit similar patterns repeatedly at fixed time intervals. These two characteristics are of paramount significance in deciphering the evolution of time series and can provide a good explanation for some phenomena in nature $\cite{hosking1984modeling}$. However, they simultaneously present two significant challenges to the study of classical stationary time series. The first challenge lies in nonstationarity, while the second stems from long-range dependence.

The main contributions of this paper are as follows: first, an estimator for the long memory parameter that is independent of periodicity is derived when $0<H<1$; second, the estimation of the periodic time-varying coefficients in the PFAR model is accomplished via a modified  GLSE; third, the spectral density of this model is computed, and the aforementioned initial estimators are refined using a One-Step estimation method. Concurrently with the derivation of these estimators, the asymptotic properties of the proposed estimators are also proven. Through the work presented above, we have addressed the challenges in estimation posed by periodicity and long memory in the PFAR model.

The time series involved in  climatology \citep{hannan1955test} , hydrology \citep{fernandez1986periodic} and economics \citep{franses2004periodic} often exhibit periodic characteristics. Unlike other varying coefficients time series,  the coefficients of periodic time series satisfy 
$$\phi(u)=\phi(u+nT), n \in\mathbb{N},u=1,2,\dots T$$
where $T$ is the period. Here, we consider p-order  PFAR model which has the recurrence
 \begin{equation}
	X_{nT+u} =\sum_{i=1}^{p}\phi_{i}(nT+u)X_{nT+u-i}+\epsilon^{H}_{nT+u}, \quad u= 1,2,...T,  \label{origneq}  
 \end{equation}
 where $\epsilon^{H}_{nT+u}$ is fractional Gaussian noise, which explains the nonperiodic fluctuations. 
	 Fractional Gaussian noise exhibits long memory when $\frac{1}{2}<H<1 $. The long memory phenomenon indicates strong autocorrelation or dependence in time series data. Unlike general fractional difference models, this model can not only exhibit long memory but also be regarded as the discretization of a rough volatility model and can be used to study some low-frequency data in financial markets. We typically say that $X_{t}$ has long memory if its covariance satisfies 
\begin{equation}
 \gamma_{j} \sim {Q}j^{2H-2},\quad j \rightarrow \infty, \label{1.2}
 \end{equation}
	the spectral density is defined by the scheme 
\begin{equation}
 f(\lambda) \sim  {V}\lambda^{1-2H},\quad \lambda \rightarrow 0^{+}, \label{1.3}
\end{equation}
	where $\frac{1}{2} < H < 1$,  {Q} and {V} are constants greater than 0. This paper not only accounts for the case where the $\frac{1}{2}<H<1$, but also performs estimation for the case when $0<H<\frac{1}{2}$, which characterizes anti-persistence and rough volatility.\\
     Traditional volatility models fail to capture volatility's non-smoothness, while rough volatility models rooted in fractional Brownian motion with Hurst exponent \( H <\frac{1}{2} \) address this gap: \citep{comte1998long} laid groundwork with the Long Memory Stochastic Volatility (LMSV) model, \citep{gloter2004stochastic} advanced theoretical understanding of fractional Brownian motion and rough volatility,\cite{corsi2021roughness} empirically confirmed volatility's roughness via S\&P 500 option data. For discrete-time applications, \citep{richard2023discrete} studies the Euler-type discrete-time schemes for the rough Heston model, proves the convergence of these schemes, and evaluates the pricing of different options under the model through numerical examples. Our PFAR model can also be used to study the rough volatility of financial markets when $0<H<\frac{1}{2}$.\\
The stationarity of a time series is a fundamental assumption for numerous statistical analyses and modeling methods. It requires that the statistical properties of the series, such as the mean, variance, and autocovariance, do not change over time. However, when a time series exhibits a periodic structure, its inherent stationarity is usually disrupted. A periodic structure implies that the series repeats similar patterns at specific time intervals, and this regularity leads to the time-varying nature of the series' statistical properties. In this paper, we aim to tackle this challenge by extracting stationarity from the PFAR model and thoroughly analyzing the difficulties its intricate structure poses to the initial estimator. \\ 
The discourse regarding the parameter estimation of the autoregressive (PAR) model predominantly centers around moment estimation, least squares estimation(\citep{jones1967time}), and maximum likelihood estimation(\citep{vecchia1985maximum}). It is a widely recognized fact that the shortcoming inherent in these methodologies is that their respective estimated values exhibit a high degree of sensitivity to outliers and minor fluctuations in the distribution pattern. 
 
 Previous studies by \citep{brouste2014asymptotic} and \citep{soltane2024asymptotic} have laid a foundation for the estimation of the parameters $\phi_{i}(u)$ in FAR models. In this paper, we use the modified GLSE  proposed by \citep{esstafa2019long} and \citep{hariz2024fast} to obtain a consistent estimator of $\phi_{i}(u)$. However, the direct use of GLSE proposed in the aforementioned references leads to biased estimation (see Remark 7 in this paper). Therefore, by taking into account the relationship between the scaling matrix in the GLSE and the periodic structure of the data, we propose an new alternative GLSE. We will prove that this estimator is asymptotically normal and consistent.

Time series models with long memory show long range dependencies between distant observations, posing challenges to traditional statistical analysis and forecasting. In the PFAR model, long memory comes from fractional Gaussian noise, where the parameter $H$ determines this characteristic. Thus, estimating $H$ is crucial. The first method for estimating $H$ was the rescaled range analysis by \citep{hurst1951long}, but its lack of a limiting distribution complicates statistical inference. Now, popular estimation techniques are the GPH estimation by \citep{geweke1983estimation} and the local Whittle estimation by \citep{robinson1995log}.   

For the estimation of the Hurst index $H$, we will adopt the Geweke Porter-Hudak (GPH) method, which exhibits a smaller bias, for an additive stationary time series derived from the samples. It is important to note that estimating $\hat{H}_{n}$ using the sequence $(X_{nT + u})_{n\in\mathbb{N}}$ may seem more straightforward. Nevertheless, this approach does not differ substantially from the method proposed in \citep{hariz2024fast}. Moreover, each $\hat{H}_{n}(u)$ fails to encapsulate information from all the available data. Given that the sequence $(X_{nT + u})_{n\in\mathbb{N},u = 1,\dots,T}$ consists of data of the same underlying nature, we make the assumption that the long memory parameter remains constant across all periods and is independent of the specific period index $u$. The disparities between different periods are solely governed by the periodic parameters. \\
To derive a unique estimate of $\hat{H}_{n}$, we aggregate the data within each cycle to form a new sequence $(Z_{n})_{n\in\mathbb{N}}$. Subsequently, we establish the stationarity of $(Z_{n})_{n\in\mathbb{N}}$, compute its spectral density, and ultimately employ the GPH method to obtain the desired estimate $\hat{H}_{n}$. 
 This improvement enables us to address the issue of parameter estimation for \(H\) in nonstationary time series with periodicity.  

After obtaining the initial estimators of $\phi_{i}(u)$ and  $H$, we modify our approach using a faster and asymptotically efficient method known as the One-Step estimator. This method, first proposed by \citep{le1956asymptotic}, has been widely applied in ergodic Markov chains \citep{kutoyants2016multi}, diffusion processes \citep{gloter2021adaptive}, and fractional autoregressive processes \citep{hariz2024fast}. In order to obtain the form of the One-Step estimation, we specifically calculated the spectral densities of sequences \((X_{nT + u})_{n\in\mathbb{N}}\) and $(Z_{n})_{n \in \mathbb{N}}$.
 
 This paper is organized as follows. Sections 2 introduces some notations and assumptions. Section 3 introduces the initial estimator of the Hurst index, $\phi(1),\phi(2)$. Section 4 discusses the One-Step estimator and its asymptotic properties. Section 5 provides numerical illustrations to demonstrate the performance of both the initial and One-Step estimators. Section 6 demonstrates the performance of the model in real-world scenarios. Section 7 provides a summary of this article. All the proofs in this paper are presented in the Appendix.

\section{Preliminaries}
\subsection{Periodic Fractional Autoregressive Process, Notations and Assumptions}
  \par From the PFAR model representation in (\ref{origneq}), we can consider the first-order model without loss of generality in this paper, $X_{nT+u}$ is said to be a PFAR(1) model if it admits the representation
  \begin{equation}
  	X_{nT+u}=\phi(nT+u)X_{nT+u-1}+\epsilon^{H}_{nT+u}, \quad u=1,2,3...T, \quad n \in \mathbb{N},   \label{eq:FormofPFAR(1)}
  \end{equation}
   where $\phi(u) = \phi(u+nT)$, $T$ represents the period length and $u$ denotes the $u$-th period of the $n$-th cycle. $\epsilon^{H}_{nT+u}=B^{H}_{nT+u+1}-B^{H}_{nT+u}$ is a stationary fractional Gaussian noise of hurst index $H$, $H \in (0,1)$. The autocovariance of sequence  $(\epsilon^{H}_{n})_{n\in \mathbb{N}}$ takes the form of 
  \begin{equation}
  	\rho_{\epsilon^{H}}(k) = \frac{1}{2}(|k+1|^{2H}-2|k|^{2H}+|k-1|^{2H}), \label{2.2}
  \end{equation}
  The spectral density of  $(\epsilon^{H}_{n})_{n\in \mathbb{N}}$ defined by
  \begin{equation}
  	{f_{\epsilon_{n}^{H}}(\lambda)} = C_{H}(1-cos(\lambda))\sum\limits_{j \in Z}\frac{1}{|\lambda+2j\pi|^{2H+1}}, \label{2.3}
  \end{equation}
where $C_H=\frac{1}{2\pi}\Gamma(2H+1)sin(\pi H)$ and $\lambda \in [-\pi,\pi]$, $\Gamma(\cdot)$ is Gamma function, 
  \par Here are some assumptions and notations bellow.
  \par {${\mathcal{A}_{0}}$: Denote \(\Theta_{u}^{l^\star}\) as a compact set with the following expression,
  \begin{center}
  	${\Theta_{u}^{l}}^{\star} = \left\{  \phi(u) \in R; \right. $ the roots of $  1-\phi(u)z=0$ have modulus$\left. \geq 1+l \right\} $
  \end{center}
  \par {We define the set  $\Theta^{l}_{u}$ as the Cartesian product  $\Theta^{l^\star}_{u} \times [d_{1},d_{2}]$, where $l$ is a positive constant and $[d_{1},d_{2}] \in (0,1)$.}
  \par {${\mathcal{A}_{1}}$: \(\phi(u) \in (-1, 1)\) and \(H \in (0, 1)\). }
 \par  {\textbf{Notation:} By $\xrightarrow{\mathcal{L}}$ and $\xrightarrow{\mathbb{P}}$, respectively, we denote convergence in law and convergence in probability. $\xrightarrow{L^{2}}$ means $L^{2}$ convergence when $n \rightarrow \infty$. Let $\boldsymbol{\phi}=(\phi(1),\phi(2),\ldots,\phi(T))$. Denote the parameters $\boldsymbol\theta=(\phi(1),\phi(2),\dots,\phi(T),H)$, $\boldsymbol\theta(u)=(\phi(u),H)$, where $\boldsymbol\theta(u)\in\mathring{\Theta}^{l}_{u}$, and $\mathring{\Theta}^{l}_{u}$ represents the interior of $\Theta^{l}_{u}$.}
 \par  {Define the parameter space $\Theta^{l}=\Theta^{l^{\star}}_{1}\times\Theta^{l^{\star}}_{2}\times\cdots\times\Theta^{l^{\star}}_{T}\times[d_{1},d_{2}]$, which encompasses all the required parameters. Given samples of size $n$, we obtain the estimators $\hat{\boldsymbol\theta}_{n}=(\hat{\phi}_{n}(1),\hat{\phi}_{n}(2),\ldots,\hat{\phi}_{n}(T),\hat{H}_{n})$ and $\hat{\boldsymbol\theta}_{n}(u)=(\hat{\phi}_{n}(u),\hat{H}_{n})$. 
}

  \par In this article, we will assume, without loss of generality, that T=2. From equation (\ref{eq:FormofPFAR(1)}), one can write
\begin{equation}
  \left(
  \begin{array}{c}
  {X_{2}}\\	
  {X_{4}}\\
  \vdots\\
  {X_{2m}}
  \end{array}
  \right)=\phi(2)  
  \left(
  \begin{array}{c}
  {X_{1}}\\	
  {X_{3}}\\
  \vdots\\
  {X_{2m-1}}
  \end{array}\right)+
   \left(
  \begin{array}{c}
  {\epsilon^{H}_{2}}\\	
  {\epsilon^{H}_{4}}\\
  \vdots\\
 {\epsilon^{H}_{2m}}
  \end{array}\right)
\end{equation}
,
\begin{equation}
 \left(
  \begin{array}{c}
  {X_{3}}\\	
  {X_{5}}\\
  \vdots\\
  {X_{2m+1}}
  \end{array}
  \right)=\phi(1)  
  \left(
  \begin{array}{c}
  {X_{2}}\\	
  {X_{4}}\\
  \vdots\\
  {X_{2m}}
  \end{array}\right)+
   \left(
  \begin{array}{c}
  {\epsilon^{H}_{3}}\\	
  {\epsilon^{H}_{5}}\\
  \vdots\\
 {\epsilon^{H}_{2m+1}}
  \end{array}\right)
\end{equation}
and 
\begin{equation}
 \left(
  \begin{array}{c}
  {X_{2}}\\	
  {X_{3}}\\
  \vdots\\
  {X_{n}}\\
  {X_{n+1}}
  \end{array}
  \right)=   
  \left(
  \begin{array}{cc}  % cc = 第一列居中，第二列居中（核心）
    0 & \ X_{1} \\  % 第一列放“0”，第二列放“逗号+X₁”，统一分隔符位置
    X_{2} & \ 0 \\   % 第一列放“X₂”，第二列放“逗号+0”，保持对称
    \vdots & \vdots \\  % 补充第二列省略号，避免上下不对称
   0 & \  X_{n-1}\\
    X_{n} & \ 0
  \end{array}\right)\left(
  \begin{array}{c}
 \phi(1)\\	
  \phi(2)
  \end{array}\right)+
   \left(
  \begin{array}{c}
  {\epsilon^{H}_{2}}\\	
  {\epsilon^{H}_{3}}\\
  \vdots\\
 {\epsilon^{H}_{n+1}}
  \end{array}\right)
  \label{expressionofalldata}
\end{equation}
  where $\epsilon^{H}_{t}$ is the fractional Gaussian noise with Hurst index $H$ and $n \in \mathbb{N}$. It is worth noting that in Equation (\ref{expressionofalldata}), the specific form of the expression may vary slightly depending on whether $n$ is odd or even. However, this distinction has no impact on the subsequent parameter estimation procedures or the validity of the theoretical proofs. Thus, without loss of generality, we assume \(n = 2m\) (where $m$ is a positive integer) to be an even number in the following analysis. One can emphasize that the random vector $(\epsilon^{H}_{2}, \epsilon^{H}_{3},\dots,\epsilon^{H}_{n})^{\top}$ (${\top}$ is the transpose of the vector) is  centered normal random vector with covariance matrix $\Omega_{n,H} = [\rho _{\epsilon^{H}}(j-i)]_{1\leq i,j \leq n} = [Cov(\epsilon^{H}_{i}, \epsilon^{H}_{j})]_{1\leq i,j \leq n}$ and $\Sigma_{n,H} = [\rho _{\epsilon^{H}}(2j-2i)]_{1\leq i,j \leq n} = [Cov(\epsilon^{H}_{2i}, \epsilon^{H}_{2j})]_{1\leq i,j \leq n}$.
\begin{remark}
   Since the definition of $\rho_{\epsilon^{H}(j-i)}$, we know that $\Omega_{n,H}$ is a symmetric real matrix that can be diagonalized by an orthogonal matrix. In clearer terms, there exists a $n \times n $ diagonal matrix $D_{n,H}$ satisfies $D_{n,H} = P^\top_{n,H}\Omega_{n,H}P_{n,H}$ where $P_{n,H}$ is an orthogonal matrix. By taking into consideration the positive definition of the matrix $\Omega_{n,H}$, we have $\Omega^{\frac{1}{2}}_{n,H} = P^\top_{n,H}D^{\frac{1}{2}}_{n,H}P_{n,H}$.
\end{remark}
\par Due to the standard form of the generalized least squares estimator, we will consider a straightforward transformation for the sequences $(\epsilon^{H}_{n})_{n\in \mathbb{N}}$, $(X_{2n})_{n \in \mathbb{N}}$  and $(X_{2n+1})_{n \in \mathbb{N}}$ in vector form. More specifically, we have\\
\begin{equation}
\textbf{X}^{(1)}_{n,H}= \Sigma^{-\frac{1}{2}}_{n,H} \left(
  \begin{array}{c}
  {X_{1}}\\	
  {X_{3}}\\
  \vdots\\
  {X_{2n-1}}
  \end{array}
  \right), \quad
 \textbf{X}^{(2)}_{n,H}= \Sigma^{-\frac{1}{2}}_{n,H}
  \left(
  \begin{array}{c}
  {X_{2}}\\	
  {X_{4}}\\
  \vdots\\
  {X_{2n}}
  \end{array}\right) 
  , \quad
 \textbf{X}_{n,H}= \Omega^{-\frac{1}{2}}_{n,H}
 \left(
  \begin{array}{cc}  % cc = 第一列居中，第二列居中（核心）
    0 & \ X_{1} \\  % 第一列放“0”，第二列放“逗号+X₁”，统一分隔符位置
    X_{2} & \ 0 \\   % 第一列放“X₂”，第二列放“逗号+0”，保持对称
    \vdots & \vdots \\  % 补充第二列省略号，避免上下不对称
    0 & \ X_{n}
  \end{array}
\right) \label{eqofX}
\end{equation}
\\
\begin{equation}
  \textbf{U}^{(1)}_{n, H} =\Sigma^{-\frac{1}{2}}_{n,H}\left(
  \begin{array}{c}
  {\epsilon^{H}_{1}}\\	
  {\epsilon^{H}_{3}}\\
  \vdots\\
 {\epsilon^{H}_{2n-1}}
  \end{array}\right),\quad
  \textbf{U}^{(2)}_{n, H} =\Sigma^{-\frac{1}{2}}_{n,H}\left(
  \begin{array}{c}
  {\epsilon^{H}_{2}}\\	
  {\epsilon^{H}_{4}}\\
  \vdots\\
 {\epsilon^{H}_{2n}}
  \end{array}\right) ,\quad
  \textbf{U}_{n, H} =\Omega^{-\frac{1}{2}}_{n,H}\left(
  \begin{array}{c}
  {\epsilon^{H}_{2}}\\	
  {\epsilon^{H}_{3}}\\
  \vdots\\
 {\epsilon^{H}_{n+1}}
  \end{array}\right),\quad
 \label{eq:formofu} 
\end{equation}
\begin{equation}
\textbf{Y}_{n,H}= \Omega^{-\frac{1}{2}}_{n,H} \left(
  \begin{array}{c}
  {X_{2}}\\	
  {X_{3}}\\
  \vdots\\
  {X_{n+1}}
  \end{array}
  \right),\quad \bm{\phi} =\left(
  \begin{array}{c}
  \phi(1)\\	
 \phi(2)
  \end{array}\right).
\end{equation}
\\
 
\subsection{\texorpdfstring{Properities of the Components of the Matrix $\Omega_{n,H}$ and $\Omega^{-1}_{n,H}$}{Properties of the components of the matrix Omega n,H and Omega inverse n,H}}
\par Thanks to the chapter 4 in \citep{esstafa2019long}. We know that the elements of $\Omega^{-1}_{n,H}$ can be expressed as a function of the spectral density of fGn. the spectral representation of 
 $(\Omega^{-1}_{n,H})_{j,k}$ implies that 
  \begin{equation}
  	(\Omega^{-1}_{n,H})_{j,k}=\frac{1}{(2\pi)^{2}}\int_{-\pi}^{\pi}\frac{1}{f_{\epsilon_{n}^{H}}(\lambda)}e^{i(k-j)\lambda}d\lambda,  
  \end{equation}
 As $\lambda \rightarrow 0$, according to the definition of fractional Gaussian noise, we have
\begin{equation}
   f_{\epsilon_{n}^{H}}(\lambda) \sim \frac{C_{H}}{2} |\lambda|^{1-2H}.  
\end{equation}
We can categorize the elements of the matrix into two types: diagonal elements and off-diagonal elements.    
  \par When $j = k$, we have
\begin{equation}
   (\Omega^{-1}_{n,H})_{j,j} = \frac{1}{(2\pi)^{2}} \int_{-\pi}^{\pi} \frac{1}{f_{\epsilon_{n}^{H}}(\lambda)} \, d\lambda = \frac{1}{2 \pi^{2}} \int_{0}^{\pi} \frac{1}{f_{\epsilon_{n}^{H}}(\lambda)} \, d\lambda. \label{2.18}
\end{equation}
One has when $\lambda \rightarrow 0$ that
\begin{equation}
   \frac{1}{f_{\epsilon_{n}^{H}}(\lambda)} = \frac{2}{C_{H}} |\lambda|^{2H-1} + o\left(\frac{2}{C_{H}} |\lambda|^{2H-1}\right).  
\end{equation}
This implies that for $l > 0$ there exists $\delta_{l} > 0$ such that for any $\lambda \in (-\delta_{l}, \delta_{l})$, we have
\begin{equation}
   (1 - l) \frac{2}{C_{H}} |\lambda|^{2H-1} \leq \frac{1}{f_{\epsilon_{n}^{H}}(\lambda)} \leq (1 + l) \frac{2}{C_{H}} |\lambda|^{2H-1}. \label{2.20}
\end{equation}
Thus, equation \eqref{2.18} and  equation \eqref{2.20} have an upper bound when $\lambda \in (-\delta_{l}, \delta_{l})$:
\begin{eqnarray}
   |(\Omega^{-1}_{n,H})_{j,j}| 
   &\leq& \frac{1 + l}{C_{H} \pi^{2}} \int_{0}^{\delta_{l}} \lambda^{2H-1} \, d\lambda + \frac{1}{2 \pi^{2}} \int_{\delta_{l}}^{\pi} \frac{1}{f_{\epsilon_{n}^{H}}(\lambda)} \, d\lambda \\
   &\leq& \frac{\delta_{l}^{2H} (1 + l)}{2H C_{H} \pi^2} + \frac{\pi - \delta_{l}}{2 \pi^{2}} \sup_{\lambda \in (\delta_{l}, \pi]} \frac{1}{f_{\epsilon_{n}^{H}}(\lambda)} \nonumber \\
   &\leq& K \nonumber, \label{XXX}
\end{eqnarray}
where $K$ is a constant.
\par When $j \neq k$, according to \citep{esstafa2019long}, there exists a positive constant $K$ and $T$ such that for any $j, k = 1, 2, \ldots$
\begin{equation}
   |(\Omega^{-1}_{n,H})_{j,k}| \leq K \left| \frac{1}{k-j} \right|^{2H}. \label{invrho}
\end{equation}
According to the equation (\ref{2.2}), we have
\begin{equation}
(\Omega_{n,H})_{i,j} = \rho_{\epsilon^{H}}(i-j).
\end{equation}
For large $k$, the asymptotic behavior of $\rho_{\epsilon^{H}}(k)$ is given by
\begin{equation}
   \rho_{\epsilon^{H}}(k) \sim H(2H-1)k^{2H-2}+o(k^{2H-2}) \label{rho}
\end{equation}
 
\subsection{\texorpdfstring{The Expression of $\textbf{X}^{(1)}_{n, H}$ and $\textbf{X}^{(2)}_{n, H}$ as Functions of $\textbf{U}^{(1)}_{n, H}$ and $\textbf{U}^{(2)}_{n, H}$}{The expression of Y n,H and Z n,H as functions of U(1) n,H and U(2) n,H}}
  Under assumption ${\mathcal{A}_{0}}$, we could write the process $\left\{X_{2n+1}\right\}_{n \in \mathbb{N}}$
  and $\left\{X_{2n+2}\right\}_{n \in \mathbb{N}}$ as a linear combination of the infinite fractional Gaussian noises $\epsilon^{H}_{n}$. More explicitly,  $X_{2n+1}$ and $X_{2n+1}$ take the following expression 
   \begin{equation}
        X_{2n+2} = \sum^{\infty}_{i=0,2,4,...}\epsilon^{H}_{2n+2-i}[\phi(2)\phi(1)]^{\frac{i}{2}}+\sum^{\infty}_{i=1,3,5,...}\epsilon^{H}_{2n+2-i}[\phi(2)\phi(1)]^{\frac{i-1}{2}}\phi(2),
        \label{X2INFINITESUMOFFGN}
        \end{equation}
        \begin{equation}
        X_{2n+1} = \sum^{\infty}_{i=0,2,4,...}\epsilon^{H}_{2n+1-i}[\phi(2)\phi(1)]^{\frac{i}{2}}+\sum^{\infty}_{i=1,3,5,...}\epsilon^{H}_{2n+1-i}[\phi(2)\phi(1)]^{\frac{i-1}{2}}\phi(1),
         \label{X1INFINITESUMOFFGN}
    \end{equation}
From equations (\ref{X2INFINITESUMOFFGN}) and (\ref{X1INFINITESUMOFFGN}), the vectors $(X_{1}, X_{2},...,X_{2n-1})$ and $(X_{2}, X_{4},...,X_{2n})$ have the following form
 \begin{equation}
  \left(
  \begin{array}{c}
  {X_{1}}\\	
  {X_{3}}\\
  \vdots\\
  {X_{2n-1}}
  \end{array}
  \right)=
  \left(
  \begin{array}{c}
   \sum^{\infty}_{i=0,2,4,...}\epsilon^{H}_{1-i}[\phi(2)\phi(1)]^{\frac{i}{2}}+\sum^{\infty}_{i=1,3,5,...}\epsilon^{H}_{1-i}[\phi(2)\phi(1)]^{\frac{i-1}{2}}\phi(1)\\	
  \sum^{\infty}_{i=0,2,4,...}\epsilon^{H}_{3-i}[\phi(2)\phi(1)]^{\frac{i}{2}}+\sum^{\infty}_{i=1,3,5,...}\epsilon^{H}_{3-i}[\phi(2)\phi(1)]^{\frac{i-1}{2}}\phi(1)\\
  \vdots\\
   \sum^{\infty}_{i=0,2,4,...}\epsilon^{H}_{2n-1-i}[\phi(2)\phi(1)]^{\frac{i}{2}}+\sum^{\infty}_{i=1,3,5,...}\epsilon^{H}_{2n-1-i}[\phi(2)\phi(1)]^{\frac{i-1}{2}}\phi(1)
    \end{array}
  \right) \nonumber
\end{equation}
 \begin{equation}
 =\sum^{\infty}_{j=0}[\phi(1)\phi(2)]^{j}\mathcal{B}^{2j}
\left(
  \begin{array}{c}
  {\epsilon^{H}_{1}}\\	
  {\epsilon^{H}_{3}}\\
  \vdots\\
 {\epsilon^{H}_{2n-1}}
  \end{array}\right)+\phi(1)
  \sum^{\infty}_{k=0}[\phi(1)\phi(2)]^{k}\mathcal{B}^{2k+1}
\left(
  \begin{array}{c}
  {\epsilon^{H}_{1}}\\	
  {\epsilon^{H}_{3}}\\
  \vdots\\
 {\epsilon^{H}_{2n-1}}
  \end{array}\right) \label{X1FORMEXPLICIT}
 \end{equation}
 and
 \begin{equation}
  \left(
  \begin{array}{c}
  {X_{2}}\\	
  {X_{4}}\\
  \vdots\\
  {X_{2n}}
  \end{array}
  \right)=
  \left(
  \begin{array}{c}
   \sum^{\infty}_{i=0,2,4,...}\epsilon^{H}_{2-i}[\phi(2)\phi(1)]^{\frac{i}{2}}+\sum^{\infty}_{i=1,3,5,...}\epsilon^{H}_{2-i}[\phi(2)\phi(1)]^{\frac{i-1}{2}}\phi(2)\\	
  \sum^{\infty}_{i=0,2,4,...}\epsilon^{H}_{4-i}[\phi(2)\phi(1)]^{\frac{i}{2}}+\sum^{\infty}_{i=1,3,5,...}\epsilon^{H}_{4-i}[\phi(2)\phi(1)]^{\frac{i-1}{2}}\phi(2)\\
  \vdots\\
   \sum^{\infty}_{i=0,2,4,...}\epsilon^{H}_{2n-i}[\phi(2)\phi(1)]^{\frac{i}{2}}+\sum^{\infty}_{i=1,3,5,...}\epsilon^{H}_{2n-i}[\phi(2)\phi(1)]^{\frac{i-1}{2}}\phi(2)
    \end{array}
  \right) \nonumber
\end{equation}
 \begin{equation}
 =\sum^{\infty}_{j=0}[\phi(1)\phi(2)]^{j}\mathcal{B}^{2j}
\left(
  \begin{array}{c}
  {\epsilon^{H}_{2}}\\	
  {\epsilon^{H}_{4}}\\
  \vdots\\
 {\epsilon^{H}_{2n}}
  \end{array}\right)+\phi(2)
  \sum^{\infty}_{k=0}[\phi(1)\phi(2)]^{k}\mathcal{B}^{2k+1}
\left(
  \begin{array}{c}
  {\epsilon^{H}_{2}}\\	
  {\epsilon^{H}_{4}}\\
  \vdots\\
 {\epsilon^{H}_{2n}}
  \end{array}\right) \label{X2FORMEXPLICIT}
 \end{equation}
where $\mathcal{B}^{i}$ is the lag shift that acts on all the fGn vector, i.e. $\mathcal{B}^{i}(\epsilon^{H}_{1}, \epsilon^{H}_{2},...,\epsilon^{H}_{n})=(\epsilon^{H}_{1-j}, \epsilon^{H}_{2-j},...,\epsilon^{H}_{n-j}).$ In view of equations (\ref{X1FORMEXPLICIT}) and (\ref{X2FORMEXPLICIT}), vectors $\textbf{Y}_{n,H}$ and $\textbf{Z}_{n,H}$ can be expressed as functions of $\textbf{U}^{(1)}_{n,H}$ and $\textbf{U}^{(2)}_{n,H}$:
\begin{equation}
  \textbf{X}^{(1)}_{n,H} =\sum^{\infty}_{j=0}[\phi(1)\phi(2)]^{j}\mathcal{B}^{2j}
  \textbf{U}^{(1)}_{n,H}+\phi(1)\sum^{\infty}_{k=0}[\phi(1)\phi(2)]^{k}\mathcal{B}^{2k+1}\textbf{U}^{(1)}_{n,H} \label{eq:formofY}
\end{equation}
\begin{equation}
  \textbf{X}^{(2)}_{n,H} =\sum^{\infty}_{j=0}[\phi(1)\phi(2)]^{j}\mathcal{B}^{2j}
  \textbf{U}^{(2)}_{n,H}+\phi(2)\sum^{\infty}_{k=0}[\phi(1)\phi(2)]^{k}\mathcal{B}^{2k+1}\textbf{U}^{(2)}_{n,H}\label{eq:formofZ}
\end{equation}
The expression of $\textbf{Y}_{n, H}$ can be directly obtained from Equations (\ref{eq:formofY}) and (\ref{eq:formofZ}) via $\epsilon^{H}_{n}$.
\section{Initial Estimators of PFAR(1) Models}

\subsection{The GPH Estimator for the Hurst Index}
  Due to the non-stationarity of $X_{n}$, obtaining an estimator for $H$ using standard semiparametric methods is not feasible. To address this, we can extract stationarity from the data by splitting the time series $(X_{n})_{n \in \mathbb{N}}$ into periodic components, resulting in $T$ stationary subsequences $\textbf{X}(u)=(X_{u},X_{T+u},...,X_{nT+u}) $ and we construct a stationary additive series defined as $Z_{n} = \sum^{T}_{u=1} X_{nT+u}$. 
  \par In this subsection, we will estimate $H$  using the log-periodogram method, specifically the GPH estimator, applied to the additive series$(Z_{n})_{n \in \mathbb{N}}$. The spectral density and stationarity properties of $(Z_{n})_{n \in \mathbb{N}}$ and $(X_{nT+u})_{u\in \mathbb{Z}}$  are outlined in the following three propositions. 
  \par Let new series $(Z_{n})_{n \in \mathbb{N}}$ be an observation sample generated via the equation (\ref{eq:FormofPFAR(1)}) and choose a suitable integer m which can decrease the mean square error of estimation, where $m < n$. we get the periodogram of ${Z_{n}}$ given by 

  \begin{equation}
  	I(\lambda)=\frac{1}{2\pi n}{|\sum_{t=1}^{n}Z_{t}exp(it\lambda)|^{2}}, 
  \end{equation}
  \begin{equation}
  	\lambda_{j}=\frac{2\pi j}{n},\quad j{\in}\left \{1,2,...m\right\},  
  \end{equation}
  \begin{equation}
  	a_{j}=\log(2sin\frac{\lambda_{j}}{2}),\quad {\overline {a}_m}=\frac {1}{m}\sum_{j=1}^{m}a_{j},\quad S_{m}=\sum_{j=1}^m(a_j-{\overline {a}_m})^2. 
  \end{equation}
  
  \noindent We estimate d by regressing $\log I({\lambda_{j}})$ with respect to $a_{j}$, such that
  \begin{equation}
  	{\hat{d}}_{n} = - \frac{1}{2S_{m}} \sum_{j=1}^m(a_j-{\overline {a}_m})  \log I(\lambda_{j}),  
  \end{equation}
  
  \noindent The estimator $\hat{H}_{n}$ is defined by
  \begin{equation}
  	\hat{H}_{n}= {\hat{d}_{n}}+\frac{1}{2},  \label{GPH} 
  \end{equation}

\begin{remark}
  There are several semi-parametric methods for estimating the long memory parameter $d$ and $H$, such as whittle estimation and R/S estimation method proposed by \citep{robinson1995log} and \citep{marinucci1998semiparametric}. These models rely on the log-periodogram approach. However, these methods tend to exhibit greater bias compared to the GPH estimator.

\end{remark}

\subsection{\texorpdfstring{The GLSE for $\phi(1)$ and $\phi(2)$}{The GLSE for phi(1) and phi(2)}}
 \par We now focus on estimating $\phi(u)$ given that the parameter $H$  has been estimated. When the noise in the periodic autoregressive model is white noise, we can easily obtain the estimator of the parameters of these models using Least Squares Estimation (LSE). However, when the noise is fractional Gaussian noise (fGn), the covariance matrix of fGn is no longer diagonal, making LSE inappropriate. Therefore, we consider using Generalized Least Squares Estimation (GLSE).
  \par To address the effect of periodic structure on parameter estimation, we apply GLSE to the sequences\\ 
$$\left( \begin{array}{cccccc}  % 改为n列，保持居中对齐
    0 & X_{2} & \cdots & 0 & X_{n} \\  % 原矩阵第一列成为转置矩阵第一行
    X_{1} & 0 & \cdots & X_{n-1} & 0  % 原矩阵第二列成为转置矩阵第二行
\end{array}
\right)^{\top},\quad \left( X_{2},X_{3} \dots ,X_{n+1} \right). $$ \\
This allows us to estimate the parameters $\phi(1),\phi(2)$, assuming the Hurst index is known.

\par The GLSE of ${\hat{{\phi}}}_{n}(1)$ and ${\hat{{\phi}}}_{n}(2)$ are defined by
\begin{equation}
 \hat{\phi}_n = \begin{pmatrix}
\hat{\phi}_n{(1)} \\
\hat{\phi}_n{(2)}
\end{pmatrix} = \left(\textbf{X} _{n,H}^\top  \textbf{X} _{n,H} \right)^{-1} \left( \textbf{X}_{n,H}^\top  \textbf{Y}_{n,H} \right)
\label{GLSEOFPHI12} 
\end{equation}
\begin{remark}
We can also construct an alternative pair of estimators for $\phi(1)$ and $\phi(2)$ by $\textbf{X}^{(1)}_{n, H}$ and $\textbf{X}^{(2)}_{n, H}$. Specifically, these estimators are given by  
\begin{equation}
 {\hat{{\phi}}}_{n}(2) =\frac{{\textbf{X}^{(1)}_{n,H}}^\top{\textbf{X}^{(2)}_{n,H}}}{{\left\lVert {\textbf{X}^{(1)}_{n,H}} \right\rVert}^{2}}, \quad
 {\hat{{\phi}}}_{n}(1) =\frac{({\textbf{X}^{(2)}_{n,H}}^\top)({{\mathcal{B}^{2}\textbf{X}^{(1)}_{n,H}}})}{{\left\lVert{\textbf{X}^{(2)}_{n,H}} \right\rVert}^{2}} \label{eq:IEESTIMATOROFPHI1PHI2}
\end{equation}
where $\mathcal{B}$ is lag operator and $||\cdot||$ is the $L_{2}$ norm.
 Moreover, in view of the equations (\ref{eq:IEESTIMATOROFPHI1PHI2}), we can deduce that
 \begin{equation}
 {\hat{{\phi}}}_{n}(2) - \phi(2) =\frac{{\textbf{X}^{(1)}_{n,H}}^\top{\textbf{U}^{(2)}_{n,H}}}{{\left\lVert {\textbf{X}^{(1)}_{n,H}} \right\rVert}^{2}}, \quad
 {\hat{{\phi}}}_{n}(1) -\phi(1) =\frac{({\textbf{X}^{(2)}_{n,H}}^\top)({{\mathcal{B}^{2}\textbf{U}^{(1)}_{n,H}}})}{{\left\lVert{\textbf{X}^{(2)}_{n,H}} \right\rVert}^{2}}.
\end{equation}
Regrettably, while this pair of estimators is asymptotically normal, they are biased. Specifically, they are unbiased when $H = 0.5$, right biased when $\frac{1}{2}<H < 1$, and left biased when $0<H < \frac{1}{2}$. Furthermore, the bias is a constant that depends on the parameters $\phi(1)$, and $\phi(2)$. We will present the simulation results of these two initial estimators in Appendix 4. The bias of the initial estimators may potentially be addressed through the Multi-Step estimator , however, this remains an open question.
\end{remark}
\subsection{Asymptotic Results of Initial Estimator}
The main asymptotic results of the initial estimator are divided into three major parts. In the first part, the stationarity of the subsequence $(X_{2n+1} )_{n \in \mathbb{N}}$, $(X_{2n+2})_{n \in \mathbb{N}}$ and the additive series  $(Z_{n})_{n \in \mathbb{N}}$ have been proved. In the second part, the spectral densities of the above three sequences are calculated. In the third part, the consistency and asymptotic normality of the estimators obtained in the second section are demonstrated.
\begin{theorem} \label{Thm:sationaryofX}
        	When T = 2, the subsequence $(X_{2n+1} )_{n \in \mathbb{N}}$ and $(X_{2n+2})_{n \in \mathbb{N}}$ satisfy $$X_{2n+2} = \phi(2)X_{2n+1}+\epsilon^{H}_{2n+2}, \quad X_{2n+1} = \phi(1)X_{2n}+\epsilon^{H}_{2n+1}, \quad n \in \mathbb{N},$$
         where $\sum^{\infty}_{j=0} |\phi(2)|^{j}$, $\sum^{\infty}_{j=0} |\phi(1)|^{j}$ is bounded. Thus, $(X_{2n+1})_{n \in \mathbb{N}}$ and $(X_{2n+2})_{n \in \mathbb{N}}$ are stationary processeses.
\end{theorem}
        \begin{remark}
        When $T \neq 2$, a more general constant, the above proposition still holds; however, the cross terms in the proof increase. Nevertheless, the property that the autocovariance function depends only on the time interval remains unchanged.
        \end{remark}
 After proving the stationarity of the sequences, we can calculate their spectral densities.
 \begin{remark}
 Since the sequences $(X_{2n+1})_{n \in \mathbb{N}}$ and $(X_{2n+2})_{n \in \mathbb{N}}$  are each stationary, it follows that the sequence $(Z_{n})_n \in \mathbb{N}$ is also stationary.
 \end{remark}
 \begin{proposition}\label{prop:sdofx}
        Let $\overline{f}_{H,\phi(u)}(\lambda)$ is the spectral density of  $\left\{X_{nT+u}\right\}_{n\in\mathbb{N}}$, we have
        	\begin{equation}
        	\overline{f}_{H,\phi(u)}(\lambda) = \left| \sum^{T-2}_{j=0}\Psi_{j}(\lambda)+\frac{1}{1-Ce^{-iT\lambda}}\right|^{2}	{f_{\epsilon_{n}^{H}}(\lambda)},   \nonumber
        	\end{equation}
        where ${f_{\epsilon_{n}^{H}}(\lambda)}$ is the spectral density of fGn, $u = 1,2,3,...,T$, $C = \phi(u)\phi(u-1)\cdots\phi(u-T+1)$,  $\Psi_{j}(\lambda) = \frac{\phi(u)\phi(u-1)\cdots\phi(u-j)}{1-Ce^{-iT\lambda}}$. 
\end{proposition} 

\begin{remark}
  This proposition holds for any finite value of $T$
. Consequently, it also holds when $T=2$.
\end{remark}

\begin{proposition}\label{prop:sdofz}
   	Let $p_{H,\phi(1),\phi(2)}(\lambda)$ be the spectral density of $\left\{Z_{n} \right\}_{n \in \mathbb{N}}$, when $T=2$, we have
   		 $$p_{H,\phi(2),\phi(1)}(\lambda) = \left|1+\frac{\phi(2)(1+\phi(1))e^{-2\lambda i}}{1-e^{-2\lambda i}\phi(2)\phi(1)}+\frac{e^{-\lambda i}(1+\phi(1))}{1-e^{-2\lambda i}\phi(2)\phi(1)}\right|^{2}f_{\epsilon^{H}_{n}}(\lambda), $$ 
   	  where $f_{\epsilon^{H}_{n}}$ is the spectral density of fractional Gaussian noise.
   \end{proposition}
    
 When studying the problem of parameter estimation for  
 $\hat{\phi}(1)$ and  $\hat{\phi}(2)$. we are committed to exploring the properties of estimators, such as consistency and asymptotic distribution, so as to provide a theoretical basis for statistical inference. Under the premises of assumptions \((\mathcal{A}_{0})\) and \((\mathcal{A}_{1})\), a series of key lemmas lays a solid foundation for the proofs of subsequent theorems.
\begin{lemma}\label{lemma:denominatortopdm}
Suppose that Assumptions $(\mathcal{A}_0)$ and $(\mathcal{A}_1)$ hold. Let $\mathbf{X}_{n,H}$ denote the matrix  associated with parameter $H$ and sample size $n$. Then, the expected value of the corresponding sample moment matrix converges as $n \to \infty$, i.e.,
\[
\mathbb{E}\left[ \frac{1}{n} \mathbf{X}_{n,H}^\top \mathbf{X}_{n,H} \right] \xrightarrow[n \to \infty]{} Q,
\]
where $Q \in \mathbb{R}^{2 \times 2}$ is a constant positive definite matrix.
\end{lemma}
\begin{remark}
Under assumption $(\mathcal{A}_{0})$ and $(\mathcal{A}_{1})$, we have
$\frac{1}{n}\left\|\mathbf{X}^{(k)}_{n,H}\right\|^2 \to C^{(k)}_{\phi(1),\phi(2)}$ for \( k = 1,2 \), where \( C^{(k)}_{\phi(1),\phi(2)} \) are constants depending on \( \phi(1) \), and \( \phi(2) \) .
\end{remark}
\begin{remark}
     Denote by $\Sigma^{(X)}_{n,H}$ the variance matrix of the random vector $(X_{2},...,X_{n+1})^\top$, i.e. for any $1 \leq i,j\leq n$, $(\Sigma^{(X)}_{n,H})_{i,j} = COV(X_{i},X_{j})$, let 
\begin{equation}
\textbf{Y}^{(X)}_{n,H} = {\Sigma^{(X)}}^{-\frac{1}{2}}_{n,H}(X_{2}, X_{3},...,X_{n+1})^\top.
\end{equation}
    We have
\begin{equation}
\frac{1}{n}\left\lVert \textbf{X}_{n,H} \right\rVert^{2} = \frac{1}{n} \textbf{X}^\top_{n,H}\textbf{X}_{n,H} = \frac{1}{n}\left( \textbf{Y}^{(X)}_{n,H}\right)^\top\left( \left( \Sigma^{(X)}_{n,H}\right)^{\frac{1}{2}} \right)^\top\Omega^{-1}_{n,H}\left( \Sigma^{(X)}_{n,H}\right)^{\frac{1}{2}}\textbf{Y}^{(X)}_{n,H}
\end{equation}
Since $\left( \left( \Sigma^{(X)}_{n,H}\right)^{\frac{1}{2}} \right)^\top\Omega^{-1}_{n,H}\left( \Sigma^{(X)}_{n,H}\right)^{\frac{1}{2}} $ is a symmetric real matrix. So $\frac{1}{n}\left\lVert \textbf{X}_{n,H} \right\rVert^{2}$ has the following form
\begin{equation}
\frac{1}{n}\left\lVert \textbf{X}_{n,H} \right\rVert^{2} = 
\frac{1}{n}\sum^{n}_{j=1}\eta^{n,H}_{j} \left[ \left(W^\top_{n,H}\textbf{Y}^{(X)}_{n,H}\right)_{j} \right]^{2},\label{YEE}
\end{equation}
where $\eta^{n,H}_{1},\eta^{n,H}_{2}...,\eta^{n,H}_{n}$ are the eigenvalues of the matrix  $\left( \left( \Sigma^{(X)}_{n,H}\right)^{\frac{1}{2}} \right)^\top\Omega^{-1}_{n,H}\left( \Sigma^{(X)}_{n,H}\right)^{\frac{1}{2}}$, $W^\top_{n,H}$ is an orthogonal matrix such that $W^\top_{n,H}\left( \left( \Sigma^{(X)}_{n,H}\right)^{\frac{1}{2}} \right)^\top\Omega^{-1}_{n,H}\left( \Sigma^{(X)}_{n,H}\right)^{\frac{1}{2}}W_{n,H}$ is a diagonal matrix. $W^\top_{n,H}\textbf{Y}^{(X)}_{n,H}$ is the standard normal random vector that will help us simplify the calculations in the following.
\end{remark}
\begin{lemma}\label{lemma:etabound}
Under assumption $(\mathcal{A}_{0})$ and $(\mathcal{A}_{1})$, for any $1\leq j \leq n$, when $n \rightarrow \infty$, we have
$$ D^{(1)}_{H,\phi(1),\phi(2)} \leq \eta^{n,H}_{j} \leq D^{(2)}_{H,\phi(1),\phi(2)},$$
where $D^{(1)}_{H,\phi(1),\phi(2)}$ and $D^{(2)}_{H,\phi(1),\phi(2)}$ are constants related to $H$, $\phi(1)$ and $\phi(2)$.
\end{lemma}
After determining the bounds of $\eta^{n,H}_{j}$, it is essential to further explore the properties of the relevant random vectors under limiting conditions. This exploration will enable us to gain a more profound understanding of the behavioral characteristics of the estimators $\hat{\phi}_{n}(1)$ and  $\hat{\phi}_{n}(2)$, leading to the following lemma.

\begin{lemma}\label{lemma:nominatortozero}
Under assumption $(\mathcal{A}_{0})$ and $(\mathcal{A}_{1})$, we have 
\[
\mathbb{E}\left[ \frac{1}{n} \mathbf{X}_{n,H}^\top\textbf{U}_{n,H} \right] \xrightarrow[n \to \infty]{} {\bf{0}},
\]
where ${\bf{0}}\in \mathbb{R}^{2 \times 1}$ denotes a zero column vector.
\end{lemma}
 \begin{remark}
 We can further show that
\[
\frac{1}{n} \mathbf{X}_{n,H}^\top \mathbf{U}_{n,H} \xrightarrow[n \to \infty]{L^2} \mathbf{0},
\]
and the proof of this result follows a similar line of reasoning to that of Lemma \ref{lemma:nominatortozero}.
 \end{remark}
\begin{lemma}\label{lemma:nominatortonormal}
Suppose Assumptions $(\mathcal{A}_0)$ and $(\mathcal{A}_1)$  hold. Let $\mathbf{X}^{\top}_{n,H}$ and $\mathbf{U}_{n,H}$ as above. Then, as the sample size $n \to \infty$, the scaled random vector satisfies the following asymptotic normality 
\[
\frac{1}{\sqrt{n}} \mathbf{X}^{\top}_{n,H} \mathbf{U}_{n,H} \xrightarrow[n \to \infty]{{\mathcal{L}}} \mathcal{N}\left( {\bf{0}}, \Lambda_{\phi_1,\phi_2} \right),
\]
where ${\bf{0}}\in \mathbb{R}^{2 \times 1}$ and $\Lambda_{\phi(1),\phi(2)} \in \mathbb{R}^{2 \times 2}$ is a constant matrix whose depend on the parameters $\phi(1)$ and $\phi(2)$. \\
%\textcolor{red}{I think there are two approaches to  proving that the joint distribution of the estimators of $\phi(1)$ and $\phi(2)$ is asymptotically normal: the first is a direct proof, because joint distribution can alternatively be expressed as a linear transformation of the one in Esstafa. And the second is to first prove that each estimator is asymptotically normal individually, then demonstrate that any linear combination of them is also asymptotically normal. Thereby demonstrating that its joint distribution is asymptotically normal is still needed to consider.}%
\end{lemma}
\begin{remark}
This lemma is intended to illustrate that the joint distribution of the estimators of $\phi(1)$ and  $\phi(2)$ is asymptotically normal. Furthermore, the estimators of $\phi(1)$ and $\phi(2)$ are also each asymptotically normal individually.
\end{remark}
\begin{remark}
 The result derived from Lemma \ref{lemma:nominatortozero} is invalid for $$\frac{1}{n}{\textbf{X}^{(1)}_{n,H}}^\top{\textbf{U}^{(2)}_{n,H}}.$$ However, Lemma \ref{lemma:nominatortonormal} holds for $ \frac{1}{\sqrt{n}}{\textbf{X}^{(1)}_{n,H}}^\top{\textbf{U}^{(2)}_{n,H}}$, except that in this case, the latter follows a non-central normal distribution.
\end{remark}
\begin{theorem}  \label{Thm：consistencyofie}
  \par  Letting  $m = [Kn^{\delta}]$ for some $\frac{1}{2}<\delta < 1$, assume that $(X_{2n+u})_{n \geq 0}$ satisfies the equation (\ref{origneq}) for $u =1,2$. Under conditions $(\mathcal{A}_{0})$ and $(\mathcal{A}_{1})$, we have  
$$\left(
  \begin{array}{c}
  \hat{H}_{n}\\
  {\hat{\phi}_{n}(1)}	\\
  {\hat{\phi}_{n}(2)}
  \end{array}
  \right)\xrightarrow [n\rightarrow \infty]{\mathbb{P}}   \left(
  \begin{array}{c}
  	H\\
  {\phi(1)}\\	
  {\phi(2)}
  \end{array}\right),
  $$
  where $\hat{\phi}_{n}(u)$ are defined by equation (\ref{GLSEOFPHI12}), $\hat{H}_{n}$ is given by equations (\ref{GPH}) and $\mathbb{P}$ stands for convergence in probability.
\end{theorem} 

\begin{theorem} \label{Thm：asynofie}
      Let $m =[n^{\delta}]$ for some $\frac{1}{2}<{\delta}<\frac{2}{3}$. $\hat{\phi}_{n}(1)$, $\hat{\phi}_{n}(2)$ and $\hat{H}_{n}$  has a three dimensions limiting normal distribution given by 	 
   $$\sqrt{m}\left(
  \begin{array}{c}
  	\hat{H}_{n}-H\\
  {\hat{\phi}_{n}(1)}-{\phi(1)} 
  \\
  {\hat{\phi}_{n}(2)}-{\phi(2)} 
  \end{array}
  \right)\xrightarrow [n\rightarrow \infty]{\mathcal{L}}  \mathcal{N}(0,\Xi_{H,\phi(1),\phi(2)}),
  $$
  The covariance matrix $\Xi_{\phi(1),\phi(2)}$ is of the form $\Xi_{\phi(1),\phi(2)} = V_{H}\widetilde{\Xi}_{\phi(1),\phi(2)}$
  	and $V_{H}$	is the asymptotic variance of $\sqrt{m}({\hat{H}}_{n}-H)$, $\widetilde{\Xi}_{\phi(1),\phi(2)}$ is a built-in singular matrix.
  \end{theorem}
 \begin{remark} {\citep{hariz2024fast} represents \(\frac{1}{2} < \delta < \frac{2}{3}\), and \cite{hurvich1998mean} states that if \(m = n^{\delta}\), where \(0 < \delta < 1\), it can ensure the asymptotic normality of  \(\hat{H}_{n}\). The condition \(\frac{1}{2}< \delta\) is to ensure that \(\hat{\phi}_{n}(1)\) and \(\hat{\phi}_{n}(2)\) are asymptotically normal. But according to \cite{kutoyants2016multi}, if \(\frac{2}{3} < \delta\), a multi-step estimator may be required, which contradicts our consideration of a One-Step estimator. Thus, we consider restricting \(\delta\) to the interval \((\frac{1}{2},\frac{2}{3})\).}
  \end{remark}
  \begin{remark}
  	These results can be extented to the PFAR(p). 
  \end{remark}
  \begin{remark}
  Additionally, the estimation of $ {\hat{\phi}_{n}(1)}$ and $ {\hat{\phi}_{n}(2)}$  can be also approached using methods from \cite{brouste2014asymptotic} and \cite{soltane2024asymptotic}.	 
  \end{remark}
  \section{One-Step Estimator of PFAR(1) Models}	
  \par In this section, we explore modifications to the initial estimator \(\hat{\boldsymbol\theta}_n\) to develop a One-Step estimator \(\tilde{\boldsymbol\theta}_n\).
  \subsection{Construction of One-Step Estimator}
  \par We assume that $(Z_{n})_{n\in\mathbb{N}}$ is stationary with a spectral density $p_{H,\phi(1),\phi(2)}(\lambda)$, as obtained in propositon 3.3. For $p_{H,\phi(1),\phi(2)}(\lambda)$ to satisfy the necessary regularity conditions (see \cite{hariz2024fast} section 2.1), the following lemma must be met.\\
    \begin{lemma} \label{Lemma：regularityofsdofz}
      Under the hypothesis on the parametric space, $p_{H,\phi(1),\phi(2)}(\lambda)$ have the following results\\
   (1) For any $H \in [0,1]$ and $j \in \left\{0,1,2,3 \right\}$, $\frac{\partial}{\partial \lambda } \frac{\partial^{j}}{\partial^{j}H}p_{H,\phi(1),\phi(2)}(\lambda)$.\\
   (2) For any $j \in \left\{0,1,2,3 \right\} $ the functions $\frac{\partial^{j}}{\partial^{j}H}p_{H,\phi(1),\phi(2)}(\lambda)$ are symmetric with respect to  $\lambda$.\\
   (3) For any $\delta > 0 $ and all $(H,\lambda) \in [0,1] \times  [-\pi,\pi]  \backslash\left\{0 \right\}$
    \par a.$	C_{1,\delta} |\lambda|^{1-2H+\delta }\leq p_{H,\phi(1),\phi(2)}(\lambda) \leq C_{2,\delta}|\lambda|^{1-2H-\delta }.$ 
    \par b.$|\frac{\partial}{\partial \lambda}p_{H,\phi(1),\phi(2)}(\lambda)| \leq C_{3,\delta}|\lambda|^{-2H-\delta}.$
    \par c.For any $j \in \left\{0,1,2,3 \right\} $, $|\frac{\partial^{j}}{\partial^{j} H}p_{H,\phi(1),\phi(2)}(\lambda)| \leq C_{4,\delta}|\lambda|^{-2H-\delta}.$
     \end{lemma}

  	We let $l_{n}$ be the log-likelihood function of a stationary process $(Z_{n})_{n\in \mathbb{N}}$. We assume that $p_{H,\phi(1),\phi(2)}(\lambda)$ satisfies the regularity conditions and let $B({\theta, R})$ is an open ball of center $\theta$ and radius R for some $R > 0$. For any $t\in B(\boldsymbol\theta,R)$, $u \in \mathbb{N}$
  	\begin{equation}
  		l_{n}(\boldsymbol\theta + \frac{t}{\sqrt{n}})-l_{n}(\boldsymbol\theta) = t\frac{\nabla l_{n}(\boldsymbol\theta)}{\sqrt{n}}-\frac{t \mathcal{I}(\boldsymbol\theta)t^\top}{2} + r_{n,\boldsymbol\theta}(t),  	
  	\end{equation}
  	\noindent where, under ${\mathbb{P}}_{\boldsymbol\theta}^{(n)}$, when $n \rightarrow \infty $, the score function         $\nabla(\cdot)$ satisfies 	
    \begin{equation}
  		\frac{\nabla l_{n}(\boldsymbol\theta)}{\sqrt{n}}\xrightarrow [n\rightarrow\infty]{\mathbb{P}} \mathcal{N}(0,\mathcal{I}(\boldsymbol\theta)),     
    \end{equation}
  	\noindent and
    \begin{equation}
  		r_{n,\boldsymbol\theta}(t) \xrightarrow [n\rightarrow\infty]{a.s.} 0,   
    \end{equation}
  	\noindent uniformly on each compact set. The Fisher information matrix is given in our case by
  	\begin{equation}
  		\mathcal{I}(\boldsymbol\theta) = \frac{1}{4\pi}\left(\int_{-\pi}^{\pi}\frac{\partial log \,   p_{H,\phi(1),\phi(2)}(\lambda)}{\partial \theta_{k}}\frac{\partial log \,   p_{H,\phi(1),\phi(2)}(\lambda)}{\partial \theta_{j}}\right)_{1 \leq k,j \leq T+1},  \label{3.3}
  	\end{equation}
 
   where $\boldsymbol\theta=(\theta_{1},\theta_{1},\theta_{3}) = (H,\phi(1),\phi(2))$. This result is a direct consequence of Theorem from \citep{cohen:hal-00638121}.
  \par Since \( p_{H,\phi(1),\phi(2)}(\lambda) \) satisfies the regularity conditions, the elements of the Fisher information matrix \(\mathcal{I}(\boldsymbol\theta)\) are finite.  After obtaining the Fisher information matrix \(\mathcal{I}(\boldsymbol\theta)\) and the log-likelihood function of \((Z_n)_{n \in \mathbb{N}}\), we can compute the One-Step estimator as follows 
\begin{equation}
    \tilde{\boldsymbol\theta}_n = \hat{\boldsymbol\theta}_n + \mathcal{I}(\hat{\boldsymbol\theta}_n)^{-1} \frac{1}{n} \nabla l_n(\hat{\boldsymbol\theta}_n), \label{OSESTIMATOR}
\end{equation}
\subsection{Asymptotic Results of One-Step Estimator}
This subsection focuses on the asymptotic results of One-Step estimators. Through a series of lemmas, the asymptotic normal distribution property of One-Step estimators is gradually derived.
First, the lemmas characterize the local properties of the information matrix by setting certain conditions. Then, they further reveal the relationship between the difference of the log-likelihood function and the information matrix. Additionally, for stochastic sequences satisfying specific conditions, the lemmas establish relevant conclusions regarding the difference of the log-likelihood function.\\ 
\begin{lemma}  \label{Lemma：regular1}
  Let $\boldsymbol\theta_{0} \in \Theta$, $\delta>0$, such that for any $\boldsymbol\theta \in B(\boldsymbol\theta_{0},\delta)$, it holds that  
    \begin{equation}
    ||\mathcal{I}(\boldsymbol\theta)-\mathcal{I}(\boldsymbol\theta_{0})||\leq K||\boldsymbol\theta-\boldsymbol\theta_{0}||,
    \end{equation}
where $K$ is some constant.
\end{lemma}
 
\begin{lemma}\label{Lemma：regular2}
     	For any $\boldsymbol\theta \in \Theta$, it follows from the distribution of the parameter $\boldsymbol\theta$ that
  \begin{equation}
     \frac{\Delta l_{n}(\boldsymbol\theta)}{\sqrt{n}}+\sqrt{n}\mathcal{I}(\boldsymbol\theta) = O_{\mathbb{P}}(1),
   \end{equation}
\end{lemma}

\begin{lemma}\label{Lemma：regular3}
     Let $\left\{ {\overline{\boldsymbol\theta}}_{n}\right\}_{n}$ be a stochastic sequence satisfying ${\overline{\boldsymbol\theta}}_{n} - \boldsymbol\theta = o_{\mathbb{P}}(1)$. Then, according to the distribution of parameter $\boldsymbol\theta$, for any $k>0$, it holds that
    \begin{equation}
      \frac{\Delta l_{n}({\overline{\boldsymbol\theta}}_{n})}{n} - \frac{\Delta l_{n}(\boldsymbol\theta)}{n} = O_{\mathbb{P}}(n^{k}({\overline{\boldsymbol\theta}}_{n}-\boldsymbol\theta)).
   \end{equation} 
\end{lemma}
Based on the groundwork laid by these lemmas, the theorem, under the given initial estimator and the requirement that the function satisfies the regularity conditions, proves that the One-Step estimator follows an asymptotic normal distribution as the sample size approaches infinity. 
\begin{theorem}\label{Thm:asnofos}
  Let ${\hat{\boldsymbol\theta}}_{n}$ denote the initial estimator of $\boldsymbol\theta$ given by \eqref{GPH} and \eqref{GLSEOFPHI12}. When $g_{H,\underline{\phi}}(\boldsymbol\lambda)$ satisfies the regularity conditions, the one-step estimator ${\tilde{\boldsymbol\theta}}_{n}$ of $\boldsymbol\theta$ satisfies
\begin{equation}
\sqrt{n}\left({\tilde{\boldsymbol\theta}}_{n}-\boldsymbol\theta\right)\xrightarrow[n\rightarrow\infty]{\mathcal{L}}\mathcal{N}\left(0,\mathcal{I}(\boldsymbol\theta)^{-1}\right),
\end{equation}
 \end{theorem} 
   \begin{remark}
    The parameter \(\boldsymbol\theta\) should not lie on the boundary of the parameter space \(\Theta^l\).
\end{remark}
\begin{remark}
   {The One-Step estimatior can be applied more generally even if the initial estimator \(\hat{\boldsymbol\theta}_n\) does not satisfy asymptotic normality. According to proposition 2.3 in Hariz (2025), if the initial estimator with convergence speed lower than $\sqrt{n}$ and the spectral density of time series meets the regular condition, then the One-Step estimator \(\tilde{\boldsymbol\theta}_n\) can still achieve asymptotic normality.}
\end{remark}
\begin{remark}
 One-step estimator can achieve Hájek's lower bound, thus it is  asymptotically efficient in the local minimax sense. We can find related conclusions in \cite{cohen:hal-00638121}.
\end{remark}
 \section{Simulation Study}   
    
    \par According to equation $Z_{n} = \sum^{T}_{u=1} X_{nT+u}$, the likelihood function based on the sample \(\boldsymbol{Z}^{(n)} = (Z_{0}, Z_{1}, \ldots, Z_{n-1})\) is given by
\begin{equation}
    l_{n}(\boldsymbol\theta) = -\frac{1}{2} \log \det \left( \Gamma^{Z}_{n}(\boldsymbol\theta) \right) - \frac{1}{2} \boldsymbol{Z}^{(n)^\top} \Gamma^{Y}_{n}(\boldsymbol\theta)  \boldsymbol{Z}^{(n)},  
\end{equation}
where \(\Gamma^{Z}_{n}(\boldsymbol\theta)\) is the covariance matrix of \( \boldsymbol{Z}^{(n)}\). For any \(K \in \mathbb{N}\),
\begin{equation}
     Cov (Z_{0}, Z_{k}) = \int_{-\pi}^{\pi} \exp(ik\lambda) \, g_{H, \boldsymbol{\phi}}(\lambda) \, d\lambda,  
\end{equation}
where \(  Cov(\cdot,\cdot) \) denotes the covariance. The score function with respect to \(\boldsymbol\theta\) is given by
\begin{equation}
    \frac{\partial l_{n}(\boldsymbol\theta)}{\partial \theta_{i}} = -\frac{1}{2}  Tr   \left( \left( \Gamma^{Z}_{n}(\boldsymbol\theta) \right)^{-1} \frac{\partial \Gamma^{Z}_{n}(\boldsymbol\theta)}{\partial \theta_{i}} \right) + \frac{1}{2}  \boldsymbol{Z}^{(n)^\top} \left( \Gamma^{Z}_{n}(\boldsymbol\theta) \right)^{-1} \frac{\partial \Gamma^{Z}_{n}(\boldsymbol\theta)}{\partial \theta_{i}} \left( \Gamma^{Z}_{n}(\boldsymbol\theta) \right)^{-1}  \boldsymbol{Z}^{(n)},  
\end{equation}
where \( Tr(\cdot)\) denotes the trace of a matrix. The Fisher information matrix (FIM) can be deduced from equation (\ref{3.3}). We simulate the spectral density and its derivatives using the method described in \citep{hariz2024fast}, then plug the FIM and score functions into equation (\ref{OSESTIMATOR}) to compute the One-Step estimator numerically.

 {For each set of parameters, specifically \((\phi(1),\phi(2),H)=(0.6,0.8,0.3)\) and \((\phi(1),\phi(2),H)=(0.6,0.7,0.6)\), we conduct \(M = 1000\) Monte Carlo simulations. The sample sizes considered are \(n = 100\), \(n = 1000\), and \(n = 2000\). }The number of Fourier frequencies for the initial estimations is set as \(m = [  n^{0.6}]\) and remains fixed throughout the simulations. 
 Without loss of generality, we assume \(T = 2\), and the spectral density of \(Z_{n}\) in this case is given by propositon 3.3.
 
   \begin{figure}[htbp]
\centering
\vspace{2\baselineskip}
\includegraphics[width=1\textwidth]{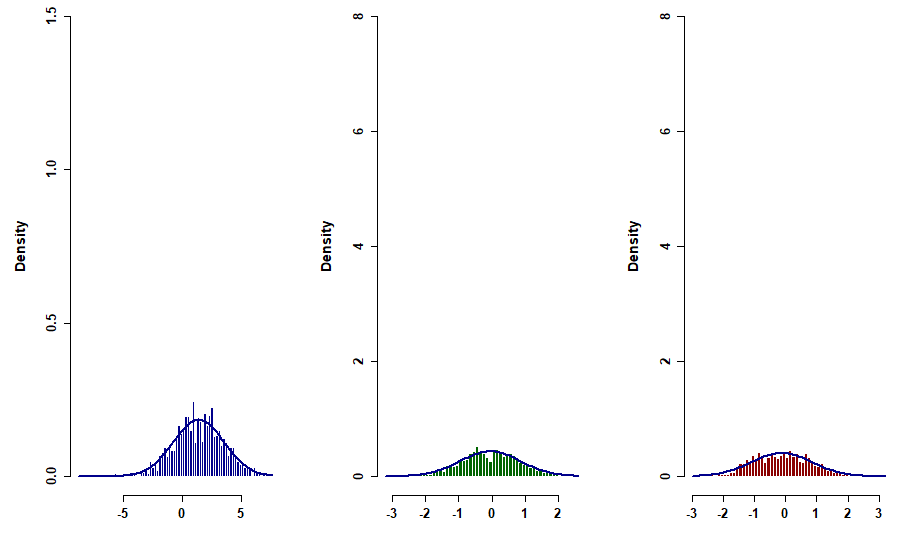}
\hspace{30pt}\bf{(a) IE of \boldmath$H$ = 0.3}
\hspace{70pt}
 \text{(b) IE of} $\boldsymbol{\phi}$(1) = 0.6
\hspace{60pt}
\text{(c) IE of} \boldmath$\mathbf{\phi}$(2) = 0.8\\
\includegraphics[width=1\textwidth]{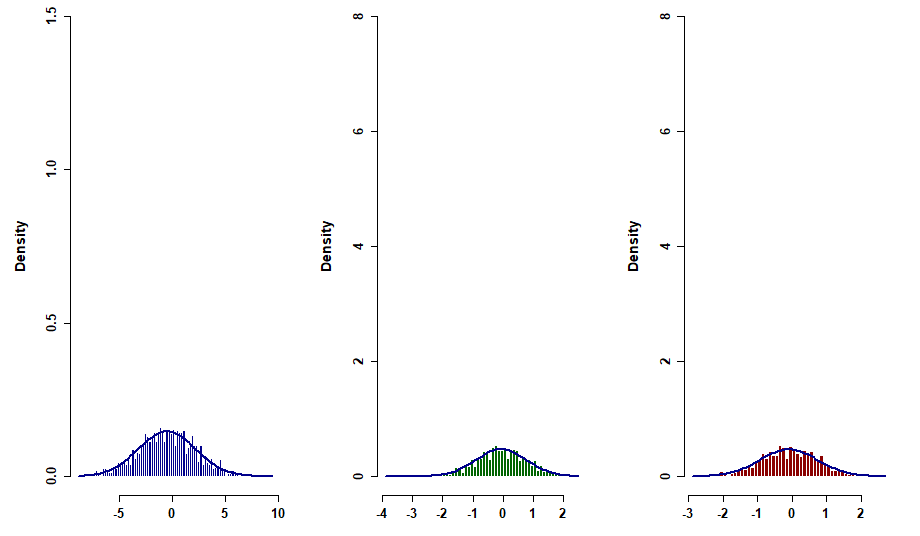}
\hspace{30pt} \bf{(d) OS of \boldmath$H$ = 0.3}
\hspace{70pt}
 \text{(e) OS of} $\boldsymbol{\phi}$(1) = 0.6
\hspace{60pt} 
\text{(f) OS of} \boldmath$\mathbf{\phi}$(2) = 0.8\\
\caption{{The simulation of initial estimator and One-Step estimator where $\theta=(0.6,0.8,0.3)$ for  $ n=100$.}}
    \label{683h}
\end{figure}

 \begin{figure}[htbp]

    \centering
     \vspace{2\baselineskip}
   \includegraphics[width=1\textwidth]{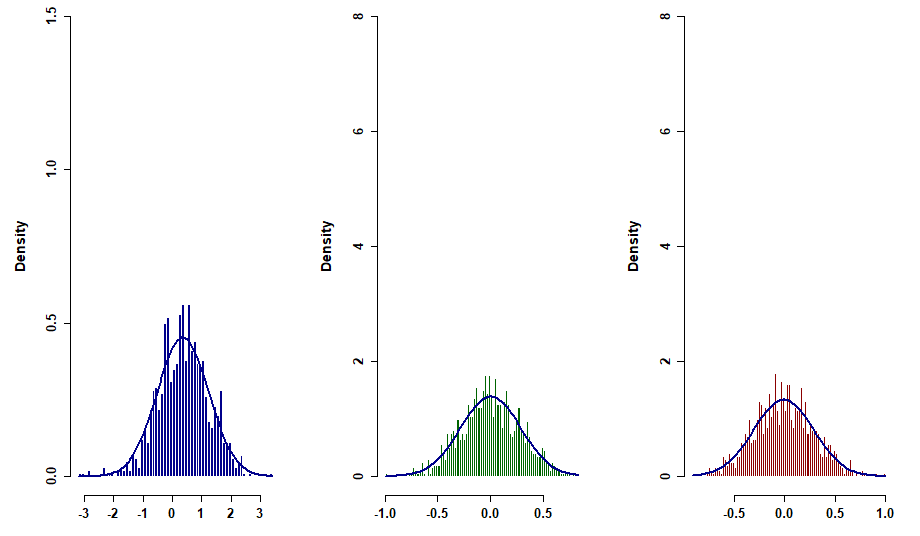}
 \hspace{30pt}\bf{(a) IE of \boldmath$H$ = 0.3}
\hspace{70pt} 
 \text{(b) IE of} $\boldsymbol{\phi}$(1) = 0.6
\hspace{60pt}
\text{(c) IE of} \boldmath$\mathbf{\phi}$(2) = 0.8\\
   \includegraphics[width=1\textwidth]{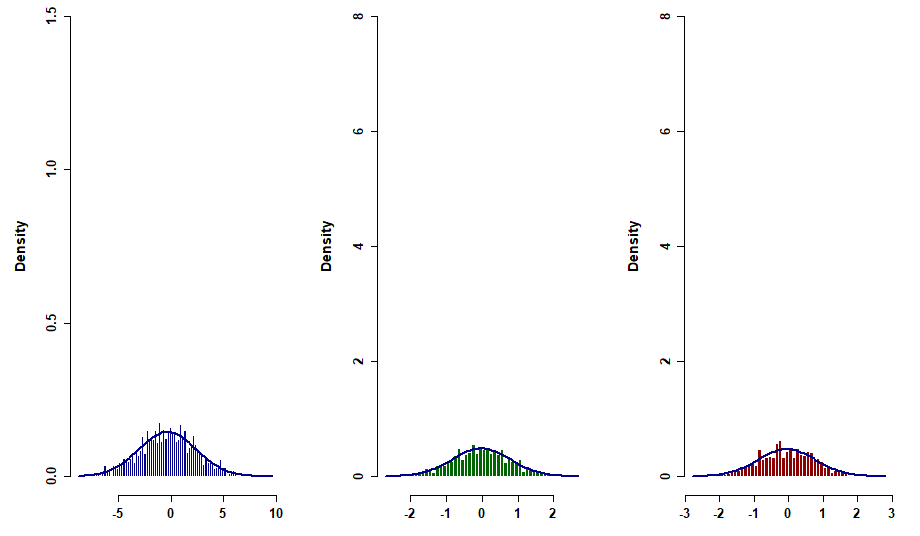}
\qquad  \bf{(d) OS of \boldmath$H$ = 0.3}
\hspace{70pt} 
 \text{(e) OS of} $\boldsymbol{\phi}$(1) = 0.6
\hspace{60pt}
\text{(f) OS of} \boldmath$\mathbf{\phi}$(2) = 0.8\\
   \caption{{The simulation of initial estimator and One-Step estimator where $\theta=(0.6,0.8,0.3)$ for $ n=1000$.}}
    \label{683t}
\end{figure}

 \begin{figure}[htbp]
    \centering
     \vspace{2\baselineskip}
   \includegraphics[width=1\textwidth]{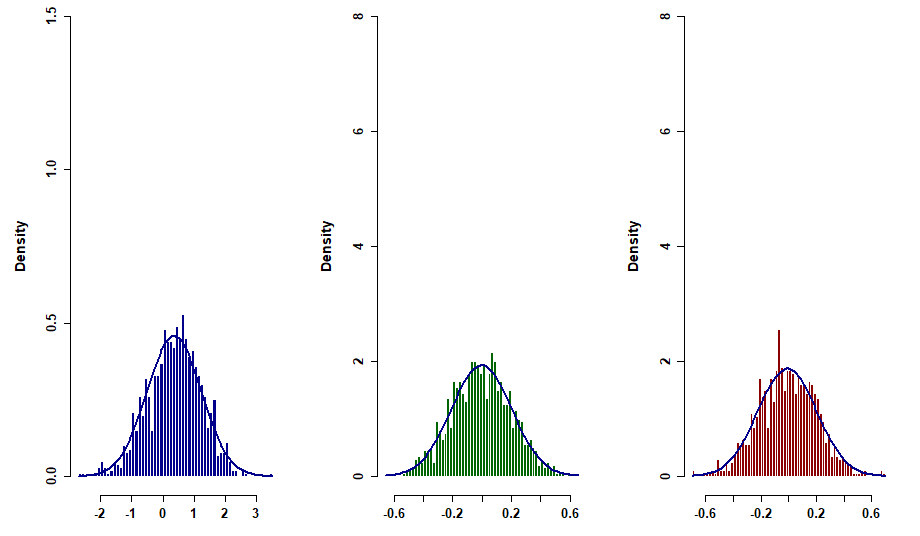}
 \hspace{30pt}\bf{(a) IE of \boldmath$H$ = 0.3}
\hspace{70pt} 
 \text{(b) IE of} $\boldsymbol{\phi}$(1) = 0.6
\hspace{60pt}
\text{(c) IE of} \boldmath$\mathbf{\phi}$(2) = 0.8\\
   \includegraphics[width=1\textwidth]{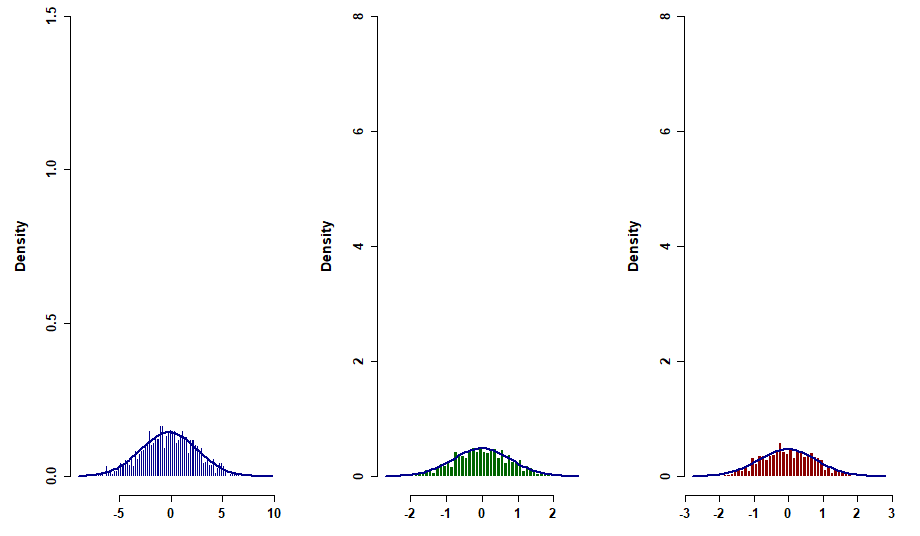}
\qquad  \bf{(d) OS of \boldmath$H$ = 0.3}
\hspace{70pt} 
 \text{(e) OS of} $\boldsymbol{\phi}$(1) = 0.6
\hspace{60pt}
\text{(f) OS of} \boldmath$\mathbf{\phi}$(2) = 0.8\\
   \caption{{The simulation of initial estimator and One-Step estimator where $\theta=(0.6,0.8,0.3)$ for $ n=2000$.}}
    \label{6832tt}
\end{figure}
 
\newpage
% 第一个表格
\begin{table}[htbp]
    \caption{The Bias and RMSE of Initial estimator and One-step estimator for $\boldsymbol\theta=(0.6,0.8,0.3)$ when $n = 100$}
    \centering
    \begin{tabular}{lcccc}
        \toprule
        \multicolumn{1}{c}{$n=100$} &  {B IE} &  {B OS} &  {RMSE IE} &  {RMSE OS} \\
        \midrule
        {$H$} &  {0.1725} &  {-0.0157} &  {0.2720} &  {0.0857} \\
        {$\phi(1)$} &  {-0.0072} &  {-0.0022} &  {0.1151} &  {0.0265} \\
        {$\phi(2)$} &  {-0.0162} &  {-0.0032} &  {0.1235} &  {0.0269} \\
        \bottomrule
    \end{tabular}
    \label{table1}
\end{table}

% 第二个表格
\begin{table}[htbp]
    \caption{The Bias and RMSE of Initial estimator and One-step estimator for $\boldsymbol\theta=(0.6,0.8,0.3)$ when $n = 1000$}
    \centering
    \begin{tabular}{lcccc}
        \toprule
        \multicolumn{1}{c}{$n=1000$} &  {B IE} &  {B OS} &  {RMSE IE} &  {RMSE OS} \\
        \midrule
        {$H$} &  {0.0454} &  {0.0110} &  {0.1109} &  {0.0869} \\
        {$\phi(1)$} &  {0.0003} &  {0.0007} &  {0.0362} &  {0.0260} \\
        {$\phi(2)$} &  {-0.0012} &  {0.0018} &  {0.0376} &  {0.0267} \\
        \bottomrule
    \end{tabular}
    \label{table2}
\end{table}

% 第三个表格
\begin{table}[htbp]
    \caption{The Bias and RMSE of Initial estimator and One-step estimator for $\boldsymbol\theta=(0.6,0.8,0.3)$ when $n = 2000$}
    \centering
    \begin{tabular}{lcccc}
        \toprule
        \multicolumn{1}{c}{$n=2000$} &  {B IE} &  {B OS} &  {RMSE IE} &  {RMSE OS} \\
        \midrule
        {$H$} &  {0.0454} &  {-0.0071} &  {0.1097} &  {0.0876} \\
        {$\phi(1)$} &  {-0.0004} &  {-0.0004} &  {0.0260} &  {0.0260} \\
        {$\phi(2)$} &  {-0.0012} &  {-0.0014} &  {0.0268} &  {0.0267} \\
        \bottomrule
    \end{tabular}
    \label{table3}
\end{table}
 
   \begin{figure}[htbp]
\centering
\vspace{2\baselineskip}
\includegraphics[width=1\textwidth]{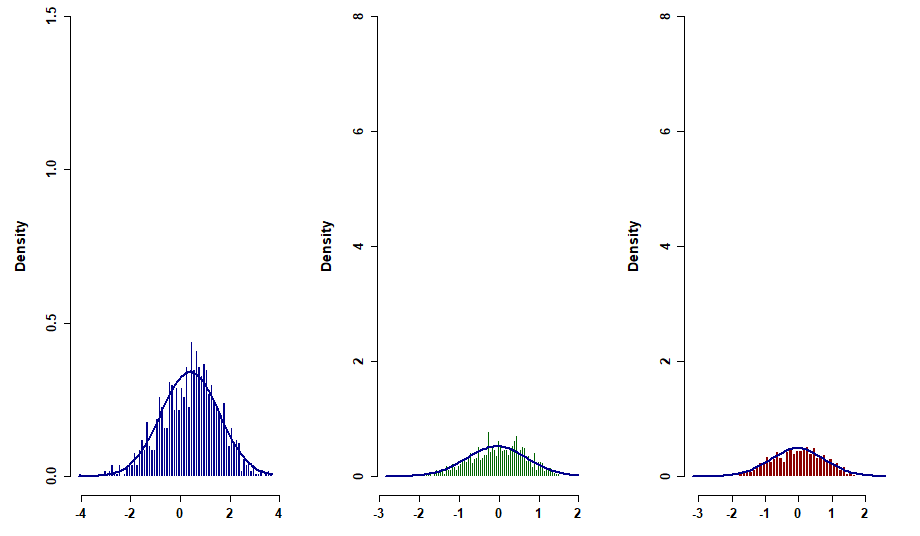}
\hspace{30pt}\bf{(a) IE of \boldmath$H$ = 0.6}
\hspace{70pt}
 \text{(b) IE of} $\boldsymbol{\phi}$(1) = 0.7
\hspace{60pt}
\text{(c) IE of} \boldmath$\mathbf{\phi}$(2) = 0.6\\
\includegraphics[width=1\textwidth]{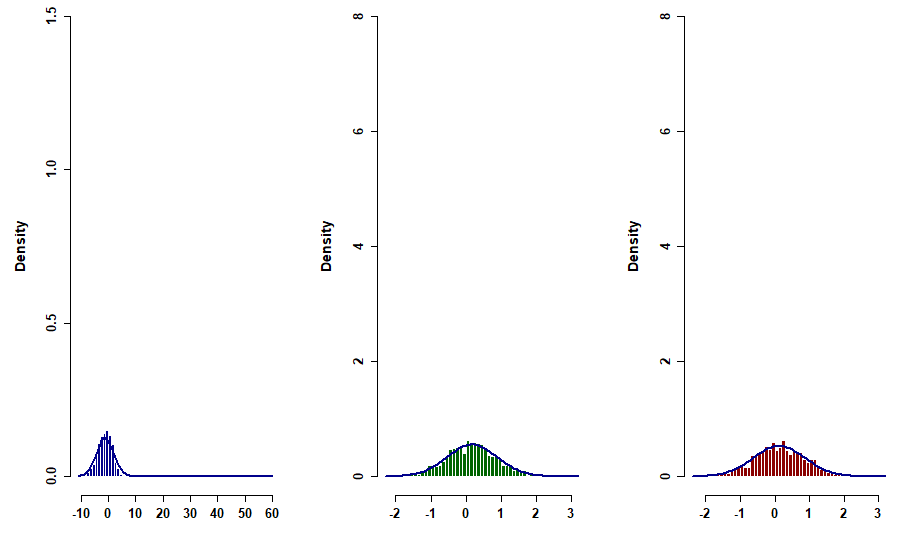}
\hspace{30pt} \bf{(d) OS of \boldmath$H$ = 0.6}
\hspace{70pt}
 \text{(e) OS of} $\boldsymbol{\phi}$(1) = 0.6
\hspace{60pt} 
\text{(f) OS of} \boldmath$\mathbf{\phi}$(2) = 0.7\\
\caption{{The simulation of initial estimator and One-Step estimator where $\theta=(0.6,0.7,0.6)$ for  $ n=100$.}}
    \label{776}
\end{figure}

 \begin{figure}[htbp]

    \centering
     \vspace{2\baselineskip}
   \includegraphics[width=1\textwidth]{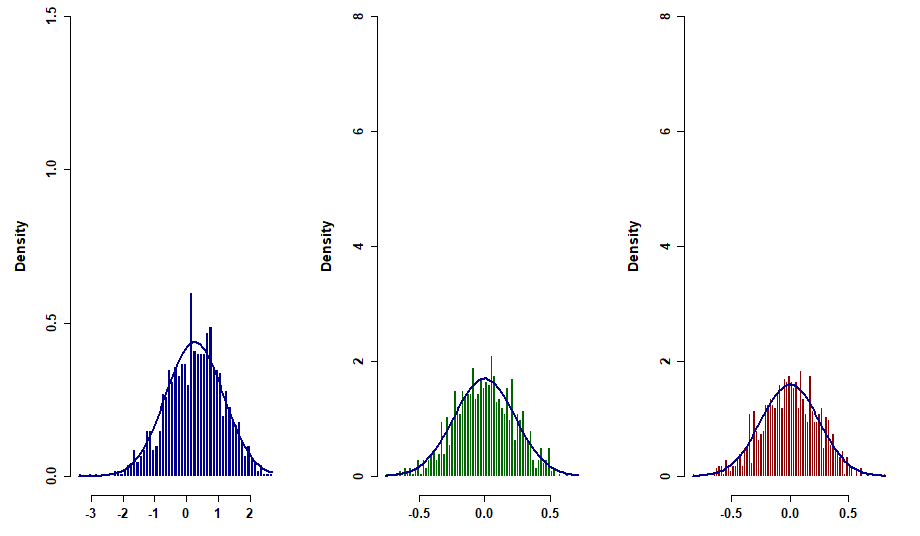}
 \hspace{30pt}\bf{(a) IE of \boldmath$H$ = 0.6}
\hspace{70pt} 
 \text{(b) IE of} $\boldsymbol{\phi}$(1) = 0.6
\hspace{60pt}
\text{(c) IE of} \boldmath$\mathbf{\phi}$(2) = 0.7\\
   \includegraphics[width=1\textwidth]{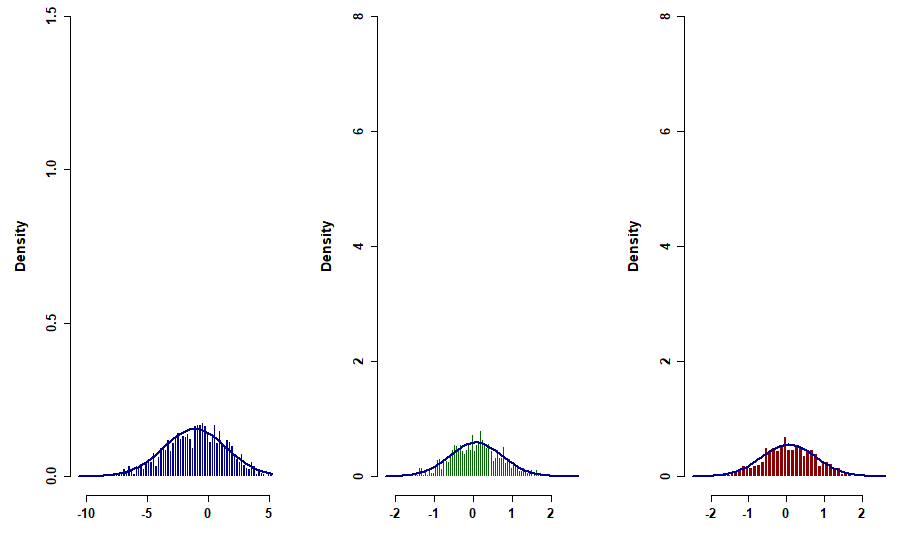}
\qquad  \bf{(d) OS of \boldmath$H$ = 0.6}
\hspace{70pt} 
 \text{(e) OS of} $\boldsymbol{\phi}$(1) = 0.6
\hspace{60pt}
\text{(f) OS of} \boldmath$\mathbf{\phi}$(2) = 0.7\\
   \caption{{The simulation of initial estimator and One-Step estimator where $\theta=(0.6,0.7,0.6)$ for $ n=1000$.}}
    \label{683t}
\end{figure}

 \begin{figure}[htbp]
    \centering
     \vspace{2\baselineskip}
   \includegraphics[width=1\textwidth]{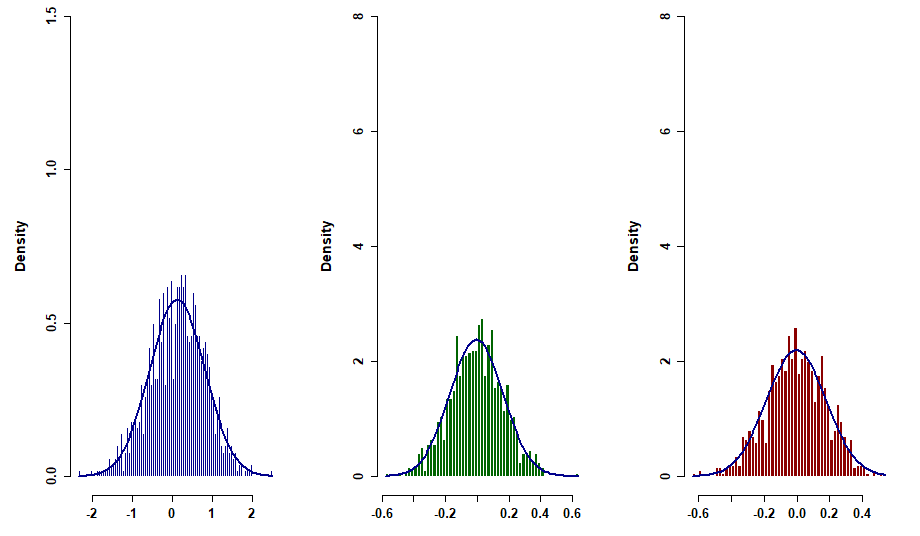}
 \hspace{30pt}\bf{(a) IE of \boldmath$H$ = 0.6}
\hspace{70pt} 
 \text{(b) IE of} $\boldsymbol{\phi}$(1) = 0.6
\hspace{60pt}
\text{(c) IE of} \boldmath$\mathbf{\phi}$(2) = 0.7\\
   \includegraphics[width=1\textwidth]{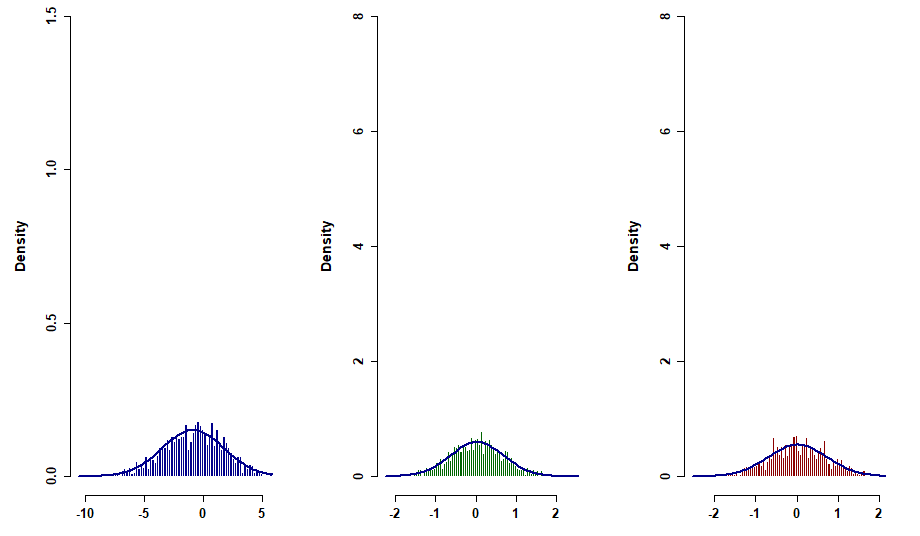}
\qquad  \bf{(d) OS of \boldmath$H$ = 0.6}
\hspace{70pt} 
 \text{(e) OS of} $\boldsymbol{\phi}$(1) = 0.6
\hspace{60pt}
\text{(f) OS of} \boldmath$\mathbf{\phi}$(2) = 0.7\\
   \caption{{The simulation of initial estimator and One-Step estimator where $\theta=(0.6,0.7,0.6)$ for $ n=2000$.}}
    \label{6832t}
\end{figure}
 
\newpage
% 第一个表格
\begin{table}[htbp]
    \caption{The Bias and RMSE of Initial estimator and One-step estimator for $\boldsymbol\theta=(0.6,0.7,0.6)$ when $n = 100$}
    \centering
    \begin{tabular}{lcccc}
        \toprule
        \multicolumn{1}{c}{$n=100$} &  {B IE} &  {B OS} &  {RMSE IE} &  {RMSE OS} \\
        \midrule
        {$H$} &  {0.0505} &  {-0.0413} &  {0.1472} &  {0.1011} \\
        {$\phi(1)$} &  {-0.0069} &  {0.0052} &  {0.1151} &  {0.0228} \\
        {$\phi(2)$} &  {-0.0050} &  {0.0038} &  {0.1235} &  {0.0238} \\
        \bottomrule
    \end{tabular}
    \label{table3}
\end{table}

% 第二个表格
\begin{table}[htbp]
    \caption{The Bias and RMSE of Initial estimator and One-step estimator for $\boldsymbol\theta=(0.6,0.7,0.6)$ when $n = 1000$}
    \centering
    \begin{tabular}{lcccc}
        \toprule
        \multicolumn{1}{c}{$n=1000$} &  {B IE} &  {B OS} &  {RMSE IE} &  {RMSE OS} \\
        \midrule
        {$H$} &  {0.0327} &  {-0.0352} &  {0.1148} &  {0.0821} \\
        {$\phi(1)$} &  {-0.0001} &  {0.0023} &  {0.0286} &  {0.0217} \\
        {$\phi(2)$} &  {0.0002} &  {0.0016} &  {0.0314} &  {0.0231} \\
        \bottomrule
    \end{tabular}
    \label{table4}
\end{table}

% 第三个表格
\begin{table}[htbp]
    \caption{The Bias and RMSE of Initial estimator and One-step estimator for $\boldsymbol\theta=(0.6,0.7,0.6)$ when $n = 2000$}
    \centering
    \begin{tabular}{lcccc}
        \toprule
        \multicolumn{1}{c}{$n=2000$} &  {B IE} &  {B OS} &  {RMSE IE} &  {RMSE OS} \\
        \midrule
        {$H$} &  {0.0158} &  {-0.0287} &  {0.0872} &  {0.0832} \\
        {$\phi(1)$} &  {-0.0006} &  {-0.0007} &  {0.0211} &  {0.0260} \\
        {$\phi(2)$} &  {-0.0005} &  {-0.0005} &  {0.0229} &  {0.0229} \\
        \bottomrule
    \end{tabular}
    \label{table5}
\end{table}

\par  {Figure 1 ,2 and 3 depict the frequency distribution of statistical errors for the initial estimatior and One-Step estimatior of the PFAR(1) model with parameters \(\phi(1) = 0.6\), \(\phi(2) = 0.8\) and \(H = 0.3\). Figure 4, 5 and 6 depict the frequency distribution of statistical errors for the initial estimatior and One-Step estimatior of the PFAR(1) model with parameters \(\phi(1) = 0.6\), \(\phi(2) = 0.7\) and \(H = 0.6\).} 
\par  {In all the tables, B stands for Bias, IE represents initial estimator, and OS denotes One-Step estimator. From the above tables, it can be seen that the OS estimator shows a significant improvement in the estimation of $H$. From these figures and the accompanying table, it is evident that the One-Step estimatior outperforms the initial estimation, with a particularly notable improvement in estimating the parameter $H$, at the same time, we found that as the sample size increases, the estimation becomes more efficient. According to \citep{hariz2024fast} and our simulations, the One-Step estimation also has a faster running speed.} 
\section{Time Series Modeling with PFAR(1) Models}
In this section, we will perform practical modeling and analysis to evaluate the application effectiveness of the PFAR model in real-world data. The data on the Standardized Precipitation Evapotranspiration Index (SPEI) in the Donner und Blitzen watershed, located in the Great Basin Desert, USA, from 1988 to 2020, selected in this study, are sourced from the public data\cite{dunham2024drought} of the United States Geological Survey. This dataset documents the monthly SPEI values from 1988 to 2020. To simplify the modeling process, we calculate the average of the data for each of the 12 months on a quarterly basis, thereby obtaining the quarterly SPEI values.  Given the distinct periodicity of SPEI, we opt to use the PFAR(1) model with \(T = 4\)  to simulate the aforementioned observations. Additionally, we compare the results with the simulation of the white noise driven periodic autoregressive (PAR) model and fractional autoregressive (FAR) model without periodic structure.

\par The periodic autoregressive model driven by white noise is as follow: 
 \begin{equation}
    X_{4n + u}=\alpha(4n + u)X_{4n+u-1}+\epsilon_{4n + u},\quad u = 1,2,3,4,\quad n\in\mathbb{Z}^{+},
\end{equation} 
and FAR(1) model is
\begin{equation}
    X_{u}=\beta X_{u-1}+\epsilon_{u}^{H}, \quad u \in\mathbb{Z}^{+},
\end{equation} 
 where $\epsilon_{4n + u}$ is white noise and $\epsilon_{ u}^{H}$ is fractional Gaussian noise.  $\alpha(1)$, $\alpha(2)$, $\alpha(3)$, $\alpha(4)$ are the model coefficients satisfying $\alpha(4n + u) = \alpha(u)$ for periodicity $T = 4$; and $\beta$ is the model coefficient of the FAR model. 
\par {We utilized the data  from 1988 to 2020 to derive the parameter estimations of the two models. Subsequently, we computed their RMSE  and MAE against the real data. The results are shown in the following table and figure.}

 \begin{table}[!htbp] %voc table result
    \centering
    \caption{{Fitting results of different models  (PAR and PFAR)}}
    \begin{tabular}{*{7}{c}}
        \toprule
        Model &  PAR & PFAR & FAR   \\
        \midrule
        Parameters    & ($\alpha(1)$, $\alpha(2)$, $\alpha(3)$, $\alpha(4)$) &($\phi(1)$, $\phi(2)$, $\phi(3)$, $\phi(4)$, $H$)  &$(\beta, H)$\\
        Values& (0.65, 0.71, 0.97, 0.90) & (0.46, 0.34, 0.53, 0.61, 0.95) & (0.48, 0.95)  \\
        RMSE  & 1.8136 & 1.5388 &1.8773 \\
        MAE &  1.4478 & 1.1877 &1.5252\\
        \bottomrule
    \end{tabular}
    \label{table4}
\end{table} 

\begin{figure}[h]
    \centering
    \includegraphics[width = 0.8\textwidth]{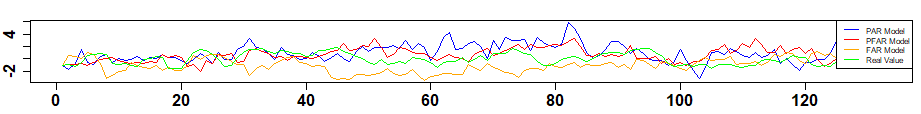}
    \caption{{Model fitting plot  (PAR, PFAR and FAR) }}
    \label{model.png}
\end{figure}
In the above figure, the values represent the parameters fitted by the PFAR and FAR models. We can observe that both models exhibit favorable fitting performance, indicating that the SPEI series has a periodicity of \(T = 4\) and long-range dependence characterized by a Hurst index of 0.95. As shown in Figures \ref{model.png}, the SPEI demonstrates stronger periodicity than long memory.  
 \section{Conclusions}
This paper presents a study on parameter estimation for a periodic fractional autoregressive process (PFAR) driven by fractional Gaussian noise. The PFAR model, a specialized varying coefficient fractional autoregressive model with periodic coefficients, is introduced to capture long memory and periodicity in time series.

The authors primarily address three key problems. First, we utilize the Generalized Least Squares Estimation (GLSE) to obtain the initial estimators of the periodic parameters for the PFAR(1) model, and then apply the GPH method. This approach yields a single estimate of $H$, rather than multiple estimates. Second, we derive asymptotic results for a specific initial estimator. Finally, an One-Step estimator is proposed to enhance the asymptotic efficiency. Theoretical analysis demonstrates that both estimators are consistent and asymptotically normal.

Monte Carlo simulations are conducted with different parameter values $(\phi(1), \phi(2), H)$ and sample sizes ($n = 100$, $n = 1000$, $n = 2000$). The simulation results show that the One-Step estimator performs better than the initial estimator, especially in estimating the Hurst index $H$. As the sample size increases, the estimation accuracy improves. Overall, this research provides effective methods for parameter estimation in PFAR models, which can be useful in various fields such as finance, meteorology, and engineering where long memory and periodicity are common features in time series data.\\
%Bibliography
\bibliographystyle{unsrt}  
\bibliography{references} 

\begin{thebibliography}{10}

\bibitem{hosking1984modeling}
Jonathan~RM Hosking.
\newblock Modeling persistence in hydrological time series using fractional differencing.
\newblock {\em Water resources research}, 20(12):1898--1908, 1984.

\bibitem{hannan1955test}
EJ~Hannan.
\newblock A test for singularities in sydney rainfall.
\newblock {\em Australian Journal of Physics, vol. 8, p. 289}, 8:289, 1955.
\newblock \doi {10.1071/ph550289}.

\bibitem{fernandez1986periodic}
Bonifacio Fernandez and Jose~D Salas.
\newblock Periodic gamma autoregressive processes for operational hydrology.
\newblock {\em Water Resources Research}, 22(10):1385--1396, 1986.
\newblock \doi {10.1029/WR022i010p01385}.

\bibitem{franses2004periodic}
Philip~Hans Franses and Richard Paap.
\newblock {\em Periodic time series models}.
\newblock OUP Oxford, 2004.
\newblock \doi {10.1093/019924202X.001.0001}.

\bibitem{comte1998long}
Fabienne Comte and Eric Renault.
\newblock Long-memory in continuous-time stochastic volatility models.
\newblock {\em Mathematical Finance}, 8(4):291--323, 1998.
\newblock \doi {10.1111/1467-9965.00057}.

\bibitem{gloter2004stochastic}
Arnaud Gloter and Marc Hoffmann.
\newblock Stochastic volatility and fractional brownian motion.
\newblock {\em Stochastic processes and their applications}, 113(1):143--172, 2004.
\newblock \doi {10.1016/j.spa.2004.03.008}.

\bibitem{corsi2021roughness}
Fulvio Corsi, Francesco Audrino, and Roberto Ren{'o}.
\newblock Roughness of volatility: Evidence from option prices.
\newblock {\em Journal of Financial Econometrics}, 19(4):1019--1056, 2021.
\newblock \doi {10.1080/24725854.2018.144429}.

\bibitem{richard2023discrete}
Alexandre Richard, Xiaolu Tan, and Fan Yang.
\newblock On the discrete-time simulation of the rough heston model.
\newblock {\em SIAM Journal on Financial Mathematics}, 14(1):223--249, 2023.
\newblock \doi {10.1137/21M1443807}.

\bibitem{jones1967time}
Richard~H Jones and William~M Brelsford.
\newblock Time series with periodic structure.
\newblock {\em Biometrika}, 54(3-4):403--408, 1967.
\newblock \doi {10.1093/biomet/54.3-4.403}.

\bibitem{vecchia1985maximum}
AV~Vecchia.
\newblock Maximum likelihood estimation for periodic autoregressive moving average models.
\newblock {\em Technometrics}, 27(4):375--384, 1985.
\newblock \doi {10.1080/00401706.1985.10488076}.

\bibitem{brouste2014asymptotic}
Alexandre Brouste, Chunhao Cai, and Marina Kleptsyna.
\newblock Asymptotic properties of the mle for the autoregressive process coefficients under stationary gaussian noise.
\newblock {\em Mathematical Methods of Statistics}, 23:103--115, 2014.
\newblock \doi {10.3103/S1066530714020021}.

\bibitem{soltane2024asymptotic}
Marius Soltane.
\newblock Asymptotic efficiency in autoregressive processes driven by stationary gaussian noise.
\newblock {\em Stochastic Models}, 40(1):70--96, 2024.
\newblock \doi {10.1080/15326349.2023.2202227}.

\bibitem{esstafa2019long}
Youssef Esstafa.
\newblock {\em Long-memory time series models with dependent innovations}.
\newblock PhD thesis, Universit{\'e} Bourgogne Franche-Comt{\'e}, 2019.

\bibitem{hariz2024fast}
Samir~Ben Hariz, Alexandre Brouste, Chunhao Cai, and Marius Soltane.
\newblock Fast and asymptotically-efficient estimation in an autoregressive process with fractional type noise.
\newblock {\em Journal of Statistical Planning and Inference}, 232:106148, 2024.
\newblock \doi {10.1016/j.jspi.2024.106148}.

\bibitem{hurst1951long}
Harold~Edwin Hurst.
\newblock Long-term storage capacity of reservoirs.
\newblock {\em Transactions of the American Society of Civil Engineers}, 116(1):770--799, 1951.
\newblock \doi {10.1061/TACEAT.0006518}.

\bibitem{geweke1983estimation}
John Geweke and Susan Porter-Hudak.
\newblock The estimation and application of long memory time series models.
\newblock {\em Journal of Time Series Analysis}, 4(4):221--238, 1983.
\newblock \doi {10.1111/j.1467-9892.1983.tb00371.x}.

\bibitem{robinson1995log}
Peter~M Robinson.
\newblock Log-periodogram regression of time series with long range dependence.
\newblock {\em The annals of Statistics}, pages 1048--1072, 1995.
\newblock \doi {10.1214/aos/1176324636}.

\bibitem{le1956asymptotic}
Lucien Le~Cam.
\newblock On the asymptotic theory of estimation and testing hypotheses.
\newblock In {\em Proceedings of the Third Berkeley Symposium on Mathematical Statistics and Probability, Volume 1: Contributions to the Theory of Statistics}, volume~3, pages 129--157. University of California Press, 1956.
\newblock \doi {10.1525/9780520313880-014}.

\bibitem{kutoyants2016multi}
Yu~A Kutoyants and Anastasia Motrunich.
\newblock On multi-step mle-process for markov sequences.
\newblock {\em Metrika}, 79:705--724, 2016.
\newblock \doi {10.1007/s00184-015-0574-4}.

\bibitem{gloter2021adaptive}
Arnaud Gloter and Nakahiro Yoshida.
\newblock Adaptive estimation for degenerate diffusion processes.
\newblock 2021.
\newblock \doi {10.1214/20-EJS1777}.

\bibitem{marinucci1998semiparametric}
Domenico Marinucci and Peter~M Robinson.
\newblock Semiparametric frequency domain analysis of fractional cointegration.
\newblock 1998.
\newblock \doi {10.1093/oso/9780199257294.003.0015}.

\bibitem{hurvich1998mean}
Clifford~M Hurvich, Rohit Deo, and Julia Brodsky.
\newblock The mean squared error of geweke and porter-hudak's estimator of the memory parameter of a long-memory time series.
\newblock {\em Journal of Time Series Analysis}, 19(1):19--46, 1998.
\newblock \doi {10.1111/1467-9892.00075}.

\bibitem{cohen:hal-00638121}
Serge Cohen, Fabrice Gamboa, C{\'e}line Lacaux, and Jean-Michel Loubes.
\newblock {LAN property for some fractional type Brownian motion}.
\newblock {\em {ALEA : Latin American Journal of Probability and Mathematical Statistics}}, 10(1):91--106, 2013.

\bibitem{dunham2024drought}
J.B. Dunham, C.E. Torgersen, J.D. Groom, and P.K. Haggerty.
\newblock Drought-related responses of stream flow, climate, and vegetation productivity from the donner und blitzen watershed, great basin desert, usa (1988-2020).
\newblock U.S. Geological Survey data release, \url{https://doi.org/10.5066/P13VBQMF}, 2024.

\bibitem{lieberman2012asymptotic}
Offer Lieberman, Roy Rosemarin, and Judith Rousseau.
\newblock Asymptotic theory for maximum likelihood estimation of the memory parameter in stationary gaussian processes.
\newblock {\em Econometric Theory}, 28(2):457--470, 2012.
\newblock \doi {10.1017/S0266466611000399}.

\end{thebibliography}
\section*{Appendices}
The proofs in Appendices are divided into four parts. The first part deals with the discussion on the stationarity of the sequence and the proof of the spectral density of the related sequence. The second part is the proof of the results related to the initial estimator. The third part is the proof of the results related to the One-Step estimator. The fourth part is the simulation of the estimator in Remark 3.
\section*{Appendix 1. Proof of Stationarity and the Spectral Density}
\subsection*{A.1.Proof of Theorem \ref{Thm:sationaryofX}}
Because the coefficients of the subsequence are absolutely summable, it can be expressed as an infinite sum of fractional Gaussian noise.
        \begin{equation}
        X_{2n+2} = \sum^{\infty}_{i=0,2,4,...}\epsilon^{H}_{2n+2-i}[\phi(2)\phi(1)]^{\frac{i}{2}}+\sum^{\infty}_{i=1,3,5,...}\epsilon^{H}_{2n+2-i}[\phi(2)\phi(1)]^{\frac{i-1}{2}}\phi(2),
        \end{equation}
        \begin{equation}
        X_{2n+1} = \sum^{\infty}_{i=0,2,4,...}\epsilon^{H}_{2n+1-i}[\phi(2)\phi(1)]^{\frac{i}{2}}+\sum^{\infty}_{i=1,3,5,...}\epsilon^{H}_{2n+1-i}[\phi(2)\phi(1)]^{\frac{i-1}{2}}\phi(1),
        \end{equation}
        \par To verify that $(X_{2n+1})_{n \in \mathbb{N}}$ are wide stationary, the following three conditions must be satisfied.\\
\textbf{(1) Constant Mean:}$\mathbb{E}(X_{2n+u}) = \mu \in \mathbb{R}$ is constant for all $n \in \mathbb{N}$ and $u = 1, 2$.\\
 Without loss of generality, we assume $\mathbb{E}(\epsilon^{H}_{n}) = 0$. For any time series \((X_{nT+u})_{n \in \mathbb{N}}\) under the monotone convergence theorem and the Cauchy-Schwarz inequality, we obtain
  \begin{eqnarray}
  	\mathbb{E}|X_{2n+u}| 
  	&\leq&
   \mathbb{E}  \left\{ \sum^{\infty}_{i=0,2,4,...}|\epsilon^{H}_{2n+1-i}[\phi(2)\phi(1)]^{\frac{i}{2}}|+\sum^{\infty}_{i=1,3,5,...}|\epsilon^{H}_{2n+1-i}[\phi(2)\phi(1)]^{\frac{i-1}{2}}\phi(u)|\right\}  \nonumber  \\
  	&\leq &
  	 \sum^{\infty}_{i=0,2,4,...}\mathbb{E}|\epsilon^{H}_{2n+1-i}[\phi(2)\phi(1)]^{\frac{i}{2}}|+\sum^{\infty}_{i=1,3,5,...}\mathbb{E}|\epsilon^{H}_{2n+1-i}[\phi(2)\phi(1)]^{\frac{i-1}{2}}\phi(u)| \nonumber\\
  	&\leq &
  	G \sum^{\infty}_{i=0,2,4,...}|[\phi(2)\phi(1)]^{\frac{i}{2}}|+|\phi(u)|\sum^{\infty}_{i=1,3,5,...}| [\phi(2)\phi(1)]^{\frac{i-1}{2}}|, \label{5} 
  \end{eqnarray}
  where $G$ is a finite constant. We know that $\phi^j(u)$ is absolutely summable, i.e., $\sum_{j=0}^{\infty} |\phi^j(u)| < \infty$. Thus, $\mathbb{E}|X_{2n+u}| < \infty$ as shown in equation (\ref{5}). By the monotone convergence theorem, $\sum_{j=0}^{\infty} \phi^j(u)$ is absolutely convergent almost surely.

\par Considering that
\begin{equation}
    \left|\sum^{\infty}_{i=0,2,4,...}\epsilon^{H}_{2n+1-i}[\phi(2)\phi(1)]^{\frac{i}{2}} \right| \leq \sum^{\infty}_{i=0,2,4,...}\left|\epsilon^{H}_{2n+1-i}[\phi(2)\phi(1)]^{\frac{i}{2}} \right|
\end{equation}

\begin{equation}
  \left|\sum^{\infty}_{i=1,3,5,...}\epsilon^{H}_{2n+1-i}[\phi(2)\phi(1)]^{\frac{i-1}{2}} \right| \leq \sum^{\infty}_{i=1,3,5,...}\left| \epsilon^{H}_{2n+1-i}[\phi(2)\phi(1)]^{\frac{i-1}{2}} \right|
\end{equation}
by the dominated convergence theorem, 
\begin{equation}
    \mathbb{E}(X_{2T+u}) = \lim_{k \to \infty} \mathbb{E} \left(\sum^{2k}_{i=0,2,4,...}\epsilon^{H}_{2n+1-i}[\phi(2)\phi(1)]^{\frac{i}{2}}+\sum^{2k+1}_{i=1,3,5,...}\epsilon^{H}_{2n+1-i}[\phi(2)\phi(1)]^{\frac{i-1}{2}}\phi(u) \right) = 0. \label{eq:MeanvalueofX}
\end{equation}
  \textbf{(2) Finite Second Moment:}$\mathbb{E}(X_{2n+u}^2) < \infty$ for $ u=1,2$.
\par From equation (\ref{eq:MeanvalueofX}), we derive
\begin{align}
    \mathbb{E}(X_{nT+u})^2 
    &= \mathbb{E} \left( \sum^{\infty}_{i=0,2,4,...}\epsilon^{H}_{2n+1-i}[\phi(2)\phi(1)]^{\frac{i}{2}}+\sum^{\infty}_{i=1,3,5,...}\epsilon^{H}_{2n+1-i}[\phi(2)\phi(1)]^{\frac{i-1}{2}}\phi(u),
        \right)^2 \nonumber \\
    &= \mathbb{E} \left( \sum_{s_{1}=0,2,4...}^{\infty} \sum_{k_{1}=0,2,4...}^{\infty} \phi^{\frac{s_{1}}{2}}(u) \phi^{\frac{k_{1}}{2}}(u) \epsilon^{H}_{2n+1-s_{1}} \epsilon^{H}_{2n+1-k_{1}} \right) \nonumber \\
     &\quad+ \phi^{2}(u)\mathbb{E} \left( \sum_{s_{2}=1,3,5...}^{\infty} \sum_{k_{2}=1,3,5...}^{\infty} \phi^{\frac{s_{2}-1}{2}}(u) \phi^{\frac{k_{2}-1}{2}}(u) \epsilon^{H}_{2n+1-s_{1}} \epsilon^{H}_{2n+1-k_{1}} \right) \nonumber \\
     &\quad+\phi(u)\mathbb{E} \left( \sum_{s_{3}=0,2,4...}^{\infty} \sum_{k_{3}=1,3,5...}^{\infty} \phi^{\frac{s_{3}}{2}}(u) \phi^{\frac{k_{3}-1}{2}}(u) \epsilon^{H}_{2n+1-s_{3}} \epsilon^{H}_{2n+1-k_{3}} \right) \nonumber \\
\end{align}
 By applying the conclusion above, we obtain
 \begin{align}
    \mathbb{E}|X_{nT+u}|^2 
    &= \sum_{s_{1}=0,2,4...}^{\infty} \sum_{k_{1}=0,2,4...}^{\infty} \left|\phi(1) \phi(2)\right|^{\frac{s_{1}+k_{1}}{2}}\mathbb{E} \left| \epsilon^{H}_{2n+1-s_{1}} \epsilon^{H}_{2n+1-k_{1}} \right| \nonumber \\
     &\quad + \phi^{2}(u) \sum_{s_{2}=1,3,5...}^{\infty} \sum_{k_{2}=1,3,5...}^{\infty} \left|\phi(1) \phi(2)\right|^{\frac{s_{2}+k_{2}}{2}-1}\mathbb{E} \left|\epsilon^{H}_{2n+1-s_{2}} \epsilon^{H}_{2n+1-k_{2}} \right|   \nonumber \\
     &\quad+\phi(u) \sum_{s_{3}=0,2,4...}^{\infty} \sum_{k_{3}=1,3,5...}^{\infty} \left| \phi (1) \phi (2) \right|^{\frac{s_{3}+k_{3}-1}{2}}\mathbb{E}\left| \epsilon^{H}_{2n+1-s_{3}} \epsilon^{H}_{2n+1-k_{3}} \right|  \label{exofxn2nd} \\
\end{align}
and the covariance of \( \epsilon^{H}_{2n+1-s} \) and \( \epsilon^{H}_{2n+1-k} \) is
\begin{equation}
    Cov(\epsilon^{H}_{2n+1-s}, \epsilon^{H}_{2n+1-k}) = {\frac{1}{2}}\left( |s-k+1|^{2H} - 2 |s-k|^{2H} + |s-k-1|^{2H} \right).
\end{equation}
 Since \( \phi^j(u) \) is absolutely summable, it is also square summable. Additionally, as \( s - k \to \infty \), \( Cov(\epsilon^{H}_{2n+1-s}, \epsilon^{H}_{2n+1-k})  \to 0 \), implying that there exists a constant \( M \) such that \( \mathbb{E} \left| \epsilon^{H}_{2n+1-s} \epsilon^{H}_{2n+1-k} \right| \leq M \). Based on the above discussion and equation (\ref{exofxn2nd}), we have established that \( \mathbb{E}(X_{nT+u})^2 < \infty \).\\
\textbf{(3) Covariance Depends Only on Lag:}
  \[
        \mathbb{E}\left[(X_{kT+u} - \mu)(X_{sT+u} - \mu)\right] = \gamma((k - s)T), \quad \forall k, s \in \mathbb{N}.
    \] 
Let us denote 
        $$\sum^{\infty}_{i=0,2,4,...}\epsilon^{H}_{2n+2-i}[\phi(2)\phi(1)]^{\frac{i}{2}} = A_{2n+2},$$
        $$\sum^{\infty}_{i=1,3,5,...}\epsilon^{H}_{2n+2-i}[\phi(2)\phi(1)]^{\frac{i-1}{2}}\phi(2) = B_{2n+2},$$
         Then, the covariance of $(X_{2n+2})_{n \in \mathbb{N}}$ is
           \begin{align}
        	\gamma^{\phi(2)}_{k}
        	&=  
        	 Cov(X_{2n},X_{2(n-k)})  \nonumber \\
        	&  = 
        	 Cov(A_{2n}+B_{2n},A_{2(n-k)}+B_{2(n-k)})\nonumber \\
        	& = 
        	 Cov(A_{2n},A_{2(n-k)})+ Cov(B_{2n},B_{2(n-k)})\nonumber \\
          & \quad + Cov(A_{2n},B_{2(n-k)}) + Cov(B_{2n},A_{2(n-k)}),
        \end{align}
         Next, we analyze $Cov(A_{2n},A_{2(n-k)})$ and the cross term $Cov(A_{2n},B_{2(n-k)})$. The properties of the other two terms on the right side can be proved similarly.
         \begin{align}
        	Cov(A_{2n},A_{2(n-k)})
        	&=  
            Cov\left(\sum^{\infty}_{i=0,2,4,...}\epsilon^{H}_{2n-i}[\phi(2)\phi(1)]^{\frac{i}{2}},\sum^{\infty}_{i=0,2,4,...}\epsilon^{H}_{2(n-k)-i}[\phi(2)\phi(1)]^{\frac{i}{2}}  \right) \nonumber \\
        	& \quad = 
        	\sum^{\infty}_{j=0,2,4,...}\sum^{\infty}_{i = 0,2,4,...}[\phi(2)\phi(1)]^{\frac{i+j}{2}}\rho_{\epsilon^{H}}(2k-i+j),
        \end{align}
           Here, $\rho_{\epsilon^{H}}(\cdot)$is the autocovariance function of fractional Gaussian noise. Observing the above expression, we find that $Cov(A_{2n},A_{2(n-k)})$ is independent of the value of $n$.
          \begin{align}
         	Cov(A_{2n},B_{2(n-k)})
         	&=  
         	Cov\left(\sum^{\infty}_{i=0,2,4,...}\epsilon^{H}_{2n-i}[\phi(2)\phi(1)]^{\frac{i}{2}},\sum^{\infty}_{i=1,3,5,...}\epsilon^{H}_{2(n-k)-i}[\phi(2)\phi(1)]^{\frac{i-1}{2}}\phi(2) \right) \nonumber \\
         	& \quad = 
         	\sum^{\infty}_{j=1,3,5,...}\sum^{\infty}_{i = 0,2,4,...}\phi(2)[\phi(2)\phi(1)]^{\frac{i+j-1}{2}}\rho(2k-i+j),
         \end{align}
          \par Similarly, it can be shown that the cross term $Cov(A_{2n},B_{2(n-k)})$ is also independent of the value of $n$. The stationarity of the sequence $(X_{2n+1})_{n \in \mathbb{N}}$ can be proved in a similar manner. Therefore, we conclude that the covariance of the subsequences $(X_{2n+1})_{n \in \mathbb{N}}$ and $(X_{2n+2})_{n \in \mathbb{N}}$ depends only on the time interval, indicating that they are stationary processes.
\subsection*{A.2.Proof of Proposition \ref{prop:sdofx}}
We can express $X_{nT+u}$ as the infinite sum of fractional Gaussian noise:
        \begin{align}
            X_{nT+u}
        	&=  
        	\epsilon^{H}_{nT+u}+\phi(u)X_{nT+u-1} \nonumber \\
        	&\quad= 
        	\epsilon^{H}_{nT+u}+\phi(u)\epsilon^{H}_{nT+u-1}+ 
        \cdots+\phi(u)\phi(u-1)\cdots\phi(u-T)X_{(n-1)T+u-1}  \nonumber \\
        	&\quad= 
        	\sum_{i = 0,T,2T,...}C^{\frac{i}{T}}\epsilon^{H}_{nT-i+u}+\sum_{i = 1,T+1,2T+1,...}\phi(u)C^{\frac{i-1}{T}}\epsilon^{H}_{nT-i+u}   \nonumber \\
        	&\qquad+ 
        	\cdots+\sum_{i =T-1,2T-1,3T-1,...}\phi(u)\phi(u-1)\cdots\phi(u-T+2)C^{\frac{i-T+1}{T}}\epsilon^{H}_{nT-i+u} \nonumber \\
        	& \quad = 
        	\sum^{\infty}_{k=0} \tilde{h}_{k}\epsilon^{H}_{nT+u-k},   
         \end{align}
         where $C = \phi(u)\phi(u-1)\cdots\phi(u+T-1)$, $\tilde{h}_{k}$ is 
     $$\tilde{h}_{k} = \left\{
      \begin{aligned}
      &C^{\frac{k}{T}}, && \text{when} k=0,T,2T,...\\
      &\phi(u)C^{\frac{k-1}{T}}, && \text{when} k= 1,1+T,1+2T,..\\
     & \qquad \vdots & \\
    &\phi(u)\phi(u-1)\cdots\phi(u-T+2)C^{\frac{k-T+1}{T}}, && \text{when} k=T-1,2T-1,3T-1,...\\	
    \end{aligned}
   \right.$$
     Thus, we can get a Polynomial of $e^{-i\lambda}$ satisfies 
\begin{align}
	H(e^{-i\lambda}) 
 &= \sum^{\infty}_{j=kT}C^{\frac{j}{T}}e^{-ij\lambda}  + \sum^{\infty}_{j=1+kT}\phi(u)C^{\frac{j-1}{T}} e^{-ij\lambda} \nonumber \\
	&\quad +\cdots+ \sum^{\infty}_{j=(k+1)T-1}\phi(u)\phi(u-1)\cdots\phi(u-T+2)C^{\frac{j-T+1}{T}}e^{-ij\lambda} \quad k = 0,1,2,\ldots
\end{align}
     Let $q = Ce^{-iT\lambda}$, and $\Psi_{j}(\lambda)$ takes the form of
    \begin{equation}
    \Psi_{j}(\lambda) = \frac{\phi(u)\phi(u-1)\cdots\phi(u-j)}{1-q}e^{-ij\lambda},
    \end{equation}
    Thus, $H(e^{-i\lambda})$ has the following expression
    \begin{equation}
    H(e^{-i\lambda})= \sum^{T-2}_{j=0}\Psi_{j}(\lambda)+\frac{1}{1-q},
    \end{equation}
 Finally, we obtain 
   \begin{equation}
   \overline{f}_{H,\phi(u)}(\lambda) = \left| \sum^{T-2}_{j=0}\Psi_{j}(\lambda)+\frac{1}{1-Ce^{-iT\lambda}}\right|^{2}	{f_{\epsilon_{n}^{H}}(\lambda)}.
    \end{equation} 
\subsection*{A.3.Proof of Proposition \ref{prop:sdofz}}
    According to the definition of $Z_{n}$, we can obtain
    \begin{align}
    Z_{n}
   	&=  
   	X_{2n}+X_{2n+1}\nonumber \\
   	&\quad= 
   	\sum_{i=0,2,4,...}\epsilon^{H}_{2n-i}[\phi(1)\phi(2)]^{\frac{i}{2}}+\sum_{i=1,3,5,...}\epsilon^{H}_{2n-i}[\phi(1)\phi(2)]^{\frac{i-1}{2}}\phi(2)  \nonumber \\
   	&\qquad+ 
   	\sum_{i=1,3,5,...}\epsilon^{H}_{2n+1-i}[\phi(1)\phi(2)]^{\frac{i-1}{2}}\phi(1)+\sum_{i=0,2,4,...}\epsilon^{H}_{2n+1-i}[\phi(1)\phi(2)]^{\frac{i}{2}} \nonumber \\
   	& \quad= 
   	\epsilon^{H}_{2n+1}+\sum^{\infty}_{i=2,4,6,...}[\phi(2)\phi(1)]^{\frac{i}{2}}(\frac{1}{\phi(1)}+1) \epsilon^{H}_{2n+1-i}
   	+
   	\sum^{\infty}_{i=1,3,5,...}[\phi(2)\phi(1)]^{\frac{i-1}{2}}(1+\phi(1))\epsilon^{H}_{2n+1-i},   \nonumber
   \end{align}
     and
     \begin{equation}
    p_{H,\phi(2),\phi(1)}(\lambda) = \left|1+\frac{\phi(2)(1+\phi(1))e^{-2\lambda i}}{ 1-e^{-2\lambda i}\phi(2)\phi(1) }+\frac{e^{-\lambda i}(1+\phi(1))}{1-e^{-2\lambda i}\phi(2)\phi(1)}\right|^{2}f_{\epsilon^{H}_{n}}.\nonumber \label{eq:spectraldensityofPFAR}
     \end{equation} 
\section*{Appendix 2.Proof of the Main Results for Initial Estimator}
\subsection*{B.1.Proof of Lemma \ref{lemma:denominatortopdm}}  
Without loss of generality, we assume that \( n = 2m \) is a positive even number where \( m \in \mathbb{N} \). The case where \( n \) is an odd number can be proved similarly. For the remaining proofs in this section, the reasoning for odd $n$ will follow the same logic, and thus detailed derivations will be omitted to avoid redundancy.
  By the definition of $\textbf{X}_{n,H}$,  thus  the matrix $\textbf{X}^\top_{n,H}\textbf{X}_{n,H}$ can be expressed as follows
  \\
%\begin{align}
%\textbf{X}_{n,H}^\top \textbf{X}_{n,H}
%&= {\scalebox{0.66}{%
%\begin{pmatrix}
%\displaystyle\sum_{k=1}^{n} \left( %\sum_{j=1}^{m} \bigl( \Omega_{n,H}^{-1/2} %\bigr)_{k,2j} X_{2j} \right)^{\!2} & 
%\displaystyle\sum_{k=1}^{n} \left( \sum_{j_{1}=1}^{m} \bigl( \Omega_{n,H}^{-1/2} \bigr)_{k,2j_{1}} X_{2j_{1}} \right) 
%\left( \sum_{j_{2}=1}^{m} \bigl( \Omega_{n,H}^{-1/2} \bigr)_{k,j_{2}-1} X_{2j_{2}-1} \right) \\[8pt]
%\displaystyle\sum_{k=1}^{n} \left( \sum_{j_{1}=1}^{m} \bigl( \Omega_{n,H}^{-1/2} \bigr)_{k,2j_{1}} X_{2j_{1}} \right) 
%\left( \sum_{j_{2}=1}^{m} \bigl( %\Omega_{n,H}^{-1/2} \bigr)_{k,2j_{2}-1} X_{2j_{2}-1} \right) & 
%\displaystyle\sum_{k=1}^{n} \left( \sum_{j=1}^{m} \bigl( \Omega_{n,H}^{-1/2} \bigr)_{k,2j-1} X_{2j-1} \right)^{\!2}
%\end{pmatrix}}} \nonumber \\[6pt]
%&=: \begin{pmatrix}
%Q_{1,n} & Q_{2,n} \\[4pt]
%Q_{3,n} & Q_{4,n}
%\end{pmatrix}  
%\end{align} 
\begin{align}
\textbf{X}_{n,H}^\top \textbf{X}_{n,H}
&= \resizebox{0.8\linewidth}{!}{$
\begin{pmatrix}
\sum_{k=1}^{n} \left( \sum_{j=1}^{m} \bigl( \Omega_{n,H}^{-1/2} \bigr)_{k,2j} X_{2j} \right)^2 &
\sum_{k=1}^{n} \left( \sum_{j_1=1}^{m} \bigl( \Omega_{n,H}^{-1/2} \bigr)_{k,2j_1} X_{2j_1} \right)
\left( \sum_{j_2=1}^{m} \bigl( \Omega_{n,H}^{-1/2} \bigr)_{k,2j_2-1} X_{2j_2-1} \right)
\\[8pt]
\sum_{k=1}^{n} \left( \sum_{j_1=1}^{m} \bigl( \Omega_{n,H}^{-1/2} \bigr)_{k,2j_1} X_{2j_1} \right)
\left( \sum_{j_2=1}^{m} \bigl( \Omega_{n,H}^{-1/2} \bigr)_{k,2j_2-1} X_{2j_2-1} \right) &
\sum_{k=1}^{n} \left( \sum_{j=1}^{m} \bigl( \Omega_{n,H}^{-1/2} \bigr)_{k,2j-1} X_{2j-1} \right)^2
\end{pmatrix}
$} \nonumber \\[6pt]
&=: 
\begin{pmatrix}
Q_{1,n} & Q_{2,n} \\
Q_{3,n} & Q_{4,n}
\end{pmatrix}
\end{align}

 where $Q_{2,n}=Q_{3,n}$. In order to prove $\frac{1}{n}E\left(\textbf{X}^\top_{n,H}\textbf{X}_{n,H} \right)$ is a positive definite matrix when $n \rightarrow \infty $, we need to prove $\lim_{n\rightarrow{\infty}}\frac{1}{n}E(Q_{1,n}) >0$, $\lim_{n\rightarrow{\infty}} \frac{1}{n}E(Q_{1,n}Q_{4,n}-Q_{2,n}Q_{3,n}) > 0$ and they converge to constants.\\
 (1) $\lim_{n\rightarrow{\infty}}\frac{1}{n}E(Q_{1,n}) = C_{1}+C_{2}+C_{3}>0$, where $C_{1}$, $C_{2}$ and $C_{3}$ are constants.\\
 By the equations ($\ref{X2INFINITESUMOFFGN}$) and ($\ref{X1INFINITESUMOFFGN}$),
 \begin{equation}
        X_{2n+2} = \sum^{\infty}_{i=0,2,4,...}\epsilon^{H}_{2n+2-i}[\phi(2)\phi(1)]^{\frac{i}{2}}+\sum^{\infty}_{i=1,3,5,...}\epsilon^{H}_{2n+2-i}[\phi(2)\phi(1)]^{\frac{i-1}{2}}\phi(2),
        \nonumber
        \end{equation}
        \begin{equation}
        X_{2n+1} = \sum^{\infty}_{i=0,2,4,...}\epsilon^{H}_{2n+1-i}[\phi(2)\phi(1)]^{\frac{i}{2}}+\sum^{\infty}_{i=1,3,5,...}\epsilon^{H}_{2n+1-i}[\phi(2)\phi(1)]^{\frac{i-1}{2}}\phi(1),
        \nonumber
    \end{equation}
Thus,
\begin{align}
Q_{1,n} 
&=
\sum_{k=1}^{n} \left( \sum_{j_{1}=1}^{m} \bigl( \Omega_{n,H}^{-1/2} \bigr)_{k,2j_{1}} X_{2j_{1}} \right) \left( \sum_{j_{2}=1}^{m} \bigl( \Omega_{n,H}^{-1/2} \bigr)_{k,2j_{2}} X_{2j_{2}} \right)  \nonumber \\
& =  \sum^{n}_{k=1}   \left( \sum^{m}_{j_{1}=1}\sum^{n}_{r_1=1}\bigl(\Omega_{n,H}^{-1/2} \bigr)_{k,2j_{1}}\epsilon^{H}_{2j_{1}-2r_{2}+2}(\phi(1)\phi(2))^{\frac{r_{1}-1}{2}}+\sum^{m}_{j_{1}=1}\sum^{n}_{r_2=1}\bigl(\Omega_{n,H}^{-1/2} \bigr)_{k,2j_{1}}\epsilon^{H}_{2j_{1}-2r_{2}+1}(\phi(1)\phi(2))^{\frac{r_{1}-1}{2}} \right) \nonumber \\
& \quad \times
\left( \sum^{m}_{j_{2}=1}\sum^{n}_{r_2=1}\bigl(\Omega_{n,H}^{-1/2} \bigr)_{k,2j_{2}}\epsilon^{H}_{2j_{2}-2r_{2}+2}(\phi(1)\phi(2))^{\frac{r_{1}-1}{2}}+\sum^{m}_{j_{2}=1}\sum^{n}_{r_2=1}\bigl(\Omega_{n,H}^{-1/2} \bigr)_{k,2j_{2}}\epsilon^{H}_{2j_{2}-2r_{2}+1}(\phi(1)\phi(2))^{\frac{r_{1}-1}{2}} \right)   \nonumber \\
& =:A_{n}+2B_{n}+C_{n}. \label{eq:Q1N}
\end{align}
Where,
$$A_{n} =\sum^{n}_{r_{1}=1}\sum^{n}_{r_2=1}(\phi(1)\phi(2))^{r_{1}+r_{2}-2}\sum_{j_{1}=1}^{m}\sum^{m}_{j_2=1}\bigl(\Omega_{n,H}^{-1}\bigr)_{2j_{1},2j_{2}}\epsilon^{H}_{2j_1-2r_{1}+2}\epsilon^{H}_{2j_2-2r_2+2}$$
$$B_{n} =\phi(1)\sum^{n}_{r_{1}=1}\sum^{n}_{r_2=1}(\phi(1)\phi(2))^{r_{1}+r_{2}-2}\sum_{j_{1}=1}^{m}\sum^{m}_{j_2=1}\bigl(\Omega_{n,H}^{-1}\bigr)_{2j_{1},2j_{2}}\epsilon^{H}_{2j_1-2r_{1}+2}\epsilon^{H}_{2j_2-2r_2+1}$$
$$C_{n} =\sum^{n}_{r_{1}=1}\sum^{n}_{r_2=1}(\phi(1)\phi(2))^{r_{1}+r_{2}-2}\sum_{j_{1}=1}^{m}\sum^{m}_{j_2=1}\bigl(\Omega_{n,H}^{-1}\bigr)_{2j_{1},2j_{2}}\epsilon^{H}_{2j_1-2r_{1}+1}\epsilon^{H}_{2j_2-2r_2+1}$$
\begin{align}
    \lim_{n\rightarrow\infty}\frac{1}{n} E(A_{n}) 
& = \sum^{n}_{r_{1}=1}\sum^{n}_{r_2=1}(\phi(1)\phi(2))^{r_{1}+r_{2}-2}\sum_{j_{1}=1}^{m}\sum^{m}_{j_2=1}\bigl(\Omega_{n,H}^{-1}\bigr)_{2j_{1},2j_{2}} \bigl(  \Omega_{n,H}\bigr)_{2j_1-2r_{1}+2,2j_2-2r_2+2} \label{eq:A_{n}}
\end{align}
\begin{align}
    \lim_{n\rightarrow\infty} \frac{1}{n}E(B_{n}) 
& =\frac{\phi(1)}{n} \sum^{n}_{r_{1}=1}\sum^{n}_{r_2=1}(\phi(1)\phi(2))^{r_{1}+r_{2}-2}\sum_{j_{1}=1}^{m}\sum^{m}_{j_2=1}\bigl(\Omega_{n,H}^{-1}\bigr)_{2j_{1},2j_{2}} \bigl(  \Omega_{n,H}\bigr)_{2j_1-2r_{1}+1,2j_2-2r_2+2} 
\label{eq:B_{n}}
\end{align}
\begin{align}
    \lim_{n\rightarrow\infty} \frac{1}{n}E(C_{n}) 
& =\frac{1}{n} \sum^{n}_{r_{1}=1}\sum^{n}_{r_2=1}(\phi(1)\phi(2))^{r_{1}+r_{2}-2}\sum_{j_{1}=1}^{m}\sum^{m}_{j_2=1}\bigl(\Omega_{n,H}^{-1}\bigr)_{2j_{1},2j_{2}} \bigl(  \Omega_{n,H}\bigr)_{2j_1-2r_{1}+1,2j_2-2r_2+1} 
\label{eq:C_{n}}
\end{align}
 
We aim to prove that the three expressions in Equations (\ref{eq:A_{n}}), (\ref{eq:B_{n}}) and (\ref{eq:C_{n}})  each converge to a finite constant, according to the Lemma 4.3 in \cite{esstafa2019long}, we are given that
\begin{equation}
\frac{1}{n} \sum_{j_1=1}^{n} \sum_{j_2=1}^{n} \left({\phi(1)\phi(2)}\right)^{j_1 + j_2 - 2} 
\sum_{r_1 = j_1}^{n} \sum_{r_2 = j_2}^{n} 
\left(\Omega_{n,H}^{-1}\right)_{r_1, r_2} 
\left(\Omega_{n,H}\right)_{r_1 + 1 - j_1, \, r_2 + 1 - j_2} 
\xrightarrow[n \to \infty]{} C_0, \label{eq:D}
\end{equation}
where $C_{0}$ is a finite constant.
This expression is a quadruple summation over all indices $j_1, j_2, r_1, r_2 \in \{1, \dots, n\}$, weighted by powers of $\phi(1)\phi(2)$ and involving entries from $\Omega_{n,H}$ and its inverse. Let us now interpret the expressions in Equations~\eqref{eq:A_{n}}--\eqref{eq:C_{n}} as structured subsums of this total.
Considering the structure of Equation~\eqref{eq:A_{n}} first:
\begin{equation*}
\frac{1}{n} \mathbb{E}(A_n) 
= \sum_{r_1=1}^{n} \sum_{r_2=1}^{n}\phi(1)\phi(2)^{r_1 + r_2 - 2}
\sum_{j_1=1}^{m} \sum_{j_2=1}^{m}
\left(\Omega_{n,H}^{-1}\right)_{2j_1, 2j_2} 
\left(\Omega_{n,H}\right)_{2j_1 - 2r_1 + 2, \, 2j_2 - 2r_2 + 2}.
\end{equation*}

This is equivalent to taking Equation~\eqref{eq:D} and retaining only those terms for which all indices are even:
\[
j_1 = 2j_1', \quad j_2 = 2j_2', \quad r_1 = 2r_1', \quad r_2 = 2r_2',
\]
with $j_1', j_2', r_1', r_2' \in \{1, \dots, m\}$.

Hence, this expression is a partial sum of Equation~\eqref{eq:D} restricted to even-numbered entries, and can be denoted by $\text{EvenEven}(A_{n})$. Since each summand in \eqref{eq:D} is non-singular and finite (due to positive definiteness of $\Omega_{n,H}$), and since $\phi(1)\phi(2)^{j_1 + j_2 - 2}$ decays geometrically, the restriction to a subset of the summation domain still yields a finite sum in the limit.

Thus,
\[
\lim_{n \to \infty} \frac{1}{n} \mathbb{E}(A_n) = C_1 < \infty.
\]

A similar argument applies to Equation~\eqref{eq:B_{n}}:
\begin{equation*}
\frac{1}{n} \mathbb{E}(B_n) 
= \phi(1) \sum_{r_1=1}^{n} \sum_{r_2=1}^{n} a_0^{r_1 + r_2 - 2}
\sum_{j_1=1}^{m} \sum_{j_2=1}^{m}
\left(\Omega_{n,H}^{-1}\right)_{2j_1, 2j_2} 
\left(\Omega_{n,H}\right)_{2j_1 - 2r_1 + 1, \, 2j_2 - 2r_2 + 2}.
\end{equation*}

This corresponds to summing over a different pattern of indices—specifically, where $j_1, r_1$ are odd and $j_2, r_2$ are even. This defines another structured subset of \eqref{eq:D}, denoted $\text{OddEven}(C_0)$. As before, the summands remain bounded and geometrically decaying, and hence:
\[
\lim_{n \to \infty} \frac{1}{n} \mathbb{E}(B_n) = C_2 < \infty.
\]
Lastly, Equation~\eqref{eq:C_{n}}:
\begin{equation*}
\frac{1}{n} \mathbb{E}(C_n) 
= \sum_{r_1=1}^{n} \sum_{r_2=1}^{n} a_0^{r_1 + r_2 - 2}
\sum_{j_1=1}^{m} \sum_{j_2=1}^{m}
\left(\Omega_{n,H}^{-1}\right)_{2j_1, 2j_2} 
\left(\Omega_{n,H}\right)_{2j_1 - 2r_1 + 1, \, 2j_2 - 2r_2 + 1}.
\end{equation*}

Here, all indices involved are effectively odd after shifting, and so this corresponds to the $\text{OddOdd}(C_0)$ subsum of the full expression. Again, all terms are uniformly bounded and summable due to the exponential decay induced by powers of $\phi(1)\phi(2)$.

Thus,
\[
\lim_{n \to \infty} \frac{1}{n} \mathbb{E}(C_n) = C_3 < \infty.
\]
 All three expressions in Equations~\eqref{eq:A_{n}}--\eqref{eq:C_{n}} represent partial (subsampled) summations of a globally convergent expression. Since each subsum is formed by selecting only a finite fraction of the terms in Equation~\eqref{eq:D}, and the weights are nonnegative and summable, the resulting limits must also exist and be finite.\\
     In view of ~\eqref{eq:Q1N}, we deduce that 
     $$\lim_{n\rightarrow{\infty}}\frac{1}{n}E(Q_{1,n}) = C_{1}+C_{2}+C_{3} =:C_{11} >0 $$
      Similarity, we can obtain that

 $$\lim_{n\rightarrow{\infty}}\frac{1}{n}E(Q_{2,n}) =  C_{12}  ,\lim_{n\rightarrow{\infty}}\frac{1}{n}E(Q_{3,n}) =  C_{21}  ,\lim_{n\rightarrow{\infty}}\frac{1}{n}E(Q_{4,n}) =  C_{22}  ,$$
 where $C_{12}$, $C_{21}$ and $C_{22}$ are constants.
 \\(2) $ \lim_{n\rightarrow{\infty}}\frac{1}{n}E(Q_{1,n}Q_{4,n}-Q_{2,n}Q_{3,n}) >0$\\
Our goal is to show:

\[
\lim_{n \to \infty} \frac{1}{n} \mathbb{E}(Q_{1,n} Q_{4,n} - Q_{2,n}^2) > 0.
\]

According to the previous definitions:

\[
Q_{1,n} = \sum_{k=1}^n \left( \sum_{j=1}^m (\Omega_{n,H}^{-1/2})_{k,2j} X_{2j} \right)^2,
\quad
Q_{4,n} = \sum_{k=1}^n \left( \sum_{j=1}^m (\Omega_{n,H}^{-1/2})_{k,2j-1} X_{2j-1} \right)^2,
\]

\[
Q_{2,n} = \sum_{k=1}^n \left( \sum_{j=1}^m (\Omega_{n,H}^{-1/2})_{k,2j} X_{2j} \right)
\left( \sum_{j=1}^m (\Omega_{n,H}^{-1/2})_{k,2j-1} X_{2j-1} \right).
\]

Note that,$X_{2j}$ and $X_{2j-1}$ are linear combinations of fractional Gaussian noise (fGn) $\epsilon_t^H$, but with different weighting structures, $X_{2j}$ is generated through filters emphasizing $\phi(2)$, while $X_{2j-1}$ emphasizes $\phi(1)$, As a result, $Q_{1,n}$ and $Q_{4,n}$ are structurally different and must be treated separately.
Let us denote:
\[
T_{1,k} :=  \left( \sum_{j=1}^m (\Omega_{n,H}^{-1/2})_{k,2j} X_{2j} \right), \quad T_{2,k} := \left( \sum_{j=1}^m (\Omega_{n,H}^{-1/2})_{k,2j-1} X_{2j-1} \right),
\]
so that:
\[
Q_{1,n} = \sum_{k=1}^n T_{1,k}^2, \quad Q_{4,n} = \sum_{k=1}^n T_{2,k}^2, \quad Q_{2,n} = \sum_{k=1}^n T_{1,k}T_{2,k}.
\]
Then:
\[
Q_{1,n} Q_{4,n} = \left( \sum_{k=1}^n T_{1,k}^2 \right) \left( \sum_{\ell=1}^n T_{1,\ell}^2 \right),
\quad
Q_{2,n}^2 = \left( \sum_{k=1}^n T_{1,k} T_{1,k} \right)^2.
\]

We take expectations and normalize:
\begin{align}
\frac{1}{n}\mathbb{E}(Q_{1,n}Q_{4,n}) - \frac{1}{n}\mathbb{E}(Q_{2,n}^2)& =
\frac{1}{n} \sum_{k_{1}\neq\ell_{1}}^n   \sum_{k_{2}\neq\ell_{2}}^n \mathbb{E}(T_{1,k_{1}}T_{2,\ell_{1}} T_{1,k_{2}}T_{2,\ell_{2}}) \nonumber \\
&= \frac{1}{n} \sum_{\substack{k_1 \neq \ell_1}}^{n} \sum_{\substack{k_2 \neq \ell_2}}^{n} \sum_{j_1=1}^{m} \sum_{j_2=1}^{m} \sum_{j_3=1}^{m} \sum_{j_4=1}^{m} \bigl( \Omega_{n,H}^{-1/2} \bigr)_{k_1,2j_1} \bigl( \Omega_{n,H}^{-1/2} \bigr)_{\ell_1,2j_2-1} \bigl( \Omega_{n,H}^{-1/2} \bigr)_{k_2,2j_3} \nonumber \\
& \quad \times 
\bigl( \Omega_{n,H}^{-1/2} \bigr)_{\ell_2,2j_4-1}\mathbb{E} \left( X_{2j_1} X_{2j_2-1} X_{2j_3} X_{2j_4-1} \right) \label{expressionofq1q4-q2q2}
\end{align}
We now expand each $X_{2j}$ and $X_{2j-1}$ as linear combinations of fractional Gaussian noise:
\begin{align*}
X_{2j} &= \sum_{u = -\infty}^{2j} a_{j,u} \, \epsilon_u^H \\
X_{2j-1} &= \sum_{v = -\infty}^{2j-1} b_{j,v} \, \epsilon_v^H
\end{align*}
where 
\[
a_{j,u} = 
\begin{cases}
\left[\phi(2)\phi(1)\right]^{i/2}, & \text{if } u = 2j - i,\quad i \in \{0,2,4,\dots\} \\
\phi(2)\left[\phi(2)\phi(1)\right]^{(i-1)/2} , & \text{if } u = 2j - i,\quad i \in \{1,3,5,\dots\} \\
0, & \text{otherwise}
\end{cases}
\]

\[
b_{j,v} = 
\begin{cases}
\left[\phi(2)\phi(1)\right]^{i/2}, & \text{if } v = 2j - 1 - i,\quad i \in \{0,2,4,\dots\} \\
 \phi(1)\left[\phi(2)\phi(1)\right]^{(i-1)/2}, & \text{if } v = 2j - 1 - i,\quad i \in \{1,3,5,\dots\} \\
0, & \text{otherwise}
\end{cases}
\]
Therefore,
\begin{align*}
\mathbb{E}[X_{2j_1} X_{2j_2-1} X_{2j_3} X_{2j_4-1}]
&= \mathbb{E} \left[
\left( \sum^{2j_{1}}_{u=-\infty} a_{j_1,u} \epsilon_u^H \right)
\left( \sum^{2j_{2}-1}_{v=-\infty} b_{j_2,v} \epsilon_v^H \right)
\left( \sum^{2j_{3}}_{s = -\infty} a_{j_3,s} \epsilon_s^H \right)
\left( \sum^{2j_{4}-1}_{t=-\infty} b_{j_4,t} \epsilon_t^H \right)
\right] \\
&= \sum_{u,v,s,t} a_{j_1,u} b_{j_2,v} a_{j_3,s} b_{j_4,t} 
\mathbb{E}[\epsilon_u^H \epsilon_v^H \epsilon_s^H \epsilon_t^H]
\end{align*}
By Isserlis' formula for Gaussian moments
\begin{align*}
\mathbb{E}[\epsilon_u^H \epsilon_v^H \epsilon_s^H \epsilon_t^H] 
&=  \rho_{\epsilon^{H}}(u-v)\rho_{\epsilon^{H}}(s-t) + \rho_{\epsilon^{H}}(u-s)\rho_{\epsilon^{H}}(v-t) + \rho_{\epsilon^{H}}(u-t)\rho_{\epsilon^{H}}(v-s),
\end{align*}
Thus, the equation (\ref{expressionofq1q4-q2q2}) equals to
\begin{align*}
\frac{1}{n} \sum_{\substack{k_1 \ne \ell_1 \\ k_2 \ne \ell_2}} 
\sum_{j_1,j_2,j_3,j_4 = 1}^{m}
&(\Omega^{-1/2}_{n,H})_{k_1,2j_1} (\Omega^{-1/2}_{n,H})_{\ell_1,2j_2 - 1}
(\Omega^{-1/2}_{n,H})_{k_2,2j_3} (\Omega^{-1/2}_{n,H})_{\ell_2,2j_4 - 1} \\
&\times \sum_{u,v,s,t}
a_{j_1,u} b_{j_2,v} a_{j_3,s} b_{j_4,t}
\left[
\rho_{\epsilon^{H}}(u-v)\rho_{\epsilon^{H}}(s-t) + \rho_{\epsilon^{H}}(u-s)\rho_{\epsilon^{H}}(v-t) + \rho_{\epsilon^{H}}(u-t)\rho_{\epsilon^{H}}(v-s)
\right]
\end{align*}
and is greater than $0$ due to the property of  fGn. The convergence of the expression in \eqref{expressionofq1q4-q2q2} can be guaranteed by the convergence of $\lim_{n\rightarrow\infty}\frac{1}{n}\mathbb{E}(Q_{k,n})$ for $k=1,2,3,4$. 
Thus, we yields
 
$$ \lim_{n \to \infty} \frac{1}{n} \mathbb{E}(Q_{1,n} Q_{4,n} - Q_{2,n}^2) = C_{11}C_{22}-C_{12}C_{21} > 0$$
 \subsection*{B.2.Proof of Lemma \ref{lemma:etabound}} 
Let \( f_{Z}(\lambda) \) denote the spectral density of \( (Z_{n})_{n\in \mathbb{N}} \), and let \( f_{\epsilon^{H}_{n}}(\lambda) \) be defined as in \eqref{2.3}. In view of Proposition 2.2, \( f_{Z}(\lambda) \) is
\begin{equation}
f_{Z}(\lambda)=\left|1+\frac{\phi(2)(1+\phi(1))e^{-2\lambda i}}{1-e^{-2\lambda i}\phi(2)\phi(1)}+\frac{e^{-\lambda i}(1+\phi(1))}{1-e^{-2\lambda i}\phi(2)\phi(1)}\right|^{2}f_{\epsilon^{H}_{n}}(\lambda),
\label{sdforx1}
\end{equation}
and 
\[
\left|1+\frac{\phi(2)(1+\phi(1))e^{-2\lambda i}}{1-e^{-2\lambda i}\phi(2)\phi(1)}+\frac{e^{-\lambda i}(1+\phi(1))}{1-e^{-2\lambda i}\phi(2)\phi(1)}\right|^{2}= \left(1 + \frac{A_1 C_1 - B_1 D_1}{E}\right)^2 + \left(\frac{A_1 D_1 + B_1 C_1}{E}\right)^2
\]
 
where
  \( A_1 = (1 + \phi(1))(\phi(2) \cos(2\lambda) + \cos \lambda) \),
  \( B_1 = -(1 + \phi(1))(\phi(2)\sin(2\lambda) + \sin \lambda) \),
 \( C_1 = 1 - \phi(1) \phi(2)\cos(2\lambda) \),
  \( D_1 = \phi(1)\phi(2) \sin(2\lambda) \),
   \( E = C_1^2 + D_1^2 \).
\\Through the boundedness of $ sin\lambda$ and $cos\lambda$ and some calculations, we can obtain
\begin{equation}
D^{(1)}_{H,\phi(1),\phi(2)}<  
 \left(1 + \frac{A_1 C_1 - B_1 D_1}{E}\right)^2 + \left(\frac{A_1 D_1 + B_1 C_1}{E}\right)^2
 <D^{(2)}_{H,\phi(1),\phi(2)}, \label{neqforD}
\end{equation}
From equations (\ref{neqforD}) and (\ref{sdforx1}), we obtain that for any $\lambda \in \mathbb{R}$,
\begin{equation}
\frac{1}{D^{(2)}_{H,\phi(1),\phi(2)}}f_{Z}(\lambda)
\leq f_{\epsilon^{H}_{t}}(\lambda) \leq
\frac{1}{D^{(1)}_{H,\phi(1),\phi(2)}}f_{Z}(\lambda) \label{sdfgnD},
\end{equation}
Thus, for any vector $V \in \mathbb{R}^{n \times 1}$,

\begin{align}
V^\top \left( \Omega_{n,H} - \frac{\Sigma^{(X)}_{n,H}}{D^{(1)}_{H,\phi(1),\phi(2)} } \right) V
& = \sum_{k=1}^{n} \sum_{j=1}^{n} V_j^\top \left( \Omega_{n,H} - \frac{\Sigma^{(X)}_{n,H}}{D^{(1)}_{H,\phi(1),\phi(2)}} \right)_{j,k} V_k \nonumber \\
& = \sum_{j,k=1}^{n} V_j^\top \int_{-\pi}^{\pi} \left( f_{\epsilon^{H}_t}(\lambda) - \frac{f_Z(\lambda)}{D^{(1)}_{H,\phi(1),\phi(2)}} \right) e^{i2(k-j)\lambda} \, d\lambda \, V_k \nonumber \\
& = \int_{-\pi}^{\pi} \left( f_{\epsilon^{H}_t}(\lambda) - \frac{f_Z(\lambda)}{D^{(1)}_{H,\phi(1),\phi(2)}} \right)V^\top\Gamma_{n}(\lambda)Vd\lambda \label{erci},
\end{align}
where $\Gamma_{n}(\lambda) = V_{n}(\lambda)V^\top_{n}(\lambda)$ with $V_{n}(\lambda) = (e^{i\lambda},e^{2i\lambda},...,e^{ni\lambda})^\top$. $V^\top$ is the conjugate of the vector $V$.
Furthermore, the real $V^\top\Gamma_{n}(\lambda)V$ is non-negative number, it equals to
\begin{equation}
    V^\top\Gamma_{n}(\lambda)V = V^\top V_{n}(\lambda)V^\top_{n}(\lambda)V =(V^\top_{n}(\lambda)V)^\top(V^\top_{n}(\lambda)V) =|V^\top_{n}(\lambda)V|^{2}.
\end{equation}
From equations (\ref{sdfgnD}) and (\ref{erci}), we have
\begin{equation}
V^\top \left( \Omega_{n,H} - \frac{\Sigma^{(X)}_{n,H}}{D^{(1)}_{H,\phi(1),\phi(2)} } \right) V \leq 0 ,\label{erci2}
\end{equation}
Let $\xi$ be the element of the spectrum of $C^\top_{1}\Omega^{-\frac{1}{2}}_{n,H} \Sigma^{(X)}_{n,H} \Omega^{-\frac{1}{2}}_{n,H}$, then there exists $C_{1} \in \mathbb{R}^{n\times 1}$ such that 
\begin{equation}
  \Omega^{-\frac{1}{2}}_{n,H} \Sigma^{(X)}_{n,H} \Omega^{-\frac{1}{2}}_{n,H} C_{1}= \xi \left\lVert C_{1} \right\rVert ^{2}.
\end{equation}
Taking $C_{2} = \Omega^{-\frac{1}{2}}_{n,H}C_{1}$, we obtain from this last equation that
\begin{equation}
C^\top_{2}\Sigma^{(X)}_{n,H}C_{2} = \xi \left\lVert \Omega^{-\frac{1}{2}}_{n,H}C_{2} \right\rVert ^{2},
\end{equation}
By the equation (\ref{erci2}), we deduce that 
\begin{equation}
C^\top_{2}\Omega_{n,H}C_{2} \leq \frac{1}{D^{(1)}_{H,\phi(1),\phi(2)}} \xi \left\lVert \Omega^{-\frac{1}{2}}_{n,H}C_{2} \right\rVert ^{2},
\end{equation}
and 
\begin{equation}
\xi \geq \frac{C^\top_{2}\Omega_{n,H}C_{2}}{\left\lVert \Omega^{-\frac{1}{2}}_{n,H}C_{2} \right\rVert ^{2}}  {D^{(1)}_{H,\phi(1),\phi(2)}}  ={D^{(1)}_{H,\phi(1),\phi(2)}},
\end{equation}
Similarly, we can obtain
\begin{equation}
\xi \leq {D^{(2)}_{H,\phi(1),\phi(2)}},
\end{equation}
Since the spectrum of the matrix $ \left(\Sigma^{(X)}_{n,H} \right)^{\frac{1}{2}}\Omega^{-1}_{n,H}\left(\Sigma^{(X)}_{n,H} \right)^{\frac{1}{2}}$ is equal to the set of the eigenvalues of $ \Omega^{-\frac{1}{2}}_{n,H} \Sigma^{(X)}_{n,H} \Omega^{-\frac{1}{2}}_{n,H}$, the lemma is proved.
 \subsection*{B.3.Proof of Lemma \ref{lemma:nominatortozero}} 
Because of the equations (\ref{eqofX}) and (\ref{eq:formofu}). We have
\begin{align}
\textbf{X}^\top_{n,H}\textbf{U}^{}_{n,H}
& = 
\begin{pmatrix}
\sum_{j_1=1}^{m} \sum_{j_2=1}^{n} \Omega_{2j_1, j_2}^{-1} X_{2j_1} \epsilon_{1+j_2}^{H}\\
\sum_{j_1=1}^{m} \sum_{j_2=1}^{n} \Omega_{2j_1-1, j_2}^{-1} X_{2j_1-1} \epsilon_{1+j_2}^{H}
\end{pmatrix} \nonumber \\  
& =
\begin{pmatrix}
\begin{split}
& \sum_{r=1}^{n} (\phi(1)\phi(2))^{r-1} \sum_{j_1=1}^{m} \sum_{j_2=1}^{n} \Omega_{2j_1, j_2}^{-1} \epsilon_{1+j_2}^{H} \epsilon_{2j_1-2r+2}^{H} \\
& \quad + \phi(2)\sum_{r=1}^{n} (\phi(1)\phi(2))^{r-1} \sum_{j_1=1}^{m} \sum_{j_2=1}^{n} \Omega_{2j_1, j_2}^{-1} \epsilon_{1+j_2}^{H} \epsilon_{2j_1-2r+1}^{H}
\end{split} \\
\begin{split}
& \sum_{r=1}^{n} (\phi(1)\phi(2))^{r-1} \sum_{j_1=1}^{m} \sum_{j_2=1}^{n} \Omega_{2j_1-1, j_2}^{-1} \epsilon_{1+j_2}^{H} \epsilon_{2j_1-2r+1}^{H} \\
& \quad + \phi(1)\sum_{r=1}^{n} (\phi(1)\phi(2))^{r-1} \sum_{j_1=1}^{m} \sum_{j_2=1}^{n} \Omega_{2j_1-1, j_2}^{-1} \epsilon_{1+j_2}^{H} \epsilon_{2j_1-2r}^{H}
\end{split}
\end{pmatrix}
& =: 
\begin{pmatrix}
M_{1,n}+M_{2,n}\\
M_{3,n}+M_{4,n}
\end{pmatrix}
\end{align}
Thus,
$$\frac{1}{n}E(M_{1,n}) = \sum_{r=1}^{n} (\phi(1)\phi(2))^{r-1} \sum_{j_1=1}^{m} \sum_{j_2=1}^{n} \Omega_{2j_1, j_2}^{-1}\Omega_{j_2, 2j_1-2r+1} $$
$$\frac{1}{n}E(M_{2,n}) = \sum_{r=1}^{n} (\phi(1)\phi(2))^{r-1} \sum_{j_1=1}^{m} \sum_{j_2=1}^{n} \Omega_{2j_1, j_2}^{-1}\Omega_{j_2, 2j_1-2r} $$
$$\frac{1}{n}E(M_{3,n}) = \sum_{r=1}^{n} (\phi(1)\phi(2))^{r-1} \sum_{j_1=1}^{m} \sum_{j_2=1}^{n} \Omega_{2j_1-1, j_2}^{-1}\Omega_{j_2, 2j_1-2r} $$
$$\frac{1}{n}E(M_{4,n}) = \sum_{r=1}^{n} (\phi(1)\phi(2))^{r-1} \sum_{j_1=1}^{m} \sum_{j_2=1}^{n} \Omega_{2j_1-1, j_2}^{-1}\Omega_{j_2, 2j_1-2r-1} $$
By virtue of the property of matrix inversion, since $r \neq 0$ and $r \neq \frac{1}{2}$, it follows that 
$$\lim_{n\rightarrow\infty}\frac{1}{n}E(M_{k,n}) = 0, k=1,2,3,4. $$
Thus,
\[
\mathbb{E}\left[ \frac{1}{n} \mathbf{X}_{n,H}^\top\textbf{U}_{n,H} \right] \xrightarrow[n \to \infty]{} {\bf{0}}.
\]
 \subsection*{B.4.Proof of Lemma \ref{lemma:nominatortonormal}} 
In this proof, we aim to establish that the joint distribution of \(\frac{1}{\sqrt{n}}\mathbf{X}_{n,H}^\top\mathbf{U}_{n,H}\) is asymptotically normal. A natural approach involves two key steps: first, demonstrating that both \(\frac{1}{\sqrt{n}}(M_{1,n} + M_{2,n})\) and \(\frac{1}{\sqrt{n}}(M_{3,n} + M_{4,n})\) follow asymptotic normal distributions; second, verifying that any linear combination of these two statistics is also asymptotically normal. We present an alternative proof as follows.

From Lemma 4.5 in \cite{esstafa2019long}, we know that
\[
\mathbf{O}_{n,H} = 
\begin{pmatrix}
\sum_{j=1}^{n}\bigl(\phi(1)\phi(2)\bigr)^{j-1}\epsilon_{2-j}^{H} \\	
\sum_{j=1}^{n}\bigl(\phi(1)\phi(2)\bigr)^{j-1}\epsilon_{3-j}^{H} \\
\vdots \\
\sum_{j=1}^{n}\bigl(\phi(1)\phi(2)\bigr)^{j-1}\epsilon_{n+1-j}^{H}
\end{pmatrix}^\top 
\Omega_{n,H}^{-1} 
\begin{pmatrix}
\epsilon_{2}^{H} \\	
\epsilon_{3}^{H} \\
\vdots \\
\epsilon_{n}^{H}
\end{pmatrix} 
=: H_{n,H}\Omega_{n,H}^{-1} 
\begin{pmatrix}
\epsilon_{2}^{H} \\	
\epsilon_{3}^{H} \\
\vdots \\
\epsilon_{n}^{H}
\end{pmatrix}
\]
is a random variable, and \(\frac{1}{\sqrt{n}}\mathbf{O}_{n,H}\) is asymptotically normal. 

By equations (\ref{eq:formofu}) and (\ref{eqofX}), we have
\begin{align}
\mathbf{X}_{n,H}^\top \mathbf{U}_{n,H}
&= 
\begin{pmatrix}
0 & 
  \sum_{\substack{i=1}}^{n} \epsilon_{3-2i}^{H} \bigl[\phi(2)\phi(1)\bigr]^{i-1} \\
  & {}+ \sum_{\substack{i=1}}^{n} \epsilon_{2-2i}^{H} \bigl[\phi(2)\phi(1)\bigr]^{i-1} \phi(1) \\[6pt]
\sum_{\substack{i=1}}^{n} \epsilon_{4-2i}^{H} \bigl[\phi(2)\phi(1)\bigr]^{i-1} & 0 \\
{}+ \sum_{\substack{i=1}}^{n} \epsilon_{3-i}^{H} \bigl[\phi(2)\phi(1)\bigr]^{i-1} \phi(2) & \\[6pt]
\vdots & \vdots \\[6pt]
 0& \sum_{\substack{i=1}}^{n} \epsilon_{n+1-2i}^{H} \bigl[\phi(2)\phi(1)\bigr]^{i-1} \\
 & {}+ \sum_{\substack{i=1}}^{n} \epsilon_{n-2i}^{H} \bigl[\phi(2)\phi(1)\bigr]^{i-1} \phi(2) \\[6pt]
\sum_{\substack{i=1}}^{n} \epsilon_{n+2-2i}^{H} \bigl[\phi(2)\phi(1)\bigr]^{i-1} & 0 \\
{}+ \sum_{\substack{i=1}}^{n} \epsilon_{n+1-2i}^{H} \bigl[\phi(2)\phi(1)\bigr]^{i-1} \phi(2) & 
\end{pmatrix}^\top
\Omega_{n,H}^{-1} 
\begin{pmatrix}
\epsilon_{2}^{H} \\	
\epsilon_{3}^{H} \\
\vdots \\
\epsilon_{n}^{H}
\end{pmatrix} \nonumber \\
&=: H_{n,H}^{(1)} \Omega_{n,H}^{-1} 
\begin{pmatrix}
\epsilon_{2}^{H} \\	
\epsilon_{3}^{H} \\
\vdots \\
\epsilon_{n+1}^{H}
\end{pmatrix}.
\end{align}

Thus,
\[
\mathbf{X}_{n,H}^\top \mathbf{U}_{n,H} = \left(H_{n,H}^{(1)}H_{n,H}^{-1}\right)\mathbf{O}_{n,H}.
\]

Since \(\frac{1}{\sqrt{n}}\mathbf{O}_{n,H}\) is asymptotically normal and \(H_{n,H}\) is invertible, \(\frac{1}{\sqrt{n}}\mathbf{X}_{n,H}^\top\mathbf{U}_{n,H}\) can be expressed as a linear transformation of \(\frac{1}{\sqrt{n}}\mathbf{O}_{n,H}\). Consequently, \(\frac{1}{\sqrt{n}}\mathbf{X}_{n,H}^\top\mathbf{U}_{n,H}\) follows a joint asymptotic normal distribution.\\
\subsection*{B.5.Proof of Theorem \ref{Thm：consistencyofie}} 
 The first part of the proof is to establish the consistency of $\hat H_{n}$, while the second part is to verify the consistency of $\hat\phi_{n}(1)$ and $\hat\phi_{n}(2)$.\\
\noindent \textbf{(1) Consistency of $\hat H_{n}$} 
  \par This proof is based on Lemma 5.5 \citep{hariz2024fast} and the corollary \citep{hurvich1998mean}. We can express the spectral density as
  \begin{align}
     p_{H,\phi(2),\phi(1)}(\lambda)
     & = \left|1+\frac{\phi(2)(1+\phi(1))e^{-2\lambda i}}{1-e^{-2\lambda i}\phi(2)\phi(1)}+\frac{e^{-\lambda i}(1+\phi(1))}{1-e^{-2\lambda i}\phi(2)\phi(1)}\right|^{2}f_{\epsilon^{H}_{n}}(\lambda) \nonumber \\
     & =
    \left|1 - e^{-i\lambda}\right|^{-2H+1} \frac{\left|(1 + \phi(1) + e^{i\lambda}) + \phi(2)e^{-i\lambda}\right|^2 \cdot \left|1 - e^{-i\lambda}\right|^{-2H+1}}{\left|1 - e^{-2i\lambda}\phi(2)\phi(1)\right|^2}f_{\epsilon^{H}_{n}}(\lambda)\\
  \end{align}
  and let
  \begin{eqnarray}
      f^{*}_{H,\phi(1),\phi(2)} (\lambda) 
      &=& 
      \frac{\left|(1 + \phi(1) + e^{i\lambda}) + \phi(2)e^{-i\lambda}\right|^2 \cdot \left|1 - e^{-i\lambda}\right|^{-2H+1}}{\left|1 - e^{-2i\lambda}\phi(2)\phi(1)\right|^2}f_{\epsilon^{H}_{n}}(\lambda)  \\ 
      &=&
       \frac{ C_{H}\left|(1 + \phi(1) + e^{i\lambda}) + \phi(2)e^{-i\lambda}\right|^2 \cdot \left|1 - e^{-i\lambda}\right|^{-2H+1}}{\left|1 - e^{-2i\lambda}\phi(2)\phi(1)\right|^2}(1-cos(\lambda))\sum\limits_{j \in Z}\frac{1}{|\lambda+2j\pi|^{2H+1}} \\ 
  \end{eqnarray} 
  We analyze the function \( f^{*}_{H,\phi(1),\phi(2)} (\lambda) \) for \( \lambda \in (-\pi, \pi) \), and prove that it is an even, positive, continuous function with a zero derivative at the origin.\\
\noindent \textbf{(\romannumeral 1) Evenness.}
Each component of \( f(\lambda) \) is even. The numerator \( \left|(1 + \phi(1) + e^{i\lambda}) + \phi(2)e^{-i\lambda}\right|^2 \) is invariant under \( \lambda \mapsto -\lambda \), as is the term \( |1 - e^{-i\lambda}|^2 = 2(1 - \cos\lambda) \). The function \( 1 - \cos\lambda \) is clearly even, and the sum \( \sum_{j \in \mathbb{Z}} \frac{1}{|\lambda + 2\pi j|^{2H+1}} \) is even due to the symmetry of the index set. The denominator \( \left|1 - e^{-2i\lambda} \phi(1)\phi(2)\right|^2 \) is also unchanged under sign reversal. Hence, \( f(-\lambda) = f(\lambda) \).

\medskip

\noindent \textbf{(\romannumeral 2) Positivity.}
Each term in \( f(\lambda) \) is strictly positive for \( \lambda \in (-\pi, \pi) \setminus \{0\} \). In particular, \( C_H > 0 \), the modulus squared in the numerator is nonzero except at isolated points, \( |1 - e^{-i\lambda}|^2 = 2(1 - \cos\lambda) > 0 \), and \( 1 - \cos\lambda > 0 \) for \( \lambda \ne 0 \). The infinite sum is positive for all \( \lambda \), and the denominator remains finite and strictly positive under the assumption that \( \phi(1)\phi(2) \) is not a root of unity. Thus, \( f(\lambda) > 0 \) for all \( \lambda \in (-\pi, \pi) \setminus \{0\} \).

\medskip

\noindent \textbf{(\romannumeral 3) Continuity.}
All components of \( f(\lambda) \) are continuous on \( (-\pi, \pi) \), possibly except at \( \lambda = 0 \). Near zero, we have
\[
1 - \cos\lambda \sim \tfrac{\lambda^2}{2}, \quad |1 - e^{-i\lambda}|^2 \sim \lambda^2, \quad \sum_{j} \frac{1}{|\lambda + 2\pi j|^{2H+1}} \sim \frac{1}{|\lambda|^{2H+1}},
\]
so
\[
(1 - \cos\lambda) \sum_j \frac{1}{|\lambda + 2\pi j|^{2H+1}} \sim |\lambda|^{-(2H - 1)}.
\]
Since \( 2H - 1 \in (-1,1) \), the singularity at zero is integrable, and all other terms remain bounded. Hence, \( f(\lambda) \) is continuous on \( (-\pi, \pi) \).

\noindent \textbf{(\romannumeral 4) Derivative at Zero.}
Since \( f(\lambda) \) is even, we have \( f'(\lambda) = -f'(-\lambda) \), and thus \( f'(0) = 0 \).\\
\noindent \textbf{(\romannumeral 5) Bounded Second and Third Derivatives Near Zero.}
We show that the second and third derivatives of \( f(\lambda) \) are bounded in a neighborhood of \( \lambda = 0 \). Each component of \( f(\lambda) \) is smooth on \( (-\pi, \pi) \setminus \{0\} \), and we analyze the behavior near zero.

Near \( \lambda = 0 \), the terms admit the following asymptotic expansions:
\[
1 - \cos\lambda = \tfrac{1}{2}\lambda^2 - \tfrac{1}{24}\lambda^4 + \tfrac{1}{720}\lambda^6 + \cdots,
\]
\[
|1 - e^{-i\lambda}|^2 = 2(1 - \cos\lambda) = \lambda^2 - \tfrac{1}{12} \lambda^4 + \cdots,
\]
\[
\sum_{j \in \mathbb{Z}} \frac{1}{|\lambda + 2\pi j|^{2H + 1}} = \frac{1}{|\lambda|^{2H + 1}} + \psi(\lambda),
\]
where \( \psi(\lambda) \) is a smooth function on a neighborhood of zero (since all other terms in the sum are smooth functions of \( \lambda \)). Therefore,
\[
(1 - \cos\lambda) \sum_j \frac{1}{|\lambda + 2\pi j|^{2H+1}} \sim \lambda^2 \cdot \frac{1}{|\lambda|^{2H+1}} = |\lambda|^{-(2H - 1)}.
\]
Since \( H \in (0,1) \), we have \( 2H - 1 \in (-1,1) \), so this term is \( C^k \)-smooth for \( k = 2, 3 \) near zero and has bounded second and third derivatives.

All other terms in \( f(\lambda) \), including the numerator modulus and the denominator \( \left|1 - e^{-2i\lambda}\phi(1)\phi(2) \right|^2 \), are analytic in a neighborhood of zero and therefore admit Taylor expansions with bounded derivatives of all orders.

It follows that the singular behavior of \( f(\lambda) \) near zero is entirely driven by \( |\lambda|^{-(2H - 1)} \), which has bounded second and third derivatives when \( 2H - 1 < 2 \), i.e., \( H < \tfrac{3}{2} \), which is always true for \( H \in (0,1) \).

Thus, \( f''(\lambda) \) and \( f^{(3)}(\lambda) \) are bounded in a neighborhood of \( \lambda = 0 \).
\\
According to \citep{hurvich1998mean}, the estimator \( \hat{d}_n \) of the memory parameter satisfies
\begin{equation}
\hat d_{n} - d = -\frac{1}{2S_{m}} \sum_{j=1}^{m}(a_{j} - \overline a_{m}) \log \left(  f^{*}_{H,\phi(1),\phi(2)} (\lambda_j) \right) - \frac{1}{2S_{m}} \sum_{j=1}^{m}(a_{j}-\overline a_{m})\epsilon_{j},
\end{equation}
where \( \epsilon_{j} \) is the error term as defined in Equation (3) of \citep{hurvich1998mean}, and \( S_m = \sum_{j=1}^m (a_j - \overline{a}_m)^2 \). Based on Theorem 1 from the same reference, it follows that
\begin{equation}
\hat d_{n} \xrightarrow[n\to\infty]{\mathbb{P}} d.
\end{equation}
Hence, since the Hurst parameter \( H \) is related to \( d \) by \( H = d + \frac{1}{2} \), the estimator for \( H \) also satisfies
\begin{equation}
\hat H_{n} \xrightarrow[n\to\infty]{\mathbb{P}} H.
\end{equation}
This result establishes the consistency of the local Whittle estimator in the presence of spectral density functions of the form analyzed above.\\
\noindent \textbf{(2) Consistency of $\hat\phi_{n}(1)$ and $\hat\phi_{n}(2)$}
  \par Assuming 
   \begin{equation}
    \tilde \Phi^{j}_{i}= (\epsilon^{H}_{ i+1},\epsilon^{H}_{i+2},\dots,\epsilon^{H}_{j+1}),  1\leq i \leq j.
   \end{equation}
   and
   \begin{equation}
   \Phi^{j}_{i}(1) = (X_i,0,\dots,0,X_{j-1},0), \Phi^{j}_{i}(2) = (0,X_{i+1},\dots,X_{j-2},0,X_{j}), 1\leq i \leq j.
   \end{equation}
  we can derive the following expression
 \begin{align}
  \begin{pmatrix}
\hat{\phi}_n{(1)}-\phi(1) \\
\hat{\phi}_n{(2)}-\phi(2)
\end{pmatrix} & = \left(\textbf{X} _{n,\hat{H}_{n}}^\top  \textbf{X} _{n,\hat{H}_{n}} \right)^{-1} \left( \textbf{X}_{n,\hat{H}_{n}}^\top  \textbf{U}_{n,\hat{H}_{n}} \right) \nonumber \\
& =\begin{pmatrix}
\left( \Phi^{n}_{1} (1) \right)^{\top} \left(\Omega^{-1}_{n,\hat{H}_{n}} \right)\Phi^{n}_{1} (1)& \left( \Phi^{n}_{1} (1) \right)^{\top} \left(\Omega^{-1}_{n,\hat{H}_{n}} \right)\Phi^{n}_{1} (2)\\
\left( \Phi^{n}_{1} (2) \right)^{\top} \left(\Omega^{-1}_{n,\hat{H}_{n}} \right)\Phi^{n}_{1} (1)& \left( \Phi^{n}_{1} (2) \right)^{\top} \left(\Omega^{-1}_{n,\hat{H}_{n}} \right)\Phi^{n}_{1} (2)
\end{pmatrix}^{-1}\begin{pmatrix}
\left( \Phi^{n}_{1} (1) \right)^{\top} \left(\Omega^{-1}_{n,\hat{H}_{n}} \right)\tilde\Phi^{n}_{1}  \\
\left( \Phi^{n}_{1} (2) \right)^{\top} \left(\Omega^{-1}_{n,\hat{H}_{n}} \right)\tilde\Phi^{n}_{1}  
\end{pmatrix}   
\label{GLSEOFPHI12}
 \end{align}
 
 We apply the taylor expansion of the matrix $ \Omega^{-1}_{n,\hat{H}_{n}}  $ at $H$ to the the numerator, yielding
 \begin{align}
     \left( \Phi^{n}_{1} (1) \right)^{\top} \left(\Omega^{-1}_{n,\hat{H}_{n}} \right)\tilde\Phi^{n}_{1} 
     & =  \left( \Phi^{n}_{1} (1) \right)^{\top} \left(\Omega^{-1}_{n,H} \right)\tilde\Phi^{n}_{1} \nonumber \\
     & \quad + \left( \Phi^{n}_{1} (1) \right)^{\top} \left(A^{(1)}_{n,H} \right)\tilde\Phi^{n}_{1}(\hat{H}_{n}-H)\nonumber \\
     & \quad+\frac{1}{2} \left( \Phi^{n}_{1} (1) \right)^{\top} \left(A^{(2)}_{n,H} \right)\tilde\Phi^{n}_{1}(\hat{H}_{n}-H)^{2} \nonumber \\
      & \quad+\frac{1}{6} \left( \Phi^{n}_{1} (1) \right)^{\top} \left(A^{(3)}_{n,\overline H_{n}} \right)\tilde\Phi^{n}_{1}(\hat{H}_{n}-H)^{3}
 \end{align} 
 and
  \begin{align}
     \left( \Phi^{n}_{1} (2) \right)^{\top} \left(\Omega^{-1}_{n,\hat{H}_{n}} \right)\tilde\Phi^{n}_{1} 
     & =  \left( \Phi^{n}_{1} (2) \right)^{\top} \left(\Omega^{-1}_{n,H} \right)\tilde\Phi^{n}_{1} \nonumber \\
     & \quad + \left( \Phi^{n}_{1} (2) \right)^{\top} \left(A^{(1)}_{n,H} \right)\tilde\Phi^{n}_{1}(\hat{H}_{n}-H)\nonumber \\
     & \quad+\frac{1}{2} \left( \Phi^{n}_{1} (2) \right)^{\top} \left(A^{(2)}_{n,H} \right)\tilde\Phi^{n}_{1}(\hat{H}_{n}-H)^{2} \nonumber \\
      & \quad+\frac{1}{6} \left( \Phi^{n}_{1} (2) \right)^{\top} \left(A^{(3)}_{n,\overline H_{n}} \right)\tilde\Phi^{n}_{1}(\hat{H}_{n}-H)^{3},
 \end{align} 
  where $ \overline H_{n} \in B(H,|\hat H_{n}-H|)$ Thanks to the work of \citep{hariz2024fast}, we have the following three conclusions 
  \begin{equation}
     \frac{1}{n}\left( \Phi^{n}_{1} (1) \right)^{\top} \left(A^{(1)}_{n,H} \right)\tilde\Phi^{n}_{1}\xrightarrow[n\rightarrow\infty]{\mathbb{P}} k^{(1)}_{H,\phi(1),\phi(2)},\quad \frac{1}{n}\left( \Phi^{n}_{1} (2) \right)^{\top} \left(A^{(1)}_{n,H} \right)\tilde\Phi^{n}_{1}\xrightarrow[n\rightarrow\infty]{\mathbb{P}} \ell^{(1)}_{H,\phi(1),\phi(2)}, \label{eq：A1}
  \end{equation}
  \begin{equation}
       \frac{1}{n}\left( \Phi^{n}_{1} (1) \right)^{\top} \left(A^{(2)}_{n,H} \right)\tilde\Phi^{n}_{1} \xrightarrow[n\rightarrow\infty]{\mathbb{P}} k^{(2)}_{H,\phi(1),\phi(2)},\quad \frac{1}{n}\left( \Phi^{n}_{1} (2) \right)^{\top} \left(A^{(2)}_{n,H} \right)\tilde\Phi^{n}_{1} \xrightarrow[n\rightarrow\infty]{\mathbb{P}} \ell^{(2)}_{H,\phi(1),\phi(2)} \label{eq：A2}
  \end{equation} 
  \begin{equation}
       n^{-\frac{3}{2}}\left( \Phi^{n}_{1} (1) \right)^{\top} \left(A^{(3)}_{n,\overline H_{n}} \right)\tilde\Phi^{n}_{1} = O_{\mathbb{P}}(1),\quad n^{-\frac{3}{2}}\left( \Phi^{n}_{1} (2) \right)^{\top} \left(A^{(3)}_{n,\overline H_{n}} \right)\tilde\Phi^{n}_{1} = O_{\mathbb{P}}(1)  \label{eq：A3}
  \end{equation}
    where $k^{(i)}_{H,\phi(1),\phi(2)}$, $\ell^{(i)}_{H,\phi(1),\phi(2)}$ are constants for $i=1,2$ and  
$A^{(1)}_{n,H}$, $A^{(2)}_{n,H}$, $A^{(3)}_{n,H}$ are
\begin{eqnarray}
A^{(1)}_{n,H} = -\Omega^{-1}_{n,H}\frac{\partial{\Omega_{n,H}}}{\partial H}\Omega^{-1}_{n,H}	 
\end{eqnarray}
\begin{align}
A^{(2)}_{n,H}
&=  \Omega^{-1}_{n,H}\frac{\partial^{2}{\Omega_{n,H}}}{\partial^{2} H}\Omega^{-1}_{n,H}   \\
& \quad + 2\Omega^{-1}_{n,H}\frac{\partial{\Omega_{n,H}}}{\partial H}\Omega^{-1}_{n,H}\frac{\partial{\Omega_{n,H}}}{\partial H}\Omega^{-1}_{n,H} \nonumber
\end{align}
\begin{align}
A^{(3)}_{n,H}
&= 
-\Omega^{-1}_{n,H}\frac{\partial^{3}{\Omega_{n,H}}}{\partial^{3} H}\Omega^{-1}_{n,H} \\
&\quad- 
3\Omega^{-1}_{n,H}\frac{\partial{\Omega_{n }(H)}}{\partial H}\Omega^{-1}_{n,H}\frac{\partial^{2}{\Omega_{n,H}}}{\partial^{2} H}\Omega^{-1}_{n,H}  \nonumber \\
&\quad- 
3\Omega^{-1}_{n - 1}(H)\frac{\partial^{2}{\Omega_{n,H}}}{\partial^{2} H}\Omega^{-1}_{n,H}\frac{\partial{\Omega_{n,H}}}{\partial H}\Omega^{-1}_{n,H} \nonumber \\
&\quad- 
6\Omega^{-1}_{n }(H)\frac{\partial{\Omega_{n,H}}}{\partial H}\Omega^{-1}_{n }(H)\frac{\partial{\Omega_{n,H}}}{\partial H} \nonumber \\
&\quad\quad \times \Omega^{-1}_{n,H}\frac{\partial{\Omega_{n - 1}(H)}}{\partial H}\Omega^{-1}_{n,H}  
\end{align} 
By equations ~\eqref{eq：A1},~\eqref{eq：A2} and ~\eqref{eq：A3}, the consistency of $\hat H_{n}$, Lemma \ref{lemma:nominatortozero} and Lemma \ref{lemma:denominatortopdm}, we can obtain that
\begin{equation}
 \lim_{n\rightarrow \infty}\frac{1}{n}\mathbb{E}\begin{pmatrix}
\left( \Phi^{n}_{1} (1) \right)^{\top} \left(\Omega^{-1}_{n,\hat{H}_{n}} \right)\tilde\Phi^{n}_{1}  \\
\left( \Phi^{n}_{1} (2) \right)^{\top} \left(\Omega^{-1}_{n,\hat{H}_{n}} \right)\tilde\Phi^{n}_{1}  
\end{pmatrix}   = \bf{0} \label{eq:pffenzi}
\end{equation}
and
\begin{equation}
  \lim_{n \rightarrow \infty}\frac{1}{n}\mathbb{E} \begin{pmatrix}
\left( \Phi^{n}_{1} (1) \right)^{\top} \left(\Omega^{-1}_{n,\hat{H}_{n}} \right)\Phi^{n}_{1} (1)& \left( \Phi^{n}_{1} (1) \right)^{\top} \left(\Omega^{-1}_{n,\hat{H}_{n}} \right)\Phi^{n}_{1} (2)\\
\left( \Phi^{n}_{1} (2) \right)^{\top} \left(\Omega^{-1}_{n,\hat{H}_{n}} \right)\Phi^{n}_{1} (1)& \left( \Phi^{n}_{1} (2) \right)^{\top} \left(\Omega^{-1}_{n,\hat{H}_{n}} \right)\Phi^{n}_{1} (2)
\end{pmatrix}^{-1}= Q, \label{eq:pffenmu}
\end{equation}
where $\bf{0}$ is a zero vector and $Q$ is a positive definite matrix.\\
Combining ~\eqref{eq:pffenmu} and ~\eqref{eq:pffenzi}, we can prove that 
 $$\mathbb{E}\begin{pmatrix}
\hat{\phi}_n{(1)}-\phi(1) \\
\hat{\phi}_n{(2)}-\phi(2)
\end{pmatrix}= \textbf{0}.$$
This completes the proof of the consistency of $H$, $\phi(1)$ and $\phi(2)$.\\
\subsection*{B.6.Proof of Theorem \ref{Thm：asynofie}} 
According to Theorem 2 \citep{hurvich1998mean}, without loss of generality, we can assume \( m = o(n^{\frac{2}{3}}) \). Here, \(\mathbb{P}\) denotes convergence in distribution. We thus have
\begin{equation}
\sqrt{m} (\hat{H}_{n} - H) \xrightarrow[n\rightarrow\infty]{\mathbb{P}} \mathcal{N}(0, V_{H}).
\end{equation}
    Building on the results from equation above, we establish that
     \begin{equation}
     \sqrt{m}(\hat \phi_{n}(2) - \phi (2))  = \frac{\tilde \Phi^{n^\top}_{2}(2) \left(\Omega^{-1}_{n-1}( H )\right) \Phi^{n-1}_{1}(1)}{\Phi^{{n-1}^\top}_{1}(1)\left(\Omega^{-1}_{n-1}(\hat{H}_{n})\right)\Phi^{n-1}_{1}(1)}+R^{(1)}_{n}  \label{OSBIAS}
  \end{equation}
  \begin{align}
    \sqrt{m}  \begin{pmatrix}
\hat{\phi}_n{(1)}-\phi(1) \\
\hat{\phi}_n{(2)}-\phi(2)
\end{pmatrix} &  =\begin{pmatrix}
\left( \Phi^{n}_{1} (1) \right)^{\top} \left(\Omega^{-1}_{n,\hat{H}_{n}} \right)\Phi^{n}_{1} (1)& \left( \Phi^{n}_{1} (1) \right)^{\top} \left(\Omega^{-1}_{n,\hat{H}_{n}} \right)\Phi^{n}_{1} (2)\\
\left( \Phi^{n}_{1} (2) \right)^{\top} \left(\Omega^{-1}_{n,\hat{H}_{n}} \right)\Phi^{n}_{1} (1)& \left( \Phi^{n}_{1} (2) \right)^{\top} \left(\Omega^{-1}_{n,\hat{H}_{n}} \right)\Phi^{n}_{1} (2)
\end{pmatrix}^{-1}\begin{pmatrix}
\left( \Phi^{n}_{1} (1) \right)^{\top} \left(\Omega^{-1}_{n,H} \right)\tilde\Phi^{n}_{1}  \\
\left( \Phi^{n}_{1} (2) \right)^{\top} \left(\Omega^{-1}_{n,H} \right)\tilde\Phi^{n}_{1}  
\end{pmatrix} \nonumber \\ 
&\quad + \begin{pmatrix}
R^{(1)}_{n} \\
R^{(2)}_{n} 
\end{pmatrix} 
\label{GLSEOFPHI12}
  \end{align}
        According to the proof of consistency and some results on \citep{esstafa2019long},the denominator of the first term on the right side of the above equation satisfies
      \begin{equation}
      	\frac{1}{n} \begin{pmatrix}
\left( \Phi^{n}_{1} (1) \right)^{\top} \left(\Omega^{-1}_{n,\hat{H}_{n}} \right)\Phi^{n}_{1} (1)& \left( \Phi^{n}_{1} (1) \right)^{\top} \left(\Omega^{-1}_{n,\hat{H}_{n}} \right)\Phi^{n}_{1} (2)\\
\left( \Phi^{n}_{1} (2) \right)^{\top} \left(\Omega^{-1}_{n,\hat{H}_{n}} \right)\Phi^{n}_{1} (1)& \left( \Phi^{n}_{1} (2) \right)^{\top} \left(\Omega^{-1}_{n,\hat{H}_{n}} \right)\Phi^{n}_{1} (2)
\end{pmatrix}^{-1}\xrightarrow[n\rightarrow\infty]{\mathbb{P}} Q,  
      \end{equation}
      The nominator converge to a normal distribution 
      \begin{equation}
      \frac{1}{\sqrt{n}} \begin{pmatrix}
\left( \Phi^{n}_{1} (1) \right)^{\top} \left(\Omega^{-1}_{n,H} \right)\tilde\Phi^{n}_{1}  \\
\left( \Phi^{n}_{1} (2) \right)^{\top} \left(\Omega^{-1}_{n,H} \right)\tilde\Phi^{n}_{1}  
\end{pmatrix}\xrightarrow[n\rightarrow\infty]{\mathcal{L}} \mathcal{N}(0,\Xi_{\phi(1),\phi(2)}),
     \end{equation}
        when $n \rightarrow \infty $,  the reminder $R^{(1)}_{n}$ and $R^{(2)}_{n}$ converge to 0.\\
	   Thus, we can rewrite equation (\ref{OSBIAS}) as follows
\begin{align}
    \sqrt{m}  \begin{pmatrix}
\hat{\phi}_n{(1)}-\phi(1) \\
\hat{\phi}_n{(2)}-\phi(2)
\end{pmatrix} &  =\frac{\sqrt{m}}{n}\begin{pmatrix}
\left( \Phi^{n}_{1} (1) \right)^{\top} \left(\Omega^{-1}_{n,\hat{H}_{n}} \right)\Phi^{n}_{1} (1)& \left( \Phi^{n}_{1} (1) \right)^{\top} \left(\Omega^{-1}_{n,\hat{H}_{n}} \right)\Phi^{n}_{1} (2)\\
\left( \Phi^{n}_{1} (2) \right)^{\top} \left(\Omega^{-1}_{n,\hat{H}_{n}} \right)\Phi^{n}_{1} (1)& \left( \Phi^{n}_{1} (2) \right)^{\top} \left(\Omega^{-1}_{n,\hat{H}_{n}} \right)\Phi^{n}_{1} (2)
\end{pmatrix}^{-1}\nonumber \\
& \times \frac{1}{n}\begin{pmatrix}
\left( \Phi^{n}_{1} (1) \right)^{\top} \left(\Omega^{-1}_{n,H} \right)\tilde\Phi^{n}_{1} \\ 
\left( \Phi^{n}_{1} (2) \right)^{\top} \left(\Omega^{-1}_{n,H} \right)\tilde\Phi^{n}_{1}  
\end{pmatrix}  + \begin{pmatrix}
R^{(1)}_{n} \\
R^{(2)}_{n} 
\end{pmatrix} \nonumber \\
 &  =\frac{\sqrt{m}}{n}Q^{-1} \begin{pmatrix}
\left( \Phi^{n}_{1} (1) \right)^{\top} \left(\Omega^{-1}_{n,H} \right)\tilde\Phi^{n}_{1} \\ 
\left( \Phi^{n}_{1} (2) \right)^{\top} \left(\Omega^{-1}_{n,H} \right)\tilde\Phi^{n}_{1}  
\end{pmatrix}  + \begin{pmatrix}
R^{(1)}_{n} \\
R^{(2)}_{n} 
\end{pmatrix} 
\label{GLSEOFPHI12}
  \end{align}
      By slutsky theorem, we can conclude that
      $\sqrt{m}  \begin{pmatrix}
\hat{\phi}_n{(1)}-\phi(1) \\
\hat{\phi}_n{(2)}-\phi(2)
\end{pmatrix} $ converges to a  joint normal distribution.  
\par Lastly, we consider expressing these results in the form of a joint normal distribution of all parameters. Based on the findings of \citep{hariz2024fast}, the asymptotic distribution of $\sqrt{n}({\hat{\phi}}_{n}(u)-\phi(u))$ can be represented as a constant multiple of the asymptotic distribution of $(\hat H_{n}-H)$. Furthermore, by the Cramer-Wold theorem, the asymptotic distribution of $\sum^{T}_{u=1}({\hat{\phi}}_{n}(u)-\phi(u))$  remains asymptotically normal. Thus, the vector
	  $$({\hat{\phi}}_{n}(1)-\phi(1),{\hat{\phi}}_{n}(2)-\phi(2),\hat H_{n}-H),$$ 
     converges to a Gaussian vector, tending towards a joint normal distribution. The covariance matrix of this vector is
      \begin{equation}
	   \tilde \Xi_{\phi(1),\phi(2)} = U_{\phi(1),\phi(2)}U^\top_{\phi(1),\phi(2)},
      \end{equation}
	  where
  \begin{equation}	
      U_{\phi(1),\phi(2)} = \left(
	 \begin{array}{c}
		1\\
		 Q^{(1)}_{\phi(1),\phi(2)} \\
		 
		 Q^{(2)}_{\phi(1),\phi(2)}
	\end{array}
	\right),
 \end{equation}
        $Q^{(1)}_{\phi(1),\phi(2)}$ and $Q^{(1)}_{\phi(1),\phi(2)}$ are the constants related to $\phi(1)$, $\phi(2)$.
\section*{Appendix 3. Proof of the Main Results for One-Step Estimatior}
\subsection*{C.1.Proof of Lemma \ref{Lemma：regularityofsdofz}} 
     We start from Assertion 3a, which states that
     \begin{equation}
       p_{H,\phi(2),\phi(1)}(\lambda) = \left|1+\frac{e^{-2\lambda i}}{\phi(1)[1-e^{-2\lambda i}\phi^{2}(2)\phi^{2}(1)]}+\frac{e^{-\lambda i}(1+\phi(1))}{1-e^{-2\lambda i}\phi^{2}(2)\phi^{2}(1)}\right|^{2}f_{\epsilon^{H}_{n}},   
    \end{equation}
     where $C_H=\frac{1}{2\pi}\Gamma(2H+1)sin(\pi H)$ and $\Gamma(.)$ denote the Gamma function. 
     According to lemma 5.4 in \citep{hariz2024fast}, we have
     \begin{equation}
     K_{1,\delta}|\lambda|^{1-2H+\delta } \leq C_{H}(1-cos(\lambda))\sum\limits_{j \in Z}\frac{1}{|\lambda+2j\pi|^{2H+1}} \leq K_{2,\delta}|\lambda|^{1-2H-\delta }, 
     \end{equation}
     and 
     \begin{equation}
       H^{(1)}_{\phi(1),\phi(2)} \leq \left|1+\frac{e^{-2\lambda i}}{\phi(1)[1-e^{-2\lambda i}\phi^{2}(2)\phi^{2}(1)]}+\frac{e^{-\lambda i}(1+\phi(1))}{1-e^{-2\lambda i}\phi^{2}(2)\phi^{2}(1)}\right|^{2}\leq H^{(2)}_{\phi(1),\phi(2)},
     \end{equation}
     where $H^{(1)}_{\phi(1),\phi(2)}= \frac{\left|\left|\phi(1)\right|-\left|\phi(1)+\phi^{2}(1)\right|-\left|1 - \phi^{3}(1)\phi^{2}(2)\right|\right|^{2}}{\left|\phi(1)\right|^{2}\left(1+\left|1 - \phi^{3}(1)\phi^{2}(2)\right|\right)^{2}}$ and $H^{(2)}_{\phi(1),\phi(2)}=\frac{\left(\left|\phi(1)\right|+\left|\phi(1)+\phi^{2}(1)\right|+\left|1 - \phi^{3}(1)\phi^{2}(2)\right|\right)^{2}}{\left|\phi(2)\right|^{2}\left|\ 1-\left|1 - \phi^{3}(1)\phi^{2}(2)\right|\right|^{2}}$. Thus, Assertion 3a has been proved and assertion 3b follows straightforwardly from Assertion 3a.\\
Next, we discuss Assertion 3c, which can be obtained directly from Lemma 5.4. The partial derivative of  $p_{H,\phi(1),\phi(2)}(\lambda)$ does not depend on $\left|1+\frac{e^{-2\lambda i}}{\phi(1)[1-e^{-2\lambda i}\phi^{2}(2)\phi^{2}(1)]}+\frac{e^{-\lambda i}(1+\phi(1))}{1-e^{-2\lambda i}\phi^{2}(2)\phi^{2}(1)}\right|^{2}$, and the modulus is bounded.
\subsection*{C.2.Proof of Lemma \ref{Lemma：regular1}} 
       Without loss of generality, let $B(\boldsymbol\theta_{0},\delta)$ be a convex set in $\mathbb{R}^{3}$. For ease of notation, $p_{H,\phi(1),\phi(2)}(\lambda)$ can be denoted as $p_{\boldsymbol\theta}(\lambda)$. According to the relevant conclusions in \citep{cohen:hal-00638121} and the discussion of regularity conditions for $p_{H,\phi(1),\phi(2)}(\lambda)$, it is known that for any $k,j\in \left\{1,2,...,d\right\}$ that the following inequality holds
        \begin{equation}
           |\frac{1}{4\pi}(\int_{-\pi}^{\pi}\frac{\partial log \,   p_{\boldsymbol\theta}(\lambda)}{\partial \boldsymbol\theta_{k}}\frac{\partial log \,   p_{\boldsymbol\theta}(\lambda)}{\partial \boldsymbol\theta_{j}}d\lambda) -\frac{1}{4\pi}(\int_{-\pi}^{\pi}\frac{\partial log \,   p_{\boldsymbol\theta_{0}}(\lambda)}{\partial \boldsymbol\theta_{0,k}}\frac{\partial log \,   p_{\boldsymbol\theta_{0}}(\lambda)}{\partial \boldsymbol\theta_{0,j}}d\lambda)|\leq
           K||\boldsymbol\theta-\boldsymbol\theta_{0}||, 
       \end{equation}
       $K$ is defined as 
      \begin{equation}
          K = \sup_{\boldsymbol\theta \in B(\boldsymbol\theta_{0},\delta)}\left\lVert(\frac{\partial}{\partial  \theta_{i}}(\int^{\pi}_{-\pi}\frac{\partial log \,   p_{\boldsymbol\theta}(\lambda)}{\partial  \theta_{k}}\frac{\partial log \,   p_{\boldsymbol\theta}(\lambda)}{\partial  \theta_{j}}d\lambda))_{1 \leq i \leq d}\right\rVert, 
      \end{equation}
         which is related to $k$ and $j$.
      Furthermore, since the conditions (A1) and (A2) of \citep{cohen:hal-00638121} hold, it follows that $K< \infty$, hence the lemma holds.
\subsection*{C.3.Proof of Lemma \ref{Lemma：regular2}} 
The Lemma 3.6 \citep{cohen:hal-00638121} implies that, from the distribution of $\boldsymbol\theta$, we have
 \begin{equation}
   E\left(\frac{\Delta l_{n}(\boldsymbol\theta)}{n}\right)\rightarrow -\mathcal{I}(\boldsymbol\theta).   
  \end{equation}
To determine the convergence rate of the above expression, Lemma 3 and Lemma 4 \citep{lieberman2012asymptotic} yield the following conclusion
  \begin{equation}
     E\left(\frac{\Delta l_{n}(\boldsymbol\theta)}{n}\right)+\mathcal{I}(\boldsymbol\theta) = O(n^{-1+\delta}),
  \end{equation}
where $\delta$ is a positive real number. Therefore,
  \begin{equation}
    E\left(\frac{\Delta l_{n}(\boldsymbol\theta)}{\sqrt{n}}\right)+\sqrt{n}\mathcal{I}(\boldsymbol\theta) = O(n^{-\frac{1}{2}+\delta}).
 \end{equation}
Furthermore, by utilizing Lemma 3.6 \citep{cohen:hal-00638121} once again, we obtain
\begin{equation}
    Var\left(\frac{\Delta l_{n}(\boldsymbol\theta)}{\sqrt{n}}\right)=O(1).
\end{equation}
   Thus, the proof is concluded.
\subsection*{C.4.Proof of Lemma \ref{Lemma：regular3}}
      Let $C_{k,\boldsymbol\theta}$ be a compact convex set depending on $k$ and $\boldsymbol\theta$, and ${\overline{\boldsymbol\theta}}_{n} \in C_{k,\boldsymbol\theta}$. According to the proof of Lemma 3.7 in \citep{cohen:hal-00638121}, we have
 \begin{equation}
  \sup_{{\overline{\boldsymbol\theta}}_{n} \in C_{k,\boldsymbol\theta}}\left|\frac{\partial^{3}}{\partial^{i_{1}}{\theta_{1}}\partial^{i_{2}}{\theta_{2}}...\partial^{i_{d}}{\theta_{d}}}\frac{l_{n}({{\overline{\boldsymbol\theta}}_{n}})}{n^{1+k}}\right| = O_{\mathbb{P}}(1),
\end{equation}
where $(i_{1},i_{2},...,i_{d}) \in \left\{0,1,2,3\right\}^{d}$, satisfying $i_{1}+i_{2}...+i_{d} =3$. 
In conclusion, for a finite positive random variable $K$, we have
\begin{equation}
     P\left(\left\lVert \frac{\Delta l_{n}({\overline{\boldsymbol\theta}}_{n})}{n} - \frac{\Delta l_{n}(\boldsymbol\theta)}{n}\right\rVert \leq Kn^{k}(||{\overline{\boldsymbol\theta}}_{n}-\boldsymbol\theta||\right) \geq P({\overline{\boldsymbol\theta}}_{n} \in C_{k,\boldsymbol\theta}),
\end{equation}
which implies $\frac{\Delta l_{n}({\overline{\theta}}_{n})}{n} - \frac{\Delta l_{n}(\theta)}{n} = O_{\mathbb{P}}(n^{k}({\overline{\theta}}_{n}-\theta))$ holds.
\subsection*{C.5.Proof of Lemma \ref{Thm:asnofos}}
 We will discuss  the consistency and asymptotic normality of One-Step estimator.\\
	 \textbf{(1) $\sqrt{n}$-Consistency of $\tilde{\theta}_{n}$} 
  \par Observing equation (\ref{OSESTIMATOR}) ,  the first and second terms on the right-hand side can be expressed as follows 
     \begin{eqnarray}
         A_{n} 
         &=& \sqrt{n}(\hat {\boldsymbol\theta}_{n}-\boldsymbol\theta)\mathcal{I}^{-1}(\hat {\boldsymbol\theta}_{n})(\mathcal{I}(\hat {\boldsymbol\theta}_{n})+\frac{\int^{1}_{0}\bigtriangleup l_{n}(\boldsymbol\theta(v))dv}{n}),    \\ 
         &=& \sqrt{n^{\delta}}(\hat {\boldsymbol\theta}_{n}-\boldsymbol\theta)\mathcal{I}^{-1}(\hat {\boldsymbol\theta}_{n})\sqrt{n^{1-\delta}}(\mathcal{I}(\hat {\boldsymbol\theta}_{n})+\frac{\int^{1}_{0}\bigtriangleup l_{n}(\boldsymbol\theta(v))dv}{n}),  \nonumber
     \end{eqnarray}
       and
    \begin{equation}
      B_{n}=\mathcal{I}^{-1}(\hat {\boldsymbol\theta}_{n})\frac{\nabla l_{n}(\boldsymbol\theta)}{\sqrt{n}}=\mathcal{I}^{-1}(\boldsymbol\theta)\frac{\nabla l_{n}(\boldsymbol\theta)}{\sqrt{n}}+(\mathcal{I}^{-1}(\hat {\boldsymbol\theta}_{n})-\mathcal{I}^{-1}(\boldsymbol\theta))\frac{\nabla l_{n}(\boldsymbol\theta)}{\sqrt{n}}.
    \end{equation}
      Fristly, we analyze the properties of $A_{n}$ and derive the following equation
   \begin{equation} \label{t5}
     \begin{aligned}
    \mathcal{I}(\hat {\boldsymbol\theta}_{n})+\frac{\int^{1}_{0}\Delta l_{n}(\boldsymbol\theta(v))dv}{n}
    &= (\mathcal{I}(\hat {\boldsymbol\theta}_{n})-\mathcal{I}(\boldsymbol\theta_{n})) \\
    &\quad +(\mathcal{I}(\boldsymbol\theta_{n})+\frac{\Delta l_{n}(\boldsymbol\theta)}{n})+\frac{1}{\sqrt{n}}\int^{1}_{0}(\frac{\Delta l_{n}(\boldsymbol\theta(v)))}{\sqrt{n}}-\frac{\Delta l_{n}(\boldsymbol\theta)}{\sqrt{n}})dv,
   \end{aligned}
   \end{equation}
     Based on equation ($\ref{t5}$) and lemmas 5.2, 5.3, and 5.4. The convergence order of $\frac{A_{n}}{\sqrt{n}}$ is
   \begin{equation}   
      \frac{A_{n}}{\sqrt{n}} = n^{-\frac{\delta}{2}}(O_{\mathbb{P}}(n^{-\frac{\delta}{2}})+O_{\mathbb{P}}(n^{-\frac{1}{2}})+O_{\mathbb{P}}(n^{k-\frac{\delta}{2}})),
 \end{equation}
        when $k-\delta < 0$, we have $\frac{A_{n}}{\sqrt{n}}\xrightarrow[n\rightarrow\infty]{\mathbb{P}}0$.\\
Secondly, we consider the property of $B_{n}$ and it has the form of
   \begin{equation}
   B_{n}=\mathcal{I}^{-1}(\hat {\boldsymbol\theta}_{n})\frac{\nabla l_{n}(\boldsymbol\theta)}{\sqrt{n}}=\mathcal{I}^{-1}(\boldsymbol\theta)\frac{\nabla l_{n}(\boldsymbol\theta)}{\sqrt{n}}+(\mathcal{I}^{-1}(\hat {\boldsymbol\theta}_{n})-\mathcal{I}^{-1}(\boldsymbol\theta))\frac{\nabla l_{n}(\boldsymbol\theta)}{\sqrt{n}}, 
   \end{equation}
   According to \citep{hariz2024fast} and theorem 1 in \citep{lieberman2012asymptotic}, we have
   \begin{equation}
    \frac{\nabla  l_{n}({\boldsymbol\theta})}{\sqrt{n}} \xrightarrow[n\rightarrow\infty]{\mathbb{P}} 0.
   \end{equation}
   When $\mathcal{I}n(\cdot)$ is a non-degenerate continuous function, as indicated by the above equation, it can be observed that both the first and second terms of $B_{n}$ tend to 0. Consequently, $\frac{B_{n}}{\sqrt{n}}$ converges in probability to 0, and naturally, it also converges in distribution to 0.
   \par Combining the above results, we can conclude the consistency of $\tilde {\boldsymbol\theta}_{n}$. \\
\textbf{(2) Asymptotic normality of $\tilde{\boldsymbol\theta}_{n}$ }
    \par According to the results of \citep{hariz2024fast}, the equation
	\begin{equation}
        \mathcal{I}^{-1}(\boldsymbol\theta)\frac{\nabla  l_{n}(\boldsymbol\theta)}{\sqrt{n}}+(\mathcal{I}^{-1}(\hat {\boldsymbol\theta}_{n})-\mathcal{I}^{-1}(\boldsymbol\theta))\frac{\nabla  l_{n}(\boldsymbol\theta)}{\sqrt{n}}, \label{6.60}
    \end{equation}
    converges in probability to a bounded limit as $n \rightarrow \infty$. Simultaneously, the second term on the right-hand side of equation ($\ref{6.60}$) converges to 0. By applying the Slutsky theorem, we can verify the asymptotic normality of $\tilde \theta_{n}$.
    \newpage
    \section*{Appendix 4. Simulation for the Initial Estimator in Remark 3}
  \begin{figure}[htbp]
    \centering
     \vspace{2\baselineskip}
   \includegraphics[width=1\textwidth]{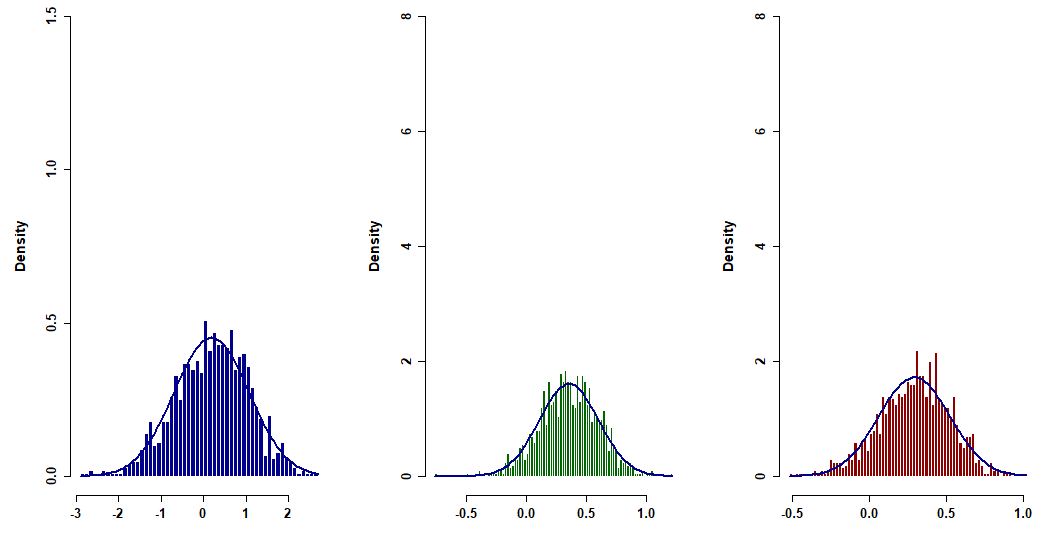}
 \hspace{30pt}\bf{(a) IE of \boldmath$H$ = 0.6}
\hspace{70pt} 
 \text{(b) IE of} $\boldsymbol{\phi}$(1) = 0.7
\hspace{60pt}
\text{(c) IE of} \boldmath$\mathbf{\phi}$(2) = 0.6\\
   \caption{{The Simulation of Initial Estimator in Remark 3 where $\theta=(0.7,0.6,0.6)$ for $ n=1000$.}}
    \label{A4-PIC-766}
\end{figure}   
\newpage
\begin{figure}[htbp]
    \centering
     \vspace{2\baselineskip}
   \includegraphics[width=1\textwidth]{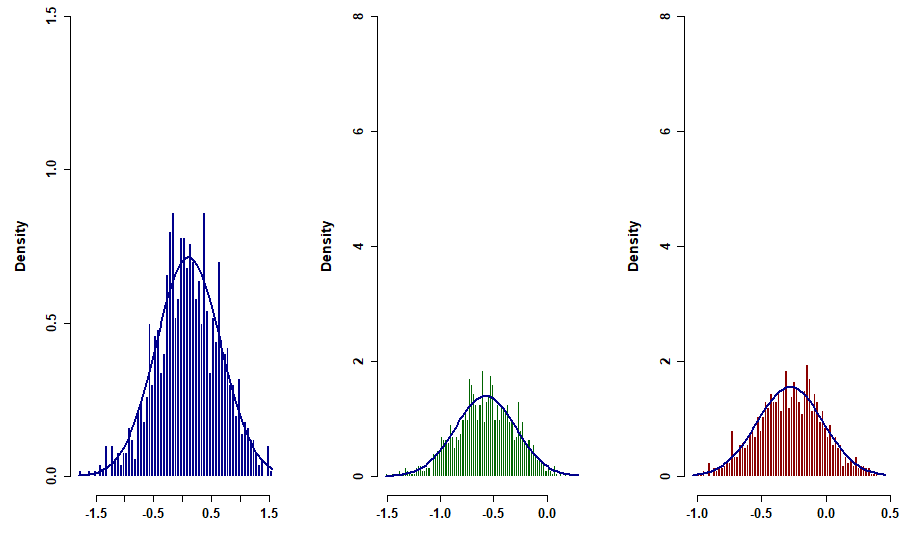}
 \hspace{30pt}\bf{(a) IE of \boldmath$H$ = 0.4}
\hspace{70pt} 
 \text{(b) IE of} $\boldsymbol{\phi}$(1) = 0.9
\hspace{60pt}
\text{(c) IE of} \boldmath$\mathbf{\phi}$(2) = 0.6\\
   \caption{{The Simulation of Initial Estimator in Remark 3 where $\theta=(0.9,0.6,0.4)$ for $ n=1000$.}}
    \label{A4-PIC-964}
\end{figure} 
\newpage
\begin{figure}[htbp]
    \centering
     \vspace{2\baselineskip}
   \includegraphics[width=1\textwidth]{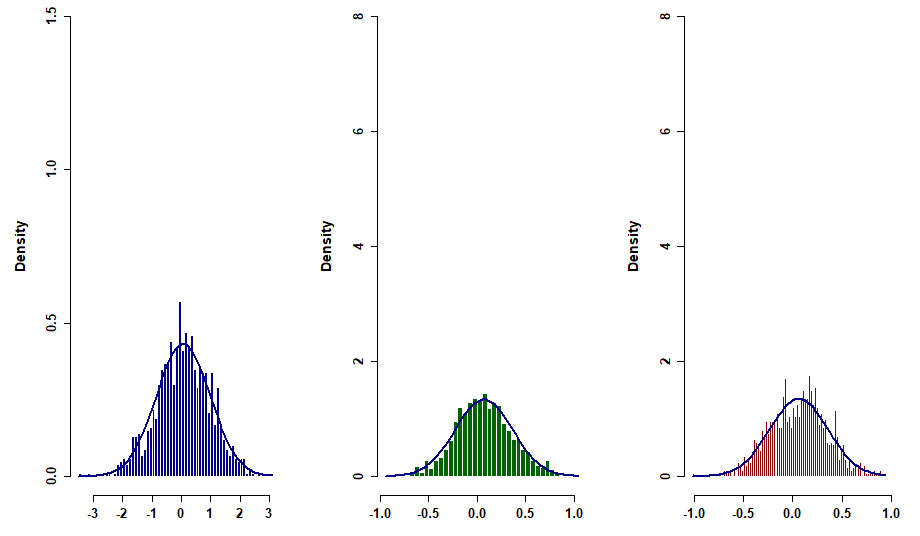}
 \hspace{30pt}\bf{(a) IE of \boldmath$H$ = 0.5}
\hspace{70pt} 
 \text{(b) IE of} $\boldsymbol{\phi}$(1) = 0.5
\hspace{60pt}
\text{(c) IE of} \boldmath$\mathbf{\phi}$(2) = 0.5\\
   \caption{{The Simulation of Initial Estimator in Remark 3 where $\theta=(0.5,0.5,0.5)$ for $ n=0000$.}}
    \label{A4-PIC-555}
\end{figure} 
\newpage
\begin{table}[htbp]
    \centering
    \footnotesize 
    \caption{{Bias and RMSE of the estimator in Remark 3}}
    \begin{tabular}{*{7}{c}}
        \toprule
        &  B (0.7,0.6,0.6)  & RMSE    (0.7,0.6,0.6)  &  B    (0.9,0.6,0.4)  & RMSE   (0.9,0.6,0.4) &  B    (0.5,0.5,0.5)  & RMSE (0.5,0.5,0.5)   \\
        \midrule
        $H$   &  0.0241 & 0.1111&  0.0130&0.0703  &0.0092  &0.1166\\
        $\phi(1)$& 0.0453 & 0.0313&-0.0727  & 0.0359 & 0.0085 &0.0375 \\
        $\phi(2)$  & 0.0368 & 0.0292& -0.0351 & 0.0323 & 0.0077 &0.0373 \\
        \bottomrule
    \end{tabular}
    \label{table5}
\end{table}
From the above figures, we find that this estimation is unbiased only when \(H=\frac{1}{2}\), it has left bias when \(0<H<\frac{1}{2}\) and it has right bias when \(\frac{1}{2}<H<1\). In our simulations, we found that repeatedly applying the one-step estimator (i.e., using the multi-step estimator) to this initial estimator leads to a reduction in bias. However, determining the optimal number of steps for the multi-step estimator remains an open problem.

\end{document}